\pdfoutput=1
\documentclass{article}
\usepackage{graphicx} 
\usepackage{amsmath}
\usepackage{amssymb}
\usepackage[font=large]{caption}

\newcommand{\N}{\mathbb{N}}
\newcommand{\Z}{\mathbb{Z}}
\newcommand{\R}{\mathbb {R}}
\newcommand{\C}{{\bf C}}
\newcommand{\Pow}{{\bf P}}

\newcommand{\s}{\subseteq}
\newcommand{\ra}{\rightarrow}
\newcommand{\es}{\emptyset}
\newcommand{\Ra}{\Rightarrow}
\newcommand{\LRa}{\Leftrightarrow}
\newcommand{\zero}{{\bf 0}}
\newcommand{\one}{{\bf 1}}

\title{Advertising  finite commutative semigroups}

\author{Marcel Wild (mwild@sun.ac.za)}

\begin{document}

\maketitle

\begin{center}
Stellenbosch University, Dept. Mathematics, South Africa
\end{center}

\begin{abstract}
Every mathematician is familiar with the beautiful structure of finite commutative groups. What is less well known is that finite commutative semigroups also have a neat and well-described structure.  We prove this in an efficient fashion. We unravel the  structural details of many concrete finite commutative semigroups.  Here ``concrete” comes in two types. First, we examine the structure of the multiplicative semigroups 
$(\Z_n,\odot)$ (more interesting than their bland siblings $(\Z_n,+)$) and show that it depends on the prime factors of $n$ in interesting ways. Second, we thoroughly treat finite commutative  semigroups defined by generators and relations. Our aim is to provide a comprehensive introduction to the area, but with some enticing directions for the expert to follow.
\end{abstract}


\vspace{3mm}

\section{Introduction}

Every teenager has  a basic understanding of how multiplication of integers $a,b,c \in\Z$ behaves, in particular  $ab=ba$ and $(ab)c=a(bc)$. Imagine shrinking $\Z$ to a finite set while keeping the properties of multiplication. If at all possible, what you get must be a so called  commutative finite  semigroup (in fact $\Z_n$, more on which in a moment).

This article has a double purpose. On the one hand, it attempts to advertise the beauty of commutative semigroups to "type 1" mathematicians that may never have gotten beyond the definition of  "semigroup", but who enjoy algebra and axioms. For this audience I inserted (inspired by Allenby [A]) little comments like "check" or "why?" throughout the text. The type 2 mathematician knows (most of) the  material but may find some novel\footnote{Novel {\it results} will be offered  as well, some proofs of which   being published elsewhere.} points of view.

In the remainder of the introduction we pin down (in Subsection 1.1)  what even the type 1 reader should master before reading on. Subsection 1.2 sketches  how the structure of commutative finite semigroups  relies on three ingredients:  semilattices, nil semigroups and groups. All three ingredients come to the fore already in the {\it multiplicative} semigroup $(\Z_n,\odot)$ which will receive special attention. Subsection 1.3 provides the detailed Section break up, and (for type 2 readers) an outline of three Open Questions (formulated carefully later on) that are hoped to stimulate research.

\vspace{5mm}
{\bf 1.1} A binary\footnote{We sometimes use dot notation $a\cdot b$ or simply concatenation $ab$.} operation $\ast$ on a set $S$ is {\it associative} if $(a\ast b)\ast c=a\ast (b\ast c)$ for all $a,b,c\in S$. One then calls $(S,\ast)$ a {\it semigroup (sgr)}. If $H\s S$ is a nonempty subset such that $a\ast b\in H$ for all $a,b\in H$, then $(H,\ast)$ is also a semigroup (why?), a so-called {\it subsemigroup} of $S$.
We will mainly focus on  {\it commutative (c.)} semigroups $S$, i.e. $a\ast b=b\ast a$ for all $a,b\in S$. Usually they are finite (f.) as well. However, if a concept can be smoothly defined for arbitrary semigroups, there is no need to impose finiteness or commutativity. A map $f:S\to T$ between semigroups is called a {\it morphism} if\footnote{When  composition of functions occurs, it is often handy to write the function symbols on the {\it right} because then $x(f\circ g)=(xf)g$, and so $f\circ g$ is the function one gets by first applying $f$ and then $g$. This is in line with the natural left-to-right direction of reading (as opposed to conventional notation where $(f\circ g)(x)=f(g(x))$).} $(xy)f=(xf)(yf)$ for all $x,y\in S$. One verifies that the {\it image} $\{xf:\ x\in S\}$ is a subsemigroup of $T$. A surjective morphism is an {\it epimorphism} and a bijective one an {\it isomorhism}. We write $S\simeq T$ if there is an isomorphism between $S$ and $T$. For semigroups $S_1,\ldots,S_t$ the direct product $S_1\times\cdots\times S_t$ becomes itself a semigroup under component-wise multiplication, i.e. $$(x_1,...,x_t)\ast (y_1,...,y_t):=(x_1y_1,...,x_ty_t)$$

We will use the shorthand "iff" for "if and only if". The quotient ring $\Z/n\Z$ we write as $\Z_n=\{0,1,...,n-1\}$.  If statements like "$5\odot 7=12$ in $\Z_{23}$" and "$51\odot 30\equiv 81\ (mod\ 23)$" perplex you, please\footnote{On two occasions we ever so briefly deviate from $(\Z_n,\odot)$ to the {\it ring structure} $(\Z_n,+,\odot)$. Nothing more than the equivalence of injectivity with the triviality of the kernel will be used.} consult  [A,Sec.2.4] before reading on. 

\vspace{5mm}

{\bf 1.2} An element $\zero$ of a semigroup (sgr) $S$ is a {\it zero} if $\zero x=x\zero=\zero$ for all $x\in S$. In particular $\zero\zero=\zero$. An element $\one$ of a sgr $S$ is an {\it identity} if $\one x=x\one=x$ for all $x\in S$. In particular $\one\one=\one$.
More generally one calls $e\in S$  an {\it idempotent} if $ee=e$. One can show that each f. sgr contains at least one idempotent.

Hence there are two extreme cases of c.f. semigroups. Those with {\it all} their elements being idempotents (so called {\bf semilattices} $Y$), and those with exactly {\it one} idempotent $e$ (so called {\bf Archimedean} sgr $A$). The latter are the topic of Section 7 and they include two natural special cases. Either $e$ is an identity (in which case $A$ is a commutative "group" - a structure likely familiar also to type 1 readers), or $e$ is a zero (in which case $A$ is a "nil" semigroup, where by definition each element has some power which is $\zero$). If $e$ is neither $\one$ nor $\zero$, then $A$ is nevertheless an elegant kind of conglomerate of a group and a nil semigroup. 

\vspace{2mm}

The structure theorem for c.f. semigroups $S$ states that $S$ is a disjoint union of Archimedean subsemigroups $A_i$ (with unique idempotent $e_i$). This yields a "local" understanding of $S$, but what happens "globally", e.g. where is $xy$ located when $x\in A_i$ and $y\in A_j\ (i\neq j)$?

For starters, due to commutativiy $e_ie_j$ is idempotent as well (why?), say $e_ie_j=e_k$ for some index $k$. It turns out that $xy$ sits in $A_k$.
This also shows that the set $Y$ of all idempotents $e_i$ is a ssgr of $S$ which hence (on its own) is a semilattice. 

\vspace{3mm}
{\bf 1.3} After discussing  1-generated semigroups and morphisms between them, we turn to nil semigroups and ideals (Sec.2), then to monoids and groups (Sec.3). While the proof of the Fundamental Theorem of finite Abelian groups is omitted, another nontrivial (and rarely proven) fact will  be given  full attention in Subsection 3.7.2. Next come closure systems (Sec.4) and semilattices (Sec.5).
The reader's possible impression that too much attention is devoted to them will hopefully be revised with hindsight.

Section 6 handles {\it relatively free commutative semigroups} $RFCS(...)$ in a painless way that (initially) avoids congruence relations. Instead the "local confluence" of a semigroup presentation will take center stage. In a nutshell, local confluence guarantees that certain  "normal forms" bijectively match the elements of $RFCS(...)$. At this stage the multiplication table (aka {\it Cayley table}) of $RFCS(..)$ could be set up, but not the "fine" structure of $RFCS(..)$. 

Having had a closer look at Archimedean semigroups in Section 7, the Structure Theorem (glimpsed in 1.2) gets proven in Section 8.
Afterwards an original five step recipe is presented to unravel the fine structure of each c.f. semigroup whose Cayley table is known. The recipe is carried out on two types of semigroups that received preliminary attention in Sections 6 and 7, i.e.  $RFCS(..)$ and $\Z_n$.

Section 9 is devoted to the Ideal Extension Problem, with emphasis on the case where the two involved semigroups are finite and cyclic.

In Section 10 (titled "Loose ends") we give more background on $RFCS(..)$. This includes congruence relations and the Church-Rosser property of digraphs. We also glimpse at {\it arbitrary} semigroups and how their structure is assessed in terms of the famous Green equivalence relations. We then point out how much of this collapses in the finite and commutative case.

\vspace{2mm}
Quoting from page 2 of [G]: {\it By well-established tradition, we regard as solved any problem which can be stated in terms of groups or semilattices (we dump it onto other unsuspecting mathematicians).}

\vspace{2mm}
I dare to break with this tradition: Neither semilattices nor Abelian groups will complain of having been neglected. As to the former, a more efficient approach to calculating semilattices (defined by generators and relations) is offered. As to the latter, the well-known\footnote{But apparently not in the MathOverflow internet community.} fact that the orders of the elements of a commutative group determine its isomorphism type, is proven in detail.

\vspace{3mm}
{\bf 1.3.1} Here, in brief, the content of the mentioned three Open Questions. Question 1 (in Sec.7) assumes that the c.f. semigroup $S$ is a direct product of cyclic semigroups, and asks in how many ways this is possible. For the special case where $S$ is a group, this is both well known and nontrivial already. 

Question 2 (also in Sec.7) considers an arbitrary finite commutative ring $(R,+,\odot)$ and asks how much is known about the structure of the semigroup $(R,\odot)$. The question is posed after we have unraveled in detail\footnote{The author is not aware that this has been done in similarly reader-friendly ways beforehand, but welcomes to be taught otherwise.} what happens for the particular case where $(R,+,\odot)$ is of type $(\Z_n,+,\odot)$.

As to Question 3, we frequently use that $\le_{\cal J}$ is a partial order for various types of semigroups (most notably nil sgr and semilattices), but only in Section 10 (=Loose ends) we give the definition of the ${\cal J}$-relation itself and, among other things, ask when $\cal J$ is a retract congruence.

\vspace{2mm}
{\bf 1.3.2} Although the paper in front of you is mainly a survey article, there are bits of original research (to the author's best knowledge), most prominently Section 9. But also Subsections 6.7, 6.8 about semilattices, and the five step recipee in 8.3. The  people most responsible for making me a semigroup aficionado, in alphabetic order, are  P. Grillet, J. Howie, J.E. Pin. More detailed credentials and a larger list of references may be given in a later version.

Readers are invited to contribute to the Open Questions, or to anything else in the realm of finite commutative semigroups. If ever the so developing arXiv-version  reaches a certain volume and maturity, one may undertake transforming it into a book(let).

\section{Cyclic semigroups, nilsemigroups, ideals}

\vspace{5mm}
 Finite cyclic semigroups $\langle a\rangle$ have  a  "body" which is a cyclic group, but additionally they may have a "tail". Hence  f. cyclic semigroups  are more complicated than f. cyclic groups (2.2). This is also reflected in Theorem 1 which characterizes the morphisms between two f. cyclic semigroups (2.3). In 2.4 to 2.7 we introduce c.f. nil semigroups and show that they are partially ordered by a natural binary relation $\le_{\cal J}$. Subsection 2.8 introduces the free commutative sgr $F_k$ and its military order. In 2.9
 we are concerned with ideals in c. semigroups $S$. In particular, when $S$ is finite, it has a  "kernel" (=smallest ideal).
 
\vspace{3mm}
{\bf 2.1} It follows from associativity [BC,p.39] that for all elements $a_1,a_2,\cdots , a_n$ in any semigroup $(S,\ast)$ the product $a_1\ast a_2\ast\cdots\ast a_n$ is well-defined, i.e. independent of the way it is bracketed.
In particular the definition $a^n:= a\ast a\ast\cdots\ast a$ ($n\ge 1$ factors $a$) is well-defined for all $a\in S$. One verifies by induction (try) that 

$$(1)\quad a^i\ast a^j=a^{i+j}\ and\ (a^i)^j=a^{ij}\ for\ all\ i,j\ge 1.$$

\noindent
If $(S,\ast)$ is commutative and\footnote{The notation $a,b,..,c$ (which is adopted from Gian-Carlo Rota's lectures) beats bothering with the subscripts of $a_1,a_2,\ldots ,a_n$. Particularly when powers of these elements are considered.}
$a,b,.., c\in S$, then there is a smallest\footnote{Thus  $\langle a,b,..,c\rangle$  is contained in every ssgr of $S$ that contains $a,b,...,c$.} ssgr $\langle a,b,..,c\rangle$ of $S$  containing these elements. In fact

$$(2)\quad \langle a,b,..,c\rangle = \{a^i\ast b^j\ast\cdots\ast c^k:\ (i,j,..,k)>(0,0,..,0)\}.$$

\noindent
One calls it the subsemigroup {\it generated } by $a,b,.., c$. Notice that\footnote{As for any partial order, also for the componentwise order $\le$ on $\N^m$ we write $x<y$ if $x\le y$ but $x\neq y$. Hence the expression $(i,j,..,k)>(0,0,..,0)$ in (2) means that {\it not all} of $i,j,..,k$ are zero.}  e.g. $a^2b^0c^7:=a^2c^7$.

\vspace{6mm}
{\bf 2.2}
Of particular importance will be the  {\it cyclic} semigroup\\ $\langle 
a\rangle=\{a^i:\ i>0\}=\{a^i:\ i\ge 1\}$. 

If additionally $|\langle a\rangle|<\infty$, then there must be $m,n\ge 1$ such that $a^m=a^{m+n}$. The smallest such $m$  will henceforth again be denoted by $m$, and the smallest $n$ (for the obtained $m$) will again be $n$. We will write $C_{m,n}$  for a cyclic sgr $\langle a\rangle$ of this type\footnote{Whether conversely for each pair $(m,n)\ge (1,1)$ there {\it exists} such a sgr $C_{m,n}$, is a puzzling question which will be answered in Sec. 6.}. Hence

$$(3)\quad C_{m,n}=\{a,a^2,\ldots,a^m,\ldots,a^{m+n-1}\}\ and\ |C_{m,n}|=m+n-1,$$

\noindent
and for all $j,k\ge 0$ it holds that

$$(4)\quad a^{m+j}=a^{m+k}\ {\it iff}\ m+j\equiv m+k\ (mod\ n).$$

\noindent
One calls $m$ the {\it index} and $n$ the {\it period} of  $C_{m,n}=\langle a\rangle$.
Furthermore, $\{a,...,a^{m-1}\}$ and $\{a^m,...,a^{m+n-1}\}$ are the {\it tail} and {\it body} of $\langle a\rangle$, respectively. The tail can be empty (if $m=1$) but  the minimum cardinality of the body $H$ is $1$ (if $n=1$). Clearly $H$ is a ssgr of $\langle a\rangle$. In 3.3 we find out whether $H$ itself  is cyclic.

\vspace{5mm}
{\bf 2.2.1} Is there some $i\le m+n-1$ such that $e:=a^i$ is idempotent (1.2), i.e. satisfies $e^2=e$? If $e$ is in the tail then either $e^2$ is in the tail itself, or in the body. Clearly in both cases $e^2\neq e$. Hence the only chance for $e$ to be idempotent is to be in the body, and so we try exponents $i$ of type $i=m+j$ where $0\le j<n$. If $e^2=e$ then necessarily $a^{2m+2j}=a^{m+j}$, hence $2m+2j\equiv m+j\ (mod\ n)$ by (4), hence 

$$m+j\equiv 0\ (mod\ n).$$

\noindent
Since $0\le j<n$, there is a unique $j_0$ with $0\le j_0<n$ satisfying $m+j_0\equiv 0\ (mod\ n)$. The argument is reversible and thus establishes the following.

\begin{itemize}
    \item[(5)] {\it If $\langle a\rangle\simeq C_{m,n}$, then there is a unique idempotent $e\in\langle a\rangle$\\
    (namely $e=a^{m+j_0}$ where $m+j_0\equiv 0\ (mod\ n)$).}
\end{itemize}

\noindent
 Let $E(S)$ be the set of all idempotents of the semigroup $S$. 
 If the semigroup $S$ is finite, then $E(S)\neq\es$. 
This is an immediate consequence of $(5)$; finiteness is crucial  viewing that $(\{1,2,3,...\},+)$ has no idempotents.

\vspace{3mm}
{\bf 2.2.2} The {\it order} of an element $x\in S$ is $o(x):=|\langle x\rangle |$. Thus if $\langle x\rangle\simeq C_{m,n}$ then $o(x)=m+n-1$. Let $\{a,b,...,c\}$ be a generating set of the c.f. sgr $S$. Refining $(2)$ it holds that the set of all elements $a^ib^j\cdots c^k$, where $(i,j,..,k)$ ranges over
$(0,0,..,0)\ <\ (i,j,..,k)\le (o(a),o(b),..,o(c))$, exhausts $S$. Therefore

$$(6)\quad |S|\ \le\ (o(a)+1)(o(b)+1)\cdots (o(c)+1)\ -1.$$

\noindent
We say that $\{a,b,..,c\}$ is {\it trimmed} if $\le$ in (6) is $=$. It then holds (why?) that $S$ is isomorphic to a direct product of cyclic semigroups matching the types of $\langle a\rangle,\langle b\rangle,...\langle c\rangle$. Does each c.f. sgr $S$ have a trimmed generating set?  In 3.6 and 7.3 we resume this issue.

\vspace{5mm}

{\bf 2.3} Each morphism $f:\langle a\rangle\to \langle b\rangle$ is determined by its value on $a$; indeed for all $i\ge 1$ it holds that 

$$af=b^k\ \Ra\ a^if=(af)(af)\cdots (af)=(b^k)^i=b^{ki}\    (using\ (1)).$$

\noindent
The converse fails in that for some "bad" exponents $k$ there might be {\it no} morphism $g$ satisfying $ag=b^k$. Thus, putting 
$a^ig:=b^{ki}$ may not be well-defined in the sense that for some $i\neq j$ one may have $a^i=a^j$, yet $b^{ki}\neq b^{kj}$.
Pleasantly, {\it if} $f$ is well-defined, then $f$ "automatically" is a morphism: 

$$(a^ia^j)f=(a^{i+j})f=b^{k(i+j)}=b^{ki+kj}=b^{ki}b^{kj}=(a^if)(a^jf)$$

\noindent
It hence suffices to uravel the conditions for well-definedness.
This is the precise state\footnote{Surprisingly I didn't find this in the literature. Theorem 1 will be a key ingredient in Section 9. The acronyms (SR1),(SR2) will be explained there as well.} of affairs:

\vspace{5mm}
{\bf Theorem 1: }{\it Let $\langle a\rangle$ and $\langle b\rangle$ be cyclic semigroups of types $C_{m,n}$ and $C_{m',n'}$ respectively, and let $k\ge 1$ be a fixed integer. Then $a^if:=b^{ki}$ is well-defined (and hence yields a unique morphism $f:\langle a\rangle\to \langle b\rangle$) iff it holds that 
\begin{itemize}
    \item[(SR1)] $m'\le km$, and
    \item[(SR2)] $n'\ divides\ kn$.
\end{itemize}   }

\vspace{5mm}
 Before we embark on the proof, let us find all $k$'s that satisfy (SR1) and (SR2) if
$C_{m,n}:=C_{2,10}$ and $C_{m',n'}:=C_{13,6}$. Since (SR1) becomes $13\le 2k$ we find that $k\ge 7$. Further (SR2) implies that $6$  divides $10k$, and so $k=3,6,9,12,...$. Together
with $k\ge 7$ this yields $k=9,12,15,18$ (dropping 21,24,... since $b^{21}=b^{15},\ b^{24}=b^{18}$ etc). Hence by Theorem 1 there are exactly four "exquisite" $k$'s, i.e. leading to morphisms $f:C_{2,10}\to C_{13,6}$; in formulas $Exq(2,10,13,6)=\{9,12,15,18\}$.

\vspace{5mm}
{\it Proof of Theorem 1.} 
Using (SR1) and (SR2) we first show that from $a^i=a^j\ (i\neq j)$ follows $b^{ki}=b^{kj}$. Indeed, $a^i=a^j$ implies $i\equiv j\ (mod\ n)$ by (4), and so (in view of $i\neq j$) there is an integer $\alpha\neq 0$ with $i-j=\alpha n$, hence $k(i-j)=k \alpha n$. By (SR2) there is $\beta\neq 0$ with $k(i-j)=\beta n'$, hence $ki\equiv kj\ (mod\ n')$.

Furthermore $a^i=a^j\ (i\neq j)$  implies $i,j\ge m$, hence $ki,kj\ge km\ge m'$ by (SR1). This together with $ki\equiv kj\ (mod\ n')$ and (4) forces $b^{ki}=b^{kj}$.

\vspace{3mm}
Conversely, suppose that (SR1) fails, i.e. $m'>km$. Then $a^m=a^{m+n}$ yet $b^{km}\neq b^{k(m+n)}$ because $km\neq k(m+n)$ and $b^{km}$ is in the tail of $\langle b\rangle$ in view of $km<m'$.

Likewise, suppose that (SR2) fails, i.e. $n'$ does not divide $kn$. Consider again $a^m=a^{m+n}$. Because $k(m+n)-km=kn$ is no multiple of $n'$ by assumption, we conclude $k(m+n)-km\not\equiv 0\ (mod\ n')$. Hence $k(m+n)\not\equiv km\ (mod\ n')$, hence $b^{k(m+n)}\neq b^{km}$ by (4). $\square$

\vspace{4mm}
{\bf 2.3.1} Does there always {\it exist} at least one morphism $f:C_{m,n}\to C_{m',n'}$? Yes, if $e$ is the unique idempotent of $C_{m',n'}$ then $(\forall i)\ a^if:=e$ evidently is a morphism. Referring to the example preceeding the proof of Theorem 1, which $k\in Exq(2,10,13,6)=\{9,12,15,18\}$ yields the idempotent $e=b^k$?

\vspace{5mm}
{\bf 2.4} Recall from 1.2 the two extreme types of idempotents  $\zero$ and $\one$. It is easy to see  that a semigroup can have {\it at most one zero} and {\it at most one identity}.
For instance $(S,\ast)=(\Z_n,\odot)$  has {\it both} $\zero$ and $\one$, but each semigroup 
$C_{m,n}$ with a tail (so $m>1$) and nontrivial body (so $n>1$) has {\it neither} $\zero$ nor $\one$.
It is also evident (check) that for each semigroup $S$ it holds that

$$ (7)\quad (\exists \zero,\one\in S\ and\ \ \zero=\one)\ \LRa\
|S|=1,$$

\noindent
in which situation $S$ is called {\it trivial}. 

\vspace{2mm}
Recall from 2.2.1 that the unique $j_0\in\{0,1,..,n-1\}$ with $m+j_0\equiv 0\ (mod\ n)$ yields the unique idempotent $e=a^{m+j_0}$ of $\langle a\rangle$. It follows from
$m+j_0+m+k\equiv m+k\ (mod\ n)$ and (4)  that $ea^{m+k}=a^{m+j_0+m+k}=a^{m+k}$ (for all $k\ge 0$), and so:

\begin{itemize}
    \item[(8)] {\it The unique idempotent $e$ of $C_{m,n}$ is an identity of
    the body of $C_{m,n}$\\ (but  not of $C_{m,n}$, unless $m=1$).}
\end{itemize}

\noindent
For the remainder of Section 2 we are concerned with $\zero$ (and its bigger brothers, ideals), while $\one$ takes the stage in Section 3.

\vspace{5mm}

{\bf 2.5} A semigroup $N$ is {\it nil} if it has a zero $\zero$ and for each $x\in N$ there is some $k\ge 1$ with $x^k=\zero$. The most extreme examples of nil sgr are the {\it zero semigroups} where  $xy=\zero$ for all $x,y\in N$.
Notice that  a cyclic sgr  $C_{m,n}$ is nil iff $n=1$. 
As another example, say $N=\langle a,b\rangle $ has $\zero$ and $a^{123}=b^{321}=\zero$. If $N$ is  {\it commutative}, then $N$ is nil (why?).

\vspace{5mm}
{\bf Theorem 2: }{\it Let $N$ be a finite semigroup ($|N|=k$) with a zero $\zero$.\\ Then $N$ is nil iff $E(N)=\{\zero\}$. In the latter case $x^k=\zero$ for all $x\in N$.}

\vspace{5mm}
{\it Proof.}  Each idempotent $e\neq \zero$ in any sgr $N$ with $\zero$ satisfies $e^t=e\neq \zero$ for all $t\ge 1$. Hence each nil sgr $N$ has $E(N)=\{\zero\}$. Conversely, assume $E(N)=\{\zero\}$. Recall that $k=|N|$. It is clear  that $x^k$ is in the body of $\langle x\rangle$ for all $x\in N$. By assumption the unique idempotent in each body is $\zero$. It follows from (7) and (8) that the body itself is $\{\zero\}$.
$\square$

\vspace{4mm}
{\bf 2.6} For nonempty subsets $X,Y$ of a sgr $S$ we put 

$$XY:=\{xy:\ x\in X,\ y\in Y\}.$$

\noindent
For instance it holds that $X$ is a ssgr of $S$ iff $X^2:=XX\s X$. The next, somewhat technical result, caters for Theorem 4 below.

\vspace{6mm}
{\bf Lemma 3: }{\it If $H$ is a commutative semigroup with $|H|=t$\\ then $H^t\s H^2 E(H)$.}
\vspace{5mm}

{\it Proof.} Take any $x\in H^t$, so $x=h_1h_2\cdots h_t$ for some (not necessarily distinct) $h_i\in H$. For all $1\le i\le t$ put $p_i:=h_1h_2\cdots h_i$.

{\it Case 1:} $p_t=e\in E(H)$. Then $x=p_t=eee\in H^2 E(H)$.

{\it Case 2:} All $p_1,...,p_t$ are distinct. Then $\{p_1,...,p_t\}=H\supseteq E(H)\neq\es$, and hence there is some $p_i=e\in E(H)$. By Case 1 we can assume $i<t$, which makes $h_{i+1}\cdots h_t$ well-defined. By commutativity $x=p_i h_{i+1}\cdots h_t=h_{i+1}\cdots h_t p_ie\in H^2 E(H)$.

{\it Case 3:} $p_i=p_j$ for some $1\le i<j\le t$. Then $p_i=p_j=p_ih_{i+1}\cdots h_j=:p_iz$, hence
$p_iz=p_iz^2,\ p_iz^2=p_iz^3$, and so forth. Because $\langle z\rangle$ contains an idempotent $e$, we conclude $p_i=p_ie$. Hence $x=p_i h_{i+1}\cdots h_t= h_{i+1}\cdots h_tee\in H^2 E(H)$. $\square$

\vspace{5mm} 
{\bf 2.7} Recall that a binary relation $R\s S\times S$ on any set $S$ is a {\it preorder} if it is transitive and reflexive. It is a {\it partial order} if additionally it is antisymmetric. A partial order $R$ is a {\it total order} if for all $a,b\in S$ it holds that $(a,b)\in R$ or $(b,a)\in R$.
A partial (or total) order $R$ is {\it strict}\footnote{Thus, "strictly speaking", a strict partial order is no partial order!} if instead of reflexive it is irreflexive (i.e. $(x,x)\not\in R$ for all $x\in S$). 

The {\it poset} (:=partially ordered set) of all divisors of 18 is rendered in Figure 1A. Divisibility also is important  for c. semigroups $S$, but in contrast to Fig.1A a "multiple" of $a$ will be {\it smaller} than $a$. So $1$ is the bottom element in Fig. 1A, but the identity $\one\in S$ is on top in Fig. 1B.

\vspace{3mm}
{\bf 2.7.1} For any commutative semigroup $S$ and all $a,b\in S$ put

 \begin{itemize}
     \item[(9)] $a<_{\cal J}b\ :\LRa\ (\exists x\in S:\ a=bx)\quad$ {\it and}   
 $\quad a\le_{\cal J}b:\ \LRa\ (a<_{\cal J}b\ or\ a=b).$
  \end{itemize}

 \noindent
 (As to $\cal J$ alone, see Section 10.5.) If $a\le_{\cal J}b$ then $a$ is a {\it  multiple} of $b$. If $a<_{\cal J}b$ then $a$ is a {\it  proper multiple} of $b$. It is evident that  $<_{\cal J}$ is transitive, and so $\le_{\cal J}$ is a preorder.
 
 Depending on the underlying sgr, $<_{\cal J}$ may enjoy extra properties.
For instance, take the ssgr  $S:=\{[0],[1],[3],[9],[10],[12]\}$ of $\Z_{18}=\{0,1,..,17\}$ in Fig. 1B (ignore the square brackets [ ] for now). Brute-force one checks that here $<_{\cal J}$ is  antisymmetric, i.e. the simultaneous occurence of $x<_{\cal J}y$ and $y<_{\cal J}x$ is possible {\it at most} when $x=y$. Actually $S$ has an identity, and so it holds for {\it all} $x\in S$ that\footnote{In semigroups $T$ without identity this may, or may not hold. For instance $T:=S\setminus\{[1]\}$ happens to be a ssgr and in it $<_{\cal J}$  satisfies $[12]<_{\cal J}[12]$, but {\it not} $[3]<_{\cal J}[3]$ (why?).} $x<_{\cal J} x$ (viewing that $\one x=x$).

\vspace{5mm}
{\bf Theorem 4: }{\it Let $N$ be a commutative finite nilsemigroup. 
\begin{itemize}
    \item[(a)]$<_{\cal J}$ is antisymmetric on $N$, and a strict partial order on $N\setminus\{\zero\}$. 
    \item[(b)]$\le_{\cal J}$ is a partial order on $N$ with  smallest element $\zero$.
    \end{itemize}  }

\vspace{6mm}

{\it Proof.} (b) immediately follows from (a). As to $<_{\cal J}$ being irreflexive on $N\setminus\{\zero\}$, take any $a\in N$ with $a<_{\cal J}a$, i.e.  $a=ax$ for some $x\in S$. Then $a=ax=aax=\cdots a^nx$ for all $n\ge 1$. Since $N$ is nil, some $a^n=\zero$, and so $a=\zero x=\zero$.

In order to show the antisymmetry of $<_{\cal J}$ on $N$ we show that from $a<_{\cal J}b$
and $b<_{\cal J}a$ follows $a=b\ (=\zero)$.  Let $t:=|N|$. There are elements $x,y\in N$ such that $b=xa$ and $a=yb$. Using commutativity this leads to

$$a=yb=yxa=y^2xb=y^2x^2a=\cdots \in N^t.$$

 But $N^t\s N^2 E(N)=N^2\{\zero\}=\{\zero\}$ by Lemma 3, and so $a=b=\zero$. $\square$

 \vspace{3mm}
 {\bf 2.7.2} It is evident that all nilsgr  with at most two elements are zero sgr. 
 We aim to show that any 3-element nilsgr $N$, which is no zero sgr, is isomorphic to $C_{3,1}$. To begin with,  $(N,\le_{\cal J})$ must be (why?) a 3-element chain, say $0<_{\cal J}y<_{\cal J}x$. If we can show that $N=\langle x\rangle$ then $N\simeq C_{3,1}$. First, $x^3=\zero$ by Theorem 2. It  remains to show that $x^2=y$. From $y<_{\cal J}x$ follows that either $y=x^2$ or $y=xy$. The latter yields the contradiction 
 \noindent
 $y=xy\Ra xy=x^2y\Ra x^2y=x^3y=\zero\Ra y=\zero$. 

\vspace{3mm}
{\bf 2.7.3} In any poset one says that $b$ is an {\it upper cover} of $a$ if $a<b$ and there is no $c$ with $a<c<b$.
  Consider $N:=\{0,2,4,6,8,10,12,14\}$, which is (why?) a nil ssgr of $(\Z_{16},\odot)$. In order to get  the diagram of the poset $(N,\le_{\cal J})$ we need, for each $x\neq 0$, the set $PM(x):=\{xy:\ y\in N\}\setminus \{x\}$ of all proper multiples of $x$. Clearly $x$ is an upper cover of $0$ in $(N,\le_{\cal J})$ iff $PM(x)=\{0\}$. One checks that only $PM(8)=\{0\}$. The upper covers of $8$ are exactly the elements $x$ with $PM(x)=\{0,8\}$; it turns out that they are $x=4$ and $x=12$. Next we need those $x$ with $PM(x)\s\{0,8,4,12\}$. These are exactly the $x$'s whose lower covers are to be found among the maximal elements $\{m_1,m_2,..\}$ of the poset so far. Here $\{m_1,m_2\}=\{4,12\}$ and the qualifying $x$'s happen to be all remaining elements, i.e. $2,6,10,14$. See Figure 1C.

 \vspace{1cm}
\includegraphics[scale=0.78]{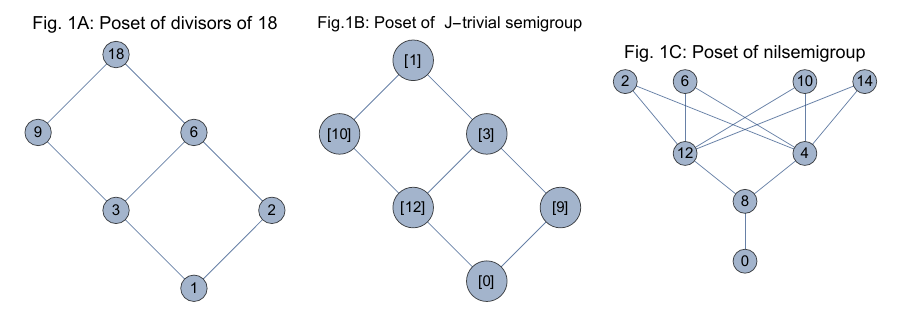}

\vspace{4mm}
{\bf 2.8} Apart from c.f. nilsemigroups, semilattices (in Sec.5), and sporadic\footnote{Fig. 1B shows some ssgr of $\Z_{18},\odot)$ (ignore the brackets [,] ), for which $\le_{\cal J}$ happens to be a partial order. More context will be provided in 10.6.} semigroups (Fig 1B), here comes another class of c. semigroups where the relation $\le_{\cal J}$ is a partial order.

The {\it free commutative semigroup}, say on three generators $a,b,c$, is defined on the infinite set $F_3$ of all  {\it words} $w:=a^{i_1}b^{i_2}c^{i_3}$ where ${\bf 0}:=(0,0,0)<(i_1,i_2,i_3)$. By definition, two words represent different elements of $F_3$ iff they differ in at least one exponent. Clearly the multiplication  

$$a^{i_1}b^{i_2}c^{i_3} \ast a^{j_1}b^{j_2}c^{j_3}:=a^{i_1+j_1}b^{i_2+j_2}c^{i_3+j_3}$$

\noindent
is associative, and the ensuing sgr $(F_3,\ast)$ is isomorphic to $(\N^3\setminus\{\zero\},+)$. So far, so obvious. 

\vspace{3mm}
{\bf 2.8.1}
Yet there's more to come.
It holds that $a^{i_1}b^{i_2}c^{i_3}\le_{\cal J} a^{j_1}b^{j_2}c^{j_3}$ iff $(i_1,i_2,i_3)\ge (j_1,j_2,j_3)$, where $\le$ is the component-wise order of $\N^3$. We  mend this inconvenience by defining

$$ a^{i_1}b^{i_2}c^{i_3}\le_c a^{j_1}b^{j_2}c^{j_3}\ :\LRa\ (i_1,i_2,i_3)\le (j_1,j_2,j_3)$$

 \noindent
 and calling $\le_c$  the {\it component-wise order} on $F_3$.
 All of this carries over to $F_k$ for any $k\ge 1$. 
 
 Apart from the component-wise order, $F_k$ carries  [FP] a certain {\it military (total) order} $\le_M$ which will be crucial in Section 6. We define it now in order not to be distracted later. For $k=3$ and $a,b,c$ instead of $a_1,a_2,a_3$ it starts like this:

$$a<_Mb<_Mc<_Ma^2<_M ab <_M ac <_M b^2 <_M bc <_M  c^2 <_M a^3 <_M a^2b<\cdots$$

\noindent
Thus for general words $v,w\in F_k$ it holds that
$v<_M w$ if either $|v|<|w|$ (smaller length), or $|v|=|w|$ but $v$ is 
lexicographic smaller\footnote{ For instance $ab<_M ac$ because in a lexicon the word $ab$ would preceed $ac$. But also $b^2<_M a^3$ although in a lexicon $a^3=aaa$ preceeds $b^2=bb$! } than $w$ (assuming that lexicographic $a_1$ comes before $a_2$, which comes before $a_3$, etc). Of course if $v<w$ component-wise, then $v<_M w$ (since $|v|<|w|$). 
By elementary combinatorics the number of words $a_1^{i_1}a_2^{i_2}..a_k^{i_k}\in F_k$ of length $n$ equals $\binom{n+k-1}{n}$; for instance $\binom{2+3-1}{2}=6$, corresponding to the 6 words $a^2,ab,ac,b^2,bc,c^2$ above.


\vspace{5mm}
{\bf 2.9}  A nonempty subset $I$ of a   semigroup $S$ is an {\it ideal} of $S$ if $IS\s I$ and $SI\s I$. In particular $II\s I$, and so each ideal is a ssgr, but not conversely. 
Each ideal $I$  of $S=(S,\cdot)$ gives rise to a smaller sgr with zero. Namely,  on the set-system

$$S/I:=\{\{x\}:\ x\in S,\ x\not\in I\}\cup\{I\}$$ 

\noindent
we define a binary operation $\ast$ as 

$$
y_1\ast y_2:=
\begin{cases}
    \{x_1\cdot x_2\}, & \text{if } y_1=\{x_1\},\ y_2=\{x_2\}\ and\ x_1\cdot x_2\not\in I \\
    I, & \text{\it otherwise }
\end{cases}
$$

\noindent
With some care (try) one verifies that $\ast$ is associative. The ensuing semigroup,  the so-called {\it Rees quotient} $(S/I,\ast)$, has the zero $\zero:=I$. Emphasizing that $\setminus$ (meaning set-complement) must not be confused with $/$ we can thus state that

$$The\ universe\ of\ S/I\ is\ \{\zero\}\cup (S\setminus I)$$

\vspace{3mm}
{\bf 2.10.1} Let $I_1\s F_3=\langle a,b,c\rangle$ be the ideal of all words $w$ which are multiples  of at least one of

$$a^3,\ b^4,\ c^5,\
 a^2b^2c^3,\ ac^4,\ b^3c^2,\ ab^3.$$

\noindent
For instance $a^2b^2c^4$ is a (proper) multiple of $a^2b^2c^3$. Hence, apart from $\zero$, the elements of the Rees quotient $F_3/I_1$ are exactly the words $w$ in $F_3\setminus I_1$, i.e. the ones that {\it simultaneously} satisfy (for visibility $\le:=\le_c$)
$$a^3\not\le w,\ ,\ b^4\not\le w,\ c^5\not\le w,\ a^2b^2c^3\not\le w,\ ac^4\not\le w,\ b^3c^2\not\le w,\ ab^3\not\le w$$

It turns out (research in progress for arbitrary $F_k\setminus I$)) that $F_3\setminus I_1$ can be rendered in a compressed format, namely (using obvious notation):

$$F_3\setminus I_1=\{a^2b^2c^{\le 2}\}\uplus \{b^{\le 2}c^{\le 4}\}\uplus \{b^3c^{\le 1}\}\uplus
\{ab^{\le 2}c^{\le 3}\}\uplus \{a^2b^{\le 1}c^{\le 3}\}$$

\noindent
It follows that $|F_3\setminus I_1|=3+(3\cdot 5-1)+2+3\cdot 4+2\cdot 4=39$. In the 40-element sgr $(F_3/I_1,\ast)$ it e.g. holds that $ab\ast bc=ab^2c$, whereas $ab^2\ast bc=\zero$. In fact  $F_3/I_1$ happens\footnote{What is a sufficient and necessary condition for $F_k/I$ to be a nilsemigroup?} to be a nilsemigroup.

\vspace{3mm}

{\bf 2.10.2} Each\footnote{This fails for the infinitely many ideals $2\Z\supseteq 4\Z\supseteq 8\Z\supseteq 16\Z\supseteq\cdots$ in $(\Z,\odot)$.} {\it finite} semigroup $S$ has a smallest ideal $K(S)$, which one calls the {\it kernel} of $S$. This hinges on the fact that, while the intersection of two ssgr can be empty (say $\{e\}\cap \{e'\}=\es$ for  distinct idempotents $e,e'$), this cannot happen for ideals $I,I'$. Indeed, if $x\in I,\ y\in I'$, then $xy\in II'\s I\cap I'$. By induction the intersection of finitely many ideals must be an ideal. In particular, the intersection $K(S)$ of {\it all} ideals of our finite semigroup $S$ is an ideal, and obviously its  smallest one.

For instance, if $S$ has a $\zero$ then $K(S)=\{\zero\}$. As to the kernel of $\langle a\rangle\simeq C_{m,n}$,  
let $I\s\langle a\rangle$ be any ideal. Pick any $a^t\in I\neq\es$. Then 

$$I\supseteq\{a^k:\ k\ge t\}\stackrel{why?}{\supseteq}\{a^{m+i}:\ i\ge 0\}=H,$$

\noindent
where $H$ is the body of $\langle a\rangle$. Since $H$ itself is an ideal of $\langle a\rangle$, we conclude $K(\langle a\rangle)=H$.

\section{ Abelian subgroups }

 Groups are important for us because many subsemigroups of  semigroups "can't help" being groups. Here a summary of  the key results. The kernel of each c.f. semigroup is an Abelian (:=commutative) subgroup. The structure of the group $(\Z_n,\odot)^{inv}$ is unraveled (while the whole of  $(\Z_n,\odot)$ has to wait until Section 7). The Fundamental Theorem for finite Abelian groups is stated but not proven. Instead, a seldom proven fact is verified in detail: Two finite Abelian $p$-groups are isomorphic iff they have the same number of elements of of each prime power order.

\vspace{3mm}
{\bf 3.1}
A {\it nontrivial} semigroup $S$ with identity $\one$ is called a {\it monoid}. In this case the set of {\it invertible} elements is defined as

$$S^{inv}:=\{\ a\in S:\ (\exists b\in S)\ ab=ba=\one\}.$$

\noindent
Clearly $\one\in S^{inv}$. It is a standard exercise (e.g. in  linear algebra courses concerned with square matrices $a,b$) to show that there is at most one $b$ satisfying $ab=ba=\one$. If $b$ exists, it hence is well-defined to write $a^{-1}$ for $b$ and call $a^{-1}$ the {\it inverse} of $a$. The argument that $S^{inv}$ is a ssgr of $S$ is also well known.
Further one verifies that the direct product of monoids is a monoid  and that

$$(10)\quad (S_1\times\cdots \times S_t)^{inv}=S_1^{inv}\times\cdots \times S_t^{inv}.$$

{\bf 3.1.1}
Suppose $S$ has a zero and $\zero\in S^{inv}$. Then $\zero\zero^{-1}=\one$, as well as $\zero\zero^{-1}=\zero$. But then $\zero=\one$, which contradicts the assumption that $S$ is nontrivial  (see (7)) .  
\noindent
We see that a   monoid $S$ with $\zero$  satisfies $S^{inv}\s S\setminus \{\zero\}$. 

One calls $x\in S\setminus\{\zero\}$ a {\it zerodivisor} if $xy=\zero$ for some $y\in S\setminus\{\zero\}$. A zerodivisor $x$ cannot be invertible because this gives  the contradiction 

$$y=\one y=(x^{-1}x)y=x^{-1}(xy)=x^{-1}\zero=\zero.$$

\noindent
Therefore, letting $NZD(S)$ be the set of all {\it non}-zerodivisors  we conclude that

$$(11)\quad S^{inv}\s NZD(S).$$

\vspace{5mm}
{\bf 3.2}
If the monoid $G$ satisfies $G^{inv}=G$, then $G$ is called a {\it group}.
Because $(S^{inv})^{inv}=S^{inv}$, each monoid $S$ yields a group $G:=S^{inv}$.

\begin{itemize}
  \item[(12)] For each semigroup $C_{m,n}=\langle a\rangle$ the following are equivalent.  
\begin{itemize}
\item[(i)] $\langle a\rangle$ is a group
\item[(ii)] $\langle a\rangle$ is a monoid
    \item [(iii)] $\langle a\rangle$ has index  $m=1$ (no tail)
    \end{itemize}    
\end{itemize}

The implication $(i)\Ra (ii)$ is trivial. As to $(ii)\Ra (iii)$, since some  $a^{m+j}$ is the {\it only} idempotent in $\langle a\rangle$, we get $a^{m+j}=\one$. Because of $a^i=a^i\one=a^{i+m+j}$ each $a^i\in\langle a\rangle$ is in the body of $\langle a\rangle$, i.e. there is no tail. As to $(iii)\Ra (i)$, if $m=1$ then it follows from $a^{m+n}=a^m$ that $a\cdot a^n=a$.
 Hence $\one:=a^n$ is an identity of $\langle a\rangle$ and $(a^i)^{-1}=a^{n-i}$ for all $1\le i\le n-1$. $\square$

 \vspace{3mm}
Recall from Theorem 2 that a finite sgr with $\zero$ is  nil iff $\zero$ is the {\it only} idempotent. In similar fashion finite groups can be characterized:

\vspace{5mm}
{\bf Theorem 5: }{\it Let $G$ be a finite monoid. Then $G$ is a group iff $E(G)=\{\one\}$.}

\vspace{5mm}
{\it Proof.} If $G$ is a group then each $e=e^2\in E(G)$ satisfies $e=e\one=e(ee^{-1})=(ee)e^{-1}=ee^{-1}=\one$.  Conversely assume that $E(G)=\{\one\}$. Then for any $a\in G$ the unique idempotent in $\langle a\rangle$ must be $\one$. So $\langle a\rangle$ contains an identity and thus by $(ii)\LRa (i)$ in (12) there is an inverse $a^{-1}$ of $a$ within the group $\langle a\rangle$. Clearly $a^{-1}$ is also an inverse of $a$ within $G$.
$\square$

\vspace{3mm}
One calls a commutative sgr {\it cancellative} if from $ab=ac$ always follows $b=c$. A large part of [RG] is dedicated to finitely generated cancellative c. sgr. $S$. It is easy to see\footnote{Show that each cyclic ssgr $\langle a\rangle\s S$ has index $m=1$ and proceed as in the proof of Thm.5.}  that each {\it finite} such $S$ must be a group.

\vspace{5mm}
{\bf 3.3} Let us show that the kernel $K$ (= its body by 2.10.2) of the cyclic sgr  $C_{m,n}=\langle a\rangle$ is a cyclic group. In fact we claim that $K=\langle ea\rangle$.
First, $ea=a^{m+j_0+1}$ belongs to $K$ (see (5)). Therefore $\{ea,ea^2,ea^3,..\}$ exhausts $K$. But $ea^k=e^ka^k=(ea)^k$, and so  $K=\langle ea\rangle$. 
From $|K|=n$ follows $K\simeq C_{1,n}$. In the sequel we put

$$C_n:=C_{1,n}$$

Whenever a ssgr $U$ of a sgr $S$ happens to be a group "on its own", one calls $U$ a {\it subgroup} of $S$. As seen above ($U=K$), the identity of $U$ need not be an identity for $S$. Generalizing $K(C_{m,n})\simeq C_n$ the following holds.

\vspace{2mm}

{\bf Theorem 6: }{\it The kernel of each commutative finite semigroup $S$ is an Abelian subgroup.}

\vspace{2mm}
{\it Proof.} Since $K:=K(S)$ is a ssgr, it contains an idempotent $e$. If we can show that $e$ is an identity of $K$, then $e$ is unique. Then Theorem 5 implies that $K$ is a group. 

We first show that $\{e\}K$ is an ideal of $S$. Thus let $ex\in \{e\}K$ and $y\in S$ be arbitrary. Then $xy\in K$, and so (by commutativity) $y(ex)=(ex)y=e(xy)\in \{e\}K$. Since $K$ is the smallest ideal, it follows from  $\{e\}K\s K$ that $\{e\}K=K$. Since each element of $K$ is of type $ex$, and  $eex=ex$, we see that $e$ is  an identity of $K$. 
$\square$

\vspace{3mm}

\noindent
{\bf 3.4} Here and in 3.5 we study  the group $\Z_n^{inv}$. It foreshadows aspects of  general 
 finite Abelian groups to be dealt with in 3.6 and 3.7. We start by strengthening (11) for $S:=\Z_n$:

\vspace{5mm}
{\bf Theorem 7:} $\Z_n^{inv}=NZD(\Z_n)$.

\vspace{5mm}
{\it Proof.}  By (11) it suffices to show that each $a\in NZD(\Z_n)$ is in $\Z_n^{inv}$. We claim that this reduces to showing that $f:\Z_n\to\Z_n:\ x\mapsto xa$ is injective. That's because injectivity implies bijectivity in view of $|\Z_n|<\infty$, and so there is $x\in\Z_n$ with $xf=\one$. Hence $x=a^{-1}$, and so $a\in\Z_n^{inv}$.

As to establishing the injectivity of $f$,
we exploit the fact that $(\Z_n,\odot)$ is a reduct of the\footnote{Despite rings being more complex than semigroups, for most readers the "rough" structure of $(\Z_n,+,\odot)$ may be more familiar than the fine structure of $(\Z_n,\odot)$ (which we are going to unravel).} ring $(\Z_n,+,\odot)$. From $xf=yf$, i.e. from $xa=ya$, follows by {\it ring distributivity} that $\zero=xa-ya=(x-y)a$.
Since $a$ is a non-zerodivisor this forces $x-y=\zero$, i.e. $x=y$, i.e. the injectivity of $f$.
$\square$

\vspace{4mm}
{\bf 3.4.1} Digging a bit deeper gives a concrete description of $NZD(\Z_n)$ in terms of the greatest common divisor $gcd(x,n)$ of two positive integers:

$$(13)\quad NZD(\Z_n)= \{1\le x\le n-1:\ gcd(x,n)=1\}$$

\noindent
In order to prove (13), we thus show that $gcd(x,n)>1$ iff $x$ is a zerodivisor. Suppose first that $gcd(x,n)=d>1$. Viewing the integers $\frac{x}{d}$ and $\frac{n}{d}$ as elements of $\Z_n$ (note $\frac{n}{d}\neq \zero$ since $d>1$) it holds that $x\cdot \frac{n}{d}=\frac{x}{d}\cdot n\equiv 0\ (mod\ n)$, so $x$ is a zerodivisor.

Conversely, if  $x\in\{1,2,...,n-1\}$ is a zerodivisor, then there is $y\not\equiv 0$ with $xy\equiv 0$. Hence $xy=kn$. Each maximal prime power $p_i^{\epsilon_i}$ dividing $n$ must divide $xy$. If each $p_i^{\epsilon_i}$ were to divide $y$,  then $y$ would be divisible by $n$, in contradiction to $y\not\equiv 0$. Hence some $p_i$ divides $x$, showing that $gcd(x,n)>1$. 
$\square$

\vspace{3mm}
The {\it Euler $\varphi$-function} is defined as

$$(14)\quad \varphi(n):=|\{1\le x\le n-1:\ gcd(x,n)=1\}|$$

\vspace{3mm}
Since $gcd(30,49)=1$, Theorem 6 combined with (13) implies that  $30\in \Z_{49}^{inv}$, but gives no indication of how to calculate $30^{-1}$! That can be done with the {\it Euclidean algorithm} [A,p.30], which in our case yields $1=18\cdot 30-11\cdot 49$. This implies $18\cdot 30\equiv 1\ (mod\ 49)$, and so $30^{-1}=18$ in $\Z_{49}$. 

\vspace{5mm}
{\bf 3.4.2}
Here comes a situation where the group $\Z_n^{inv}$ (whose fine structure is complicated, see 3.5.2) clarifies an aspect of the bland group $C_n$. 

\begin{itemize}
    \item[(15)] Let $C_n=\langle a\rangle=\{a,a^2,\ldots, a^n=\one\}$ be a cyclic group, and fix $i\in\{1,..,n\}$. Then $a^i$ generates $\langle a\rangle$ iff $gcd(i,n)=1$.
    In particular, $C_n$ has $\varphi(n)$ generators.
\end{itemize}

\noindent
To prove (15), note that $\langle a^i\rangle=\langle a\rangle$ iff
$(a^i)^j=a$ for some $j$. Hence iff $ij\equiv 1\ (mod\ n)$, hence iff
$i\in(\Z_n,\odot)^{inv}$. But $\Z_n^{inv}=\{1\le i\le n-1:\ gcd(i,n)=1\}$ by Theorem 6 and (13). $\square$

\vspace{3mm}
{\bf 3.5} For a deeper understanding of $\Z_n^{inv}$ we need again  the
  ring $(\Z_n,+,\odot)$. 

$$(16)\quad If\ n=p_1^{\gamma_1}p_2^{\gamma_2}\cdots p_t^{\gamma_t}\ then\ \Z_n \simeq \Z_{p_1^{\gamma_1}}\times\cdots\times  Z_{p_t^{\gamma_t}} \ as\ rings$$

{\it Proof.}  We take $n=60$ but it will be clear that the arguments generalize\footnote{ Generalization not only takes place from $\Z_{60}$ to $\Z_n$. Fact (16) is a special case of the so-called Chinese Remainder Theorem (CRT) that can be formulated for arbitrary rings. In this  setting establishing the surjectivity of $f$ becomes the hardest part of the proof; see [Co,p.102].}. It is easy to see that $[m]_{60}f:=([m]_3,[m]_4,[m]_5)$ yields a well-defined morphism
 $f:\Z_{60}\to\Z_3\times\Z_4\times\Z_5$. For instance
 $[29]_{60}f=([29]_3,[29]_4,[29]_5)=([2]_3,[1]_4,[4]_5)$.
If $[m]_{60}f=([0]_3,[0]_4,[0]_5)$ then $m$ is divisible by 3,4, and 5, whence by 60. It follows that $[m]_{60}=[0]_{60}$. As is well known, this implies the injectivity of $f$. As previously, "injective" implies "bijective" in view of $|\Z_{60}|=|\Z_3\times\Z_4\times\Z_5|$.   $\square$ 

\vspace{3mm}
{\bf 3.5.1} While the bijectivity  of $f$ was cheap, finding $f^{-1}(x)$ for concrete $x$ is nontrivial. The following  method (that comes in handy in 7.5) becomes the more economic the higher the number of $x$'s to be handled. For instance, for $f$ as above, find $e_1:=f^{-1}(([1]_3,[0]_4,[0]_5))$! Notice that $e_1'=[20]_{60}$ is "close" to $e_1$ in that at least $[20]_4=[0]_4$ and $[20]_5=[0]_5$ are correct. However $[20]_3=[2]_3\neq [1]_3$. Fortunately, it is immediate that   $e_1:=2{e_1'}=[40]_{60}$ does\footnote{In general  one needs to multiply $e_i'$ with its inverse in $\Z_{p_i^{\gamma_i}}$, the latter being calculated with the Euclidean algorithm.} the job. Similarly one finds that $e_2:=f^{-1}(([0]_3,[1]_4,[0]_5))=[45]_{60}$ and $e_3:=f^{-1}(([0]_3,[0]_4,[1]_5))=[36]_{60}$. With the "basis" $e_1,\ e_2,\ e_3$ it is easy to find all $f^{-1}(x)$ fast. 
Namely, from 

$$x=([a]_3,[b]_4,[c]_5)=([40a]_3,[45b]_4,[36c]_5)
\stackrel{why?}{=}$$
$$([40a+45b+36c]_3,[40a+45b+36c]_4,[40a+45b+36c]_5)$$ 

\noindent
follows that $f^{-1}(x)=[40a+45b+36c]_{60}$.

\vspace{3mm}
{\bf 3.5.2} Obviously "isomorphic as rings" implies "isomorphic as (multiplicative) monoids", and so (10) and (16) imply that

$$(17)\quad \Z_n^{inv} \simeq \Z_{p_1^{\gamma_1}}^{inv}\times\cdots\times  Z_{p_t^{\gamma_t}}^{inv} \ as\ Abelian\ groups$$

\vspace{3mm}
As to the "fine structure" of $\Z_n^{inv}$, by (17) it suffices to unravel the prime power case $\Z_{p^\gamma}^{inv}$. For starters, Theorem 6 and (13),(14) imply that $|\Z_{p^\gamma}^{inv}|=\varphi(p^\gamma)$, and clearly $\varphi(p^\gamma)=p^{\gamma}-p^{\gamma-1}=(p-1)p^{\gamma-1}$. This settles the cardinality of $\Z_{p^\gamma}^{inv}$, but what type of group is it?
Somewhat boring, it is mostly\footnote{For instance  $\Z_{16}^{inv}=\{1,3,5,7,9,11,13,15\}=\{1,3,-3^3,-3^2,3^2,3^3,-3,-1\}$,\\
and generally $\Z_{2^\gamma}^{inv}\simeq\langle -1\rangle\times\langle 3\rangle$. For German speaking folks [B,p.109ff] is recommended for a proof of (18).} cyclic. Specifically:

$$(18)\quad \Z_{p^\gamma}^{inv}\simeq C_{(p-1)p^{\gamma-1}}\ and\
\Z_{2^\gamma}^{inv}\simeq C_2\times C_{2^{\gamma-2}}\ (\gamma\ge 2)$$
 
\noindent 
Take $n=504$. Then it follows from $504=7\cdot 8\cdot 9$ and (17),(18) that

$$ \Z_{504}^{inv}\simeq \Z_{7}^{inv}\times \Z_{8}^{inv}\times \Z_{9}^{inv}\simeq
C_6\times (C_2\times C_2)\times C_6.$$

Similarly to the way that $(\Z_n,\odot)^{inv}$ (for varying $n$) covers many types of  {\it Abelian groups}, even more so $(\Z_n,\odot)$ will exhibit (in Sec.7) many concepts that pervade general  {\it commutative semigroups}.

\vspace{3mm}
{\bf 3.5.3} Here comes another, more "extreme" example of an Abelian group. The symmetric difference of sets, i.e. $A\triangle B:=(A\setminus B)\cup (B\setminus A)$, is a binary operation which is commutative (clear) and associative (not so clear). It hence yields a commutative semigroup $(\Pow(X),\triangle)$, which in fact  is a group since $\one:=\es$ is an identity and $A^{-1}=A$ for all $A\in\Pow(X)$. If $t:=|X|$, then $(\Pow(X),\triangle)\simeq (C_2)^t$ (why?).

\vspace{5mm}
{\bf 3.6} If $G$ is a finite group then (12) implies (how?) that each cyclic ssgr $C_{m,n}$ is of type $C_{1,n}=C_n$. Thus if $\langle x\rangle\simeq C_n$, then  $o(x)=n$. In the remainder of Section 3 we focus on {\it Abelian}  finite groups.

Particularly, in  3.6 we resume the trimmed generating sets $\{a,b,..,c\}$ of 2.2.2 in the scenario of Abelian f. groups $G$. The cumbersome formula (6) then simplifies to the extent that $\{a,b,..,c\}$ being trimmed amounts to $|G|=o(a)o(b)\cdots o(c)$. 

To fix ideas, consider $G':=C_{30}=\langle a\rangle$.
Apart from $\{a\}$, are there other trimmed generating sets? What about $\{a^{16},a^{25}\}$? It is a generating set since $(a^{16})^3a^{25}=a^{73}=a^{13}$, and $a^{13}$ generates $C_{30}$ by (15). Yet it is {\it not} trimmed\footnote{Here's a concrete example of non-unique generation: $a^{16}\cdot a^{25}=a^{11}$, as well as $(a^{16})^6(a^{25})^5=a^{221}=a^{11}$.} since $o(a^{16})o(a^{25})=15\cdot 6>30$. What about 
$X:=\{a^6,a^{10},a^{15}\}$?  Here the orders behave, i.e. $o(a^6)o(a^{10})o(a^{15})=5\cdot 3\cdot 2=30$, but is $X$ generating in the first place?
Yes\footnote{This generalizes: If $o(x)o(y)\cdots o(z)=|G|$, and the orders are pairwise coprime, then $\{x,y,..,z\}$ is a trimmed generating set. (However, in this case $G=\langle xy\cdots z\rangle$ is even 1-generated!) In contrast, 2 and 4 are not coprime, and things indeed go wrong:
$b^2,b^4\in G:=C_8=\langle b\rangle$ satisfy $o(b^2)o(b^4)=4\cdot 2=|G|$ yet do not generate $G$.} it is: $a^6\cdot a^{10}\cdot a^{15}=a^{31}=a$. Hence $G'\simeq\langle a^6\rangle\times\langle a^{10}\rangle\times\langle a^{15}\rangle\simeq C_5\times C_3\times C_2$.

\vspace{3mm}
 {\bf 3.6.1.  The Fundamental Theorem of finite Abelian groups} $G$ states\footnote{An elementary proof can be found in [Ar].} that each such $G$  has (possibly many) trimmed generating sets. An equivalent and more common phrasing is that $G$ is isomorphic to a direct product $C_{n_1}\times\cdots \times C_{n_t}$ of  cyclic groups. These cyclic groups need not be uniquely determined, not even their number $t$. Nevertheless, for the the minimum  occuring $t=t_{min}$, and for the maximum occuring $t=t_{max}$, there is (up to the order of factors) only one type of direct product. 

 Before we further discuss $t_{max}$ (3.6.2) and $t_{min}$ (3.6.3), observe that for each f. Abelian group
 $G$ and each prime $p$ the set $G_p:=\{x\in G:\ (\exists i\ge 0)\ o(x)=p^i\}$ is a subgroup\footnote{Any f. Abelian group $H$ of this type is called a {\it $p$-group}. An equivalent definition is that $|H|$ is a power of $p$.} of $G$. Furthermore $G$ is easliy seen to be isomorphic to the direct product of all its $p$-subgroups. The hard part of the Fundamental Theorem is to prove that each $p$-group is a direct product of {\it cyclic} $p$-groups.
 
 \vspace{3mm}
{\bf 3.6.2} To fix ideas, suppose that $G'\simeq G'_2\times G'_3\times G'_5\times G'_7\times G'_{11}$ (and that $|G_p'|=1$ for all $p\ge 13$). Suppose further that $|G_2'|=2^{15}$. That still leaves many\footnote{Specifically, there are 176 options because this is the number of {\it partitions} of the integer 15 (examples being 8+7 or 4+3+3+2+2+1).} options for its structure, say $G_2'\simeq C_{2^8}\times C_{2^7}=C_{256}\times C_{128}$ or $G_2'\simeq C_2\times (C_{4})^2\times (C_8)^2\times C_{16}$. Let's assume the latter, as well as $G'_3\simeq (C_3)^4\times C_{27}$ and $G'_5\simeq C_5\times (C_{25})^2$ and $G'_7\simeq C_7\times C_{49}$ and $G'_{11}\simeq (C_{11})^2$. It follows that $G'$ is a direct product of $6+5+3+2+2=18$ cyclic groups. and one can show that this is the maximum number of factors achievable. So $G'$ has $t_{max}=18$.

\vspace{3mm}
{\bf 3.6.3} In order to find $t_{min}$ for $G'$, let us write the cardinalities of the cyclic groups entering $G_2'$ as column 1 in the table below. Likewise $G_3'$ up to $G_{11}'$ determine the columns 2 to 5. The product of the numbers in the first row equals $n_6:=5821200$. Likewise the second row yields $n_5=46200$. It is clear that $n_5$ divides $n_6$ (and generally in such a scenario $n_{i-1}$ divides $n_i$). It holds (clear) that $G'\simeq C_{n_1}\times\cdots C_{n_6}$ and (less clear) 6 is the minimum number of factors possible. So $G'$ has $t_{min}=6$.

$$\begin{matrix}
    16 & 27 & 25 & 49 & 11 & \Ra & n_6= &5821200\\
    8 & 3 & 25 & 7 & 11 & \Ra & n_5= &46200\\
    8 & 3 & 5 &  &  & \Ra & n_4= &120\\
    4 & 3 &  &  &  & \Ra & n_3= &12\\
    4 & 3 &  &  &  & \Ra & n_2= &12\\
    2 &  &  &  &  & \Ra & n_1= &2
\end{matrix}$$

\vspace{4mm}
{\bf 3.7} Given the Cayley table of an Abelian group $G$, what can be done with it? Three questions spring to mind. The answers are easy (3.7.1), medium (3.7.2) and hard (3.7.3).

 \vspace{2mm} 
 {\bf 3.7.1} How to get a "reasonably small" generating set $\{x,y,..\}$ of $G$ ? The answer (not optimal but not bad either):  Let $x\in G$ be an element of maximum order. If $\langle x\rangle\neq G$, let $y\in (G\setminus\langle x\rangle)$ be such that $\langle y\rangle\setminus\langle x\rangle$ is large. 
If $\langle x,y\rangle\neq G$, let $z\in (G\setminus\langle x,y\rangle)$ be such that $\langle z\rangle\setminus\langle x,y\rangle$ is large. And so forth until $G$ is exhausted.

\vspace{3mm}
{\bf 3.7.2} How can one determine the isomorphism type of an Abelian $p$-group from its Cayley table?
Interestingly it turns out  that the number of elements of each order uniquely determines the isomorphism type of a $p$-group. For starters, we leave it to the reader to prove: 

 $$(19)\quad |\{x\in C_{p^e}:\ o(x)\le p^\alpha\}|=p^\alpha\quad (1\le\alpha\le e)$$

 As a hint, if say $C_{p^e}=C_{2^7}=\{a,a^2,...,a^{128}=\one\}$, then $o(x)\le 2^3$ iff $x=a^i$ for $i\in\{k\cdot 2^4:\ 1\le k\le 2^3\}$. For instance $o(a^{16})=8,\ o(a^{32})=4,\ o(a^{48})=8,\ o(a^{64})=2$.
 Furthermore, the following is clear:

 \begin{itemize}
     \item[(20)] Each $x=(x_1,...,x_t)\in C_{p^{e_1}}\times \cdots\times C_{p^{e_t}}=:G$ has $o(x)=max\{o(x_1),...,o(x_t)\}$. Consequently $o(x)\le p^\alpha$ iff $(\forall 1\le i\le t)\ o(x_i)\le p^\alpha$.
 \end{itemize}

Suppose again that $G$ in (20) has $|G|=p^{22}$. Further let it be known that $G$ has exactly $p^8$ elements $x$ of order $\le p$ (evidently $o(x)<p$ implies $x=\one$). Putting $\alpha=1$ in (20) it holds that $x\in G$ has $o(x)\le p$ iff $o(x_i)\le p$ for all $1\le i\le t$. Since $o(x_i)\le p$ by (19) occurs for exactly $p$ many $x_i\in C_{p^{e_i}}$, the number of $x\in G$ with $o(x)\le p$ is $p^t$. Hence by assumption $t=8$, and so 
$G=C_{p^{e_1}}\times \cdots \times C_{p^{e_8}}$, where wlog $e_1\ge\cdots \ge e_8\ge 1$.

Next  suppose  there are exactly\footnote{Hence  $p^{14}-p^8$ elements $x$ have $o(x)=p^2$. Yet the precise order (i.e. $o(x)=$ instead of $o(x)\le$) would be  distracting in the present argument.} $p^{14}$ elements $x\in G$ with $o(x)\le p^2$. Let $e_1\ge\cdots\ge e_i\ge 2$ and $e_{i+1}=\cdots =e_t=1$. Then by (19),(20) $o(x)\le p^2$ occurs exactly $(p^2)^i\cdot p^{8-i}=p^{8+i}$ times. Hence by assumption $i=6$, and so 
$G=C_{p^{e_1}}\times \cdots\times C_{p^{e_6}}\times (C_p)^2$ with $e_1\ge\cdots \ge e_6\ge 2$.

Suppose that $o(x)\le p^3$ occurs exactly $p^{19}$ times. Writing\\
$G=C_{p^{e_1}}\times \cdots\times C_{p^{e_i}}\times(C_{p^2})^{6-i}\times (C_p)^2$  with $e_1\ge\cdots\ge e_i\ge $3, it follows again from (19), (20) that 
 $o(x)\le p^3$ occurs exactly $(p^3)^i(p^2)^{6-i}(p)^2=p^{14+i}$ times. Hence $i=5$, and so $G=C_{p^{e_1}}\times \cdots\times C_{p^{e_5}}\times C_{p^2}\times (C_p)^2$  with $e_1\ge\cdots\ge e_5\ge 3$.

 Suppose that $o(x)\le p^4$ occurs exactly $p^{21}$ times. Writing\\
 $G=C_{p^{e_1}}\times \cdots\times C_{p^{e_i}}\times (C_{p^3})^{5-i}\times C_{p^2}\times (C_p)^2$  with $e_1\ge\cdots\ge e_i\ge 4$
it follows from (19), (20) that 
 $o(x)\le p^4$ occurs exactly $(p^4)^i(p^3)^{5-i}p^2 (p)^2=p^{19+i}$ times. Hence $i=2$, and so
 $G=C_{p^{e_1}}\times C_{p^{e_2}}\times (C_{p^3})^3\times C_{p^2}\times (C_p)^2$  with $e_1\ge e_2\ge 4$.

 From $|G|=p^{22}$ follows ad hoc that $G=C_{p^5}\times C_{p^4}\times (C_{p^3})^3\times C_{p^2}\times (C_p)^2$.

 \vspace{3mm}
 {\bf 3.7.3} The fact that the order statistics of a finite Abelian group determine its structure, is posed as Exercise 15 in [MKS,p.151]. Since the author (not a group theorist despite [W2] ) failed to google a proof, he had to invest a couple of hours to solve the exercise and communicate it in 3.7.2 in (hopefully) readible fashion.

I later learned an elegant short proof from Andrew Sutherland (google "groupprops, abelian, order statistics"). Specifically, the following is shown. Let $s_k$ be the number of cyclic group factors of order $p^k$, and let $t_k$ be the logarithm to base $p$ of the number of elements of order dividing $p^k$. It then holds that

$$ s_k=2t_k-t_{k+1}-t_{k-1}.$$

\noindent
Concerning the group in 3.7.2 it e.g. holds
that

$$s_3=2t_3-t_4-t_2=2\cdot 19-21-14=3.$$

\vspace{4mm}
{\bf 3.7.4} Knowing the Cayley table {\it and} the isomorphism of an Abelian $p$-group is still a long shot  from {\it finding} a trimmed generating set. For instance, even finding a generator of a cyclic group (with elements called 1,2,..,n) is nontrivial if only its Cayley table is known (try). In general, the state of the art concerning "finding" seems to be [S], which heavily relies on so called discrete logarithms. See also Section 8.8.

\section{Closure systems and implications}

Starting with closure systems (4.1), we turn to closure operators (4.2), and then to implications (4.3). Most of this caters for an aspect of semilattices discussed in Subsections 6.7 and 6.8.
\vspace{5mm}

{\bf 4.1} A {\it closure system} on a  set $X$ is a subset $\C$ of the powerset $\Pow(X)$ such that\footnote{Here $ \bigcap{\cal S}$  is the intersection of all sets in ${\cal S}$. Hence if ${\cal S}=\{A,B\}$ then $\bigcap{\cal S}=A\cap B$. In fact, if $|X|<\infty$ (which for us is the normal case), then (why?) $(\forall {\cal S}\s\C)\ \bigcap{\cal S}\in \C$ is equivalent to $(\forall A,B\in \C)\ A\cap B\in \C$.}  

$$(21)\quad X\in\C\quad and\quad (\forall {\cal S}\s\C)\ \bigcap{\cal S}\in \C.$$

 \noindent
 For $X=\{a,b,c,d\}$ a closure system $\C\s\Pow(X)$ is shown in Figure 2A. For say $A=\{a\}$ and $B=\{b,d\}$ in $\C$ it holds indeed that $A\cap B=\es\in\C$, but note that $A\cup B\not\in\C$.
For the time being ignore the labels 1,2,...,8 in Figure 2A.

\vspace{1cm}
\includegraphics[scale=0.9]{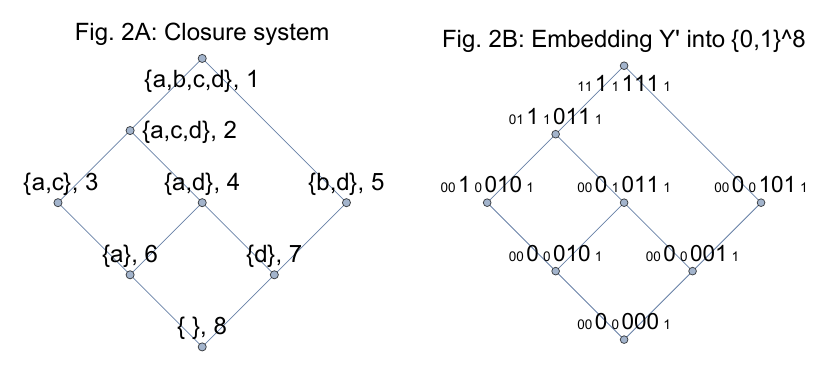}
\vspace{5mm}

{\bf 4.2} A {\it closure operator} on a set $X$ is a map $cl:\Pow(X)\to\Pow(X)$ which is extensive, idempotent and monotone, i.e. for all $U,V\in\Pow(X)$ with $U\s V$ it holds that

$$(22)\quad U\s cl(U)\ and\ cl(cl(U))=U\ and\ cl(U)\s cl(V).$$

\noindent
It is well-known [CLM,p.79] that closure systems $\C$ and closure operators $cl$ are two sides of the same coin. Specifically, put

$$(23)\quad cl_{\C}(U):=\bigcap\{V\in\C:\ V\supseteq Y\}\ and\ \C_{cl}:=\{U\in\Pow(X):\ cl(U)=U\}.$$

\noindent
Then $cl_{\C}$ is a closure operator and $\C_{cl}$ is a closure system. Furthermore   $\C_{cl_{\C}}=\C$ and $cl_{\C_{cl}}=cl$.

\vspace{5mm}
{\bf 4.3} Here comes a particular way to obtain a closure operator. For any $A,B\in\Pow(X)$ the ordered pair $(A,B)$ is called an {\it implication}. It will be more intuitive to write $A\to B$ instead of $(A,B)$, and to call $A$ the {\it premise}, and $B$ the {\it conclusion} of the implication $A\to B$. Let 

$$(24)\quad \Sigma:=\{A_1\ra B_1,\ldots, A_t\to B_t\}$$

\noindent
be a family of implications. Take any $U\s X$. If say $A_2, A_5$ are contained in $U$, we replace $U$ by $U':=U\cup B_2\cup B_5$. If now $A_3\s U'$, put $U'':=U'\cup B_3$. And so forth until we get a  set $V=U^{''\cdots''}$ that is stable in the sense that $V'=V$.
 This set we call the {\it $\Sigma$-closure} of $U$ and put $cl({\Sigma},U):=V$. It is evident that
 this yields a closure operator. The corresponding (see (23)) closure system $\C(\Sigma)$ hence consists of
all $\Sigma${\it -closed} sets $U$ in the sense that  $A_i\s U\Ra B_i\s U$ for all $1\le i\le t$. 

\vspace{3mm}
\noindent
For instance, if $X:=\{a,b,c,d\}$ and
 
$$\Sigma_1:=\Bigl\{\{a,b\}\ra \{c\},\ \{c\}\ra\{a\},\ \{b\}\ra\{d\}\Bigr\},$$

\noindent
then $cl({\Sigma_1},\{a,d\})=\{a,d\}$ and $cl({\Sigma_1},\{a,b\})=X$. One checks brute-force that $\C(\Sigma_1)$ is the closure system shown in Figure 2A.

\vspace{3mm}
{\bf 4.3.1} Recall that the set $\{0,1\}^t$ of all length $t$ bistrings naturally matches the powerset $\Pow(\{1,2,..,t\}$. For instance, if $t=7$, then $(0,1,0,1,1,0,1)$
goes to $\{2,4,5,7\}$. Bitstrings (=01-rows) generalize to {\it 012-rows}
such as\\

$$r:=(2,0,0,2,2,1,1,2,0,2).$$

\noindent
Here "2" is a don't-care symbole which can freely assume the value 0 or 1. Consequently $r$ can be identified with a {\it set-system} that contains exactly $2^4$ sets; for instance 

$$({\bf 1},0,0,{\bf 1},{\bf 0},1,1,{\bf 0},0,{\bf 1})\ "="\ \{1,4,6,7,10\}\in r.$$

\noindent
The boldface entries above are the ones that arose from the don't-cares. The article [W3] describes an algorithm which represents any closure system of type $\C(\Sigma)$ as a disjoint union of 012-rows. Consider say 

$$\Sigma_2:=\Bigl\{ \{a\}\ra \{b\},\ \{b,c\}\ra \{e\},\ \{a,e\}\ra \{b,d\},\ \{d\}\ra \{c\}\Bigr\}$$. 

\noindent
If we e.g. identify $(1,0,1,0,1)$ with $\{a,c,e\}$, and e.g. write 00122 for (0,0,1,2,2), then mentioned algorithm yields 

$$\C(\Sigma_2)=00122\ \uplus\ 00002\ \uplus\ 01002\ \uplus\ 01121\ \uplus\ 11000\ \uplus\ 11111.$$

\noindent
For instance, why is it that each $U\in (0,0,1,2,2)$ is $\Sigma_2$-closed? Because of $a,b\not\in U$, the first three implications in $\Sigma_2$ (vacuously) hold in $U$ since none of the three premises is contained in $U$. As to $\{d\}\ra\{c\}$, this implication holds in $U$ since $c\in U$.

\section{Semilattices}

By definition an (algebraic) {\it semilattice} is a commutative semigroup $Y$ such that $E(Y)=Y$.
In 5.1 we  show that $(Y,<_{\cal J})$ is a poset for each semilattice $Y$. This poset enjoys a crucial property (5.2). Conversely (5.3)  each semilattice can be {\it defined} as a certain poset. Subsection 5.4 indicates why semilattices are an important ingredient for arbitrary semigroups.

 \vspace{3mm}
{\bf 5.1}
Recall from 2.7.1 that for any semigroup $S$ and any $a,b\in S$ one says that $a$ is a {\it proper multiple} of $b$ (written $a<_{\cal J}b$) if there is $c\in S$ with $a=bc$. Futher we defined $a\le_{\cal J}b\ :\LRa\ (a<_{\cal J}b)\ or\ a=b)$. In Subsection 2.7 we met some semigroups for which the preorder $\le_{\cal J}$ actually is a partial order. It gets better:

\begin{itemize}
    \item[(25)]{\it  For each semilattice $Y$ even $<_{\cal J}$ is a partial order. Moreover for all  $e,f\in Y$ it holds that $e<_{\cal J} f \LRa e=ef$.}
\end{itemize}

{\it Proof of (25).} Since $<_{\cal J}$ is transitive, it remains to show reflexivity and antisymmetry.
Let us first show the additional claim that $e<_{\cal J} f \LRa e=ef$. The direction $\Leftarrow$ being obvious, assume that $e<_{\cal J}f$, i.e. that $e=fg$ for some $g\in S$. Then $ef=fgf=fg=e$.
 
 Now reflexivity is evident: $e<_{\cal J}e$ since $ee=e$. And so is antisymmetry: $(e<_{\cal J} f\ and\  f<_{\cal J} e)\Ra (e=ef\ and\ f=fe)\Ra e=f$.
$\square$

\vspace{3mm}
\noindent
Each semilattice $Y$ therefore comes with an associated partial order $(Y,<_{\cal J})$. In particular  for semilattices $<_{\cal J}$ is reflexive\footnote{As opposed to c.f. nilsemigroups and free c. semigroups $F_n$ in Section 2.}.

\vspace{5mm}
{\bf 5.2} In (26) we show that the partial order $<_{\cal J}$ has an exquisite property. A few preliminaries are in order. 

The elements $2$ and $6$ in the poset of Fig. 1C have two maximal {\it common lower bounds (clb)}, namely $4$ and $12$, but no largest clb $x$ (i.e. such that {\it all} clb's are $\le x$). A largest clb of $a,b$, if it  exists, is unique (why?). In this case it is called the {\it meet} of $a,b$ and written as $a\wedge b$.

\vspace{3mm}
{\bf 5.2.1} For instance, if $(P,\le)$ is any poset and $a,b\in P$ are comparable, say $a\le b$, then $a\wedge b=a$ exists. Consequently each {\it chain} (=totally ordered set) is a poset in which any two elements have a meet.

In a more general vein, for elements $c<d$ of a poset $(P,\le)$  the {\it interval} determined by $c,d$ is $[c,d]:=\{x\in P:\ c\le x\le d\}$.
If $P$ has a smallest element $\zero$ then $a,b\in P$ have a meet  iff $[\zero,a]\cap [\zero,b]$ has a largest element. 

By definition a {\it tree} is a poset $(Y,\le)$ with smallest element $\zero$ such that $[\zero,a]\cap [\zero,b]$ is a chain for all $a,b\in Y$; see Fig 3B. Hence  trees generalize chains in that 
 any two elements have a meet.

\vspace{3mm}
{\bf 5.2.2}  Each closure system $\C\s\Pow(X)$ yields
 a semilattice $(\C,\cap)$ since $\cap$ is associative, commutative, and idempotent. If $(\C,<_{\cal J})$ is the induced poset,
 then $A<_{\cal J} B\LRa A=A\cap B$ by (25), which amounts to $A\s B$. Pick any two $A,B$ in $(\C,\s)$. If $D\in\C$ is any clb of $A$ and $B$, then $D\s A$ and $X\s B$, hence $D\s A\cap B$. Since $A\cap B$ is itself a clb of $A$ and $B$, we find that $A\wedge B=A\cap B$.

\vspace{3mm}
{\bf 5.2.3} We just saw that in $(\C,<_{\cal J})$ any two elements have a meet. This generalizes:

\begin{itemize}
    \item[(26)] {\it If $Y$ is a semilattice then in the poset $(Y,<_{\cal J})$  any two elements $e,f$  have a meet $e\wedge f$. In fact $e\wedge f=ef$.}
\end{itemize}

{\it Proof of (26).} Let $x$ be any clb of $e,f$, so $x<_{\cal J}e$ and $x<_{\cal J} f$. From (25) we get $x=ex$ and $x=fx$, which yields $x=xx=(ex)(fx)=efx$, i.e. $x<_{\cal J} ef$. On the other hand $ef$ is itself a clb of $e,f$ (being a proper multiple of both). Therefore $ef$ is the largest clb of $e,f$. $\square$

\vspace{3mm}
{\bf 5.3} The above examples motivate the following definition.  Suppose  $(Y,\le)$ is {\it any} poset such that {\it all}
$a,b\in Y$  possess a meet $a\wedge b$. It makes a nice exercise (carried out in [Gr,p.9]) to show that the binary operation $\wedge$ is associative. Since idempotency is trivial, $(Y,\wedge)$ is a semilattice, which is called a {\it meet-semilattice}.

The bottom line is this. "Algebraic" semilattices (= c. sgr $Y$ with $E(Y)=Y$) and meet-semilattices (= certain  posets $(Y,\le)$) are  two sides of the same coin. We can take either view at our digression.

\vspace{3mm}
{\bf 5.3.1}
If $(Y,\le)$ is a meet-semilattice and $F\s Y$ is finite then the meet $\bigwedge F$ is well-defined (by associativity and induction). If $Y$ itself is finite, then $\zero:=\bigwedge Y$ is the smallest element of $(Y,\le)$, and simultaneously a zero of $(Y,\wedge)$. Similarly, {\it if} a largest element ${\bf I}$ of $(Y,\le)$ exists, then it is an identity of $(Y,\wedge)$.

\vspace{3mm}
{\bf 5.3.2} The smallest (nontrivial) meet-semilattice is the 2-element chain $\{\zero,{\bf I} \}$ (with $\zero< {\bf I}$). Direct products of meet-semilattices are meet-semilattices.
Interestingly each finite meet-semilattice $(Y,\wedge)$ occurs as a ssgr of $\{\zero,{\bf I}\}^t$. Specifically, let $Y=\{y_1,...,y_t\}$ and consider the direct product
$\{ \zero,{\bf I}\}\times\cdots\times\{\zero, {\bf I}\}=\{\zero, {\bf I}\}^t$, which is (4.3.1) isomorphic to $(\Pow(\{1,...,t\},\cap)$. It suffices to verify that $f(y_p):=\{i\le t:\ y_i\le y_p\}$ is an injective homomorphism from $Y$ to $\Pow(\{1,...,t\})$. Indeed,

\begin{itemize}
\item[] $f(y_p\wedge y_q)=\{i\le t:\ y_i\le y_p\wedge y_q\}\stackrel{why?}{=} \{i\le t:\ y_i\le y_p\ and\ y_i\le y_q\}$
\item[]
$=\{i\le t:\ y_i\le y_p\}\cap \{i\le t:\ y_i\le  y_q\}=f(y_p)\cap f(y_q).$
\end{itemize}

\vspace{3mm}

As to injectivity, from $y_p\neq y_q$ follows by antisymmetry that $y_p\not\le y_q$ or $y_q\not\le y_p$, say the latter. Then $q\not\in f(y_p)$ but $q\in f(y_q)$, and so $f(y_p)\neq f(y_q)$.

\vspace{3mm}

 Relabeling the meet semilattice $Y':=\C$ in Figure 2A with $1,2,...,8$ as indicated, the embedding of $Y'$ in $\{\zero,{\bf I}\}^8$ is spelled out in Fig. 2B.

\vspace{3mm}
{\bf 5.4} Let $(Y,\wedge)$ be any  meet-semilattice $Y$ and suppose there are disjoint semigroups $S_\alpha$ indexed by the elements of $Y$. We set $S:=\biguplus_{\alpha\in Y}S_\alpha$  and strive to make $S$ a semigroup in such a way that 
the semigroups $S_\alpha$ become subsemigroups of $S$.

Here come the details. For each $\alpha\ge\beta$ from $Y$ one needs a morphism $\sigma_{\alpha,\beta}: S_\alpha\to S_\beta$ such that 

\begin{itemize}
    \item[(27a)] $\sigma_{\alpha,\alpha}$ is the identity on $S_\alpha$;
    \item[(27b)] $\sigma_{\alpha,\beta}\circ\sigma_{\beta,\gamma}=\sigma_{\alpha,\gamma}$ for all $\alpha\ge\beta\ge\gamma$.
\end{itemize}

\noindent
On $S:=\biguplus_{\alpha\in Y}S_\alpha$ we define a binary operation $\ast$ as follows:

\begin{itemize}
    \item[(28)] If $a_\alpha\in S_\alpha$ and $b_\beta\in S_\beta$, then $a_\alpha\ast b_\beta:=(a_\alpha\sigma_{\alpha,\alpha\wedge\beta})\cdot(b_\beta\sigma_{\beta,\alpha\wedge\beta})\in S_{\alpha\wedge\beta}$.
\end{itemize}

\noindent
Here the dot $\cdot$ indicates multiplication within $S_{\alpha\wedge\beta}$. In view of (27a) it is clear  that  $a_\alpha\ast b_\alpha=a_\alpha\cdot b_\alpha$ for all $a_\alpha, b_\alpha\in S_\alpha$. Notice that for $\alpha>\beta$ we have $\alpha\wedge\beta=\beta$. Taking into account (27a) we conclude:

$$(28') If\ \alpha>\beta,\ then\ a_\alpha\ast b_\beta=(a_\alpha\sigma_{\alpha,\beta})\cdot b_\beta$$

\noindent
We leave it as an exercise (spelled out in [H]) to show that generally {\it both} $(a_\alpha\ast b_\beta)\ast c_\gamma$ and $a_\alpha\ast (b_\beta\ast c_\gamma)$ coincide with 

$$a_\alpha\sigma_{\alpha,\alpha\wedge\beta\wedge \gamma}\cdot b_\beta\sigma_{\beta,\alpha\wedge\beta\wedge \gamma}\cdot c_\gamma\sigma_{\gamma,\alpha\wedge\beta\wedge \gamma}$$

\noindent
for all $a_\alpha\in S_\alpha,\ b_\beta\in S_\beta,\ c_\gamma\in S_\gamma$.
Hence the operation $\ast$ is associative. One calls $S$ a {\it strong semilattice $Y$ of semigroups} $S_\alpha\ (\alpha\in Y)$.

Let us discuss two easy kinds of strong semilattices. The first (5.4.1) restricts  the semilattice $Y$, the second (5.4.2) restricts the  semigroups $S_\alpha$.

\vspace{4mm}
{\bf 5.4.1} In the first type the meet-semilattice  is a tree $Y=T$ with smallest element $\zero$. For each $\alpha\in T$ we let $\sigma_{\alpha,\alpha}$ be the identity on $S_\alpha$, and for each covering $\alpha\succ\beta$ in $T$ we choose an {\it arbitrary} morphism $\sigma_{\alpha,\beta}: S_{\alpha}\to S_{\beta}$.
Generally, when $\alpha>\delta$, there is a unique path $\alpha\succ\beta\succ\cdots\succ\gamma\succ\delta$, and accordingly we put $\sigma_{\alpha,\delta}:=\sigma_{\alpha,\beta}\circ\cdots\circ\sigma_{\gamma,\delta}$. It is then clear that (27a) and (27b) are satisfied.

\vspace{4mm}

{\bf 5.4.2} The second type does not restrict $Y$ but demands that all sgr $S_\alpha$ be {\it cyclic}, say of type $C_{m_\alpha,n_\alpha}$. As in 5.4.1 we first look at all coverings $\alpha\succ\beta$ and choose a morphism $\sigma_{\alpha,\beta}: C_{m_\alpha,n_\alpha}\to C_{m_\beta,n_\beta}$. Recall that Theorem 1 pinpoints the degree of freedom for $\sigma_{\alpha,\beta}$ in terms of exquisite integers $k_{\alpha, \beta}$ in
$\{1,2,...,m_\beta+n_\beta-1\}$. For general $\alpha>\delta$ we wish to repeat the definition of $\sigma_{\alpha,\delta}$ used for the tree case.

For the sake of notation we only consider a short sequence $\alpha\succ\beta\succ\delta$  of coverings in $Y$. Let $C_{m_\alpha,n_\alpha}=\langle a\rangle,\ C_{m_\beta,n_\beta}=\langle b\rangle,
\ C_{m_\delta,n_\delta}=\langle d\rangle,$ and let $k_1:=k_{\alpha,\beta},\ k_2:=k_{\beta,\delta}$ be coupled to $\sigma_{\alpha,\beta},\ \sigma_{\beta,\delta}$ respectively. 
Since $\sigma_{\alpha,\delta}:=\sigma_{\alpha,\beta}\circ \sigma_{\beta,\delta}$ is a morphism, it determines the exquisite integer $k:=k_{\alpha,\delta}$. One may suspect that $k=k_1k_2$. Let's see. By definition $a\sigma_{\alpha,\delta}=d^k$.
On the other hand

$$ a\sigma_{\alpha,\delta}=(a\sigma_{\alpha,\beta})\sigma_{\beta,\delta}
=(b^{k_1})\sigma_{\beta,\delta}=(b\sigma_{\beta,\delta})^{k_1}
=(d^{k_2})^{k_1}=d^{k_2k_1}$$

\noindent
Therefore $k\equiv k_1k_2\ (mod\ n_\delta)$ (but not necessarily $k=k_1k_2$).
Trouble is, there may also be $\gamma\neq\beta$ with $\alpha\succ\gamma\succ\delta$. This forces us to pick $\sigma_{\alpha,\gamma}$ and $\sigma_{\gamma,\delta}$ in such a way that the coupled integers $k_3:=k_{\alpha,\gamma}$ and $k_4:=k_{\gamma,\delta}$ also satisfy $k\equiv k_3k_4\ (mod\ n_\delta)$.

\vspace{3mm}
{\bf 5.4.3} To dig deeper, let us shrink $Y$ to the unique 4-element semilattice $Y:=\{\alpha,\beta,\gamma,\delta\}$ which is not a tree (and which has identity $\alpha$ and zero $\delta$). Further let

\begin{itemize}
\item[] $C_{m_\alpha,n_\alpha}=C_{2,4}=\langle a\rangle$
\item[] $C_{m_\beta,n_\beta}=C_{4,1}=\langle b\rangle$
\item[] $C_{m_\gamma,n_\gamma}=C_{1,6}=\langle c\rangle$
\item[] $C_{m_\delta,n_\delta}=C_{5,3}=\langle d\rangle$
\end{itemize}

\noindent
Arguing as in 2.3.1 one finds that

\begin{itemize}
    \item[] $Exq(\alpha,\beta):=Exq(2,4,4,1)=\{2,3,4\}$
     \item[] $Exq(\beta,\delta):=Exq(4,1,5,3)=\{3,6\}$
      \item[] $Exq(\alpha,\gamma):=Exq(2,4,1,6)=\{3,6\}$
       \item[] $Exq(\gamma,\delta):=Exq(1,6,5,3)=\{5,6,7\}$
        \item[] $Exq(\alpha,\delta):=Exq(2,4,5,3)=\{3,6\}$
\end{itemize}

\noindent
By  considering  not all morphisms $\sigma: \langle a\rangle\to \langle d\rangle$, but only those $\sigma$ that factor through $\langle b\rangle$,  the set  
 $Exq(\alpha,\delta)=\{3,6\}$ will shrink to some subset  $Exq(\alpha,\beta,\delta)$. By the above we know how to calculate the latter:

 $$Exq(\alpha,\beta,\delta)=\{2,3,4\}\cdot\{3,6\}=\{6,12,9,18,12,24\}=\{6\}$$

 \noindent
 The last $"="$ is due to the fact that all integers $\ge m_\delta$ are again reduced modulo $n_\delta$ to numbers lying in $\{m_\delta,m_\delta+1,...,m_\delta+n_\delta-1\}=\{5,6,7\}$. Similarly

 $$Exq(\alpha,\gamma,\delta)=\{3,6\}\cdot\{5,6,7\}=\{15,18,21,30,36,42\}=\{6\}$$

\vspace{1cm}
\includegraphics[scale=0.66]{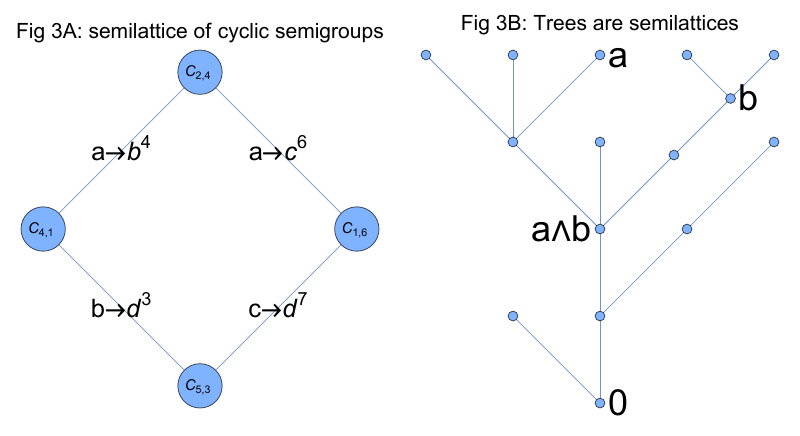}

\noindent
We see that whatever exquisite $k_1$ and $k_2$ we pick, composing the corresponding morphisms yields (incidently)  six times  the same result. Likewise for $k_3,k_4$. Consequently there are exactly 36 nonisomorphic "$Y$-frame" semilattices $S$ of $\{C_{2,4},\ C_{4,1},\ C_{1,6},\ C_{5,3}\}$. One of these is shown in Figure 3A. We will resume this matter in Section 9.

\vspace{3mm}
{\bf 5.5} Dually to the meet,  elements $a,b$ of a poset $P$ have a {\it join} $a\vee b$ if the latter is the least common upper bound of $a,b$. One calls $P$ a {\it join-semilattice} if any two elements of $P$ possess a join.
Then, akin to $(P,\wedge)$, also $(P,\vee)$ is a semilattice in the algebraic sense of 5.1. Note that $x\le y\LRa x\vee y=y$ in  each join-semilattice $P$. Further, if $P$ has a smallest element $\zero$, then the poset $P\setminus\{\zero\}$ remains a join-semilattice.

A poset which is both a meet-semilattice and a join-semilattice, is called a {\it lattice}.
Often lattices arise as follows. Let $\C\s\Pow(X)$ be a closure system. Recall that the poset $(\C,\s)$ is a meet-semilattice with $A\wedge B=A\cap B$. One checks that $(\C,\s)$ also is a join-semilattice with $A\vee B=cl_{\C}(A\cup B)$.

\section{Generators and relations}

An elegant way  to define c.f. semigroups $S$ is by generators and relations. Some of the technicalities (local confluence in digraphs, congruence relations) being deferred to Section 9, in Section 6 we concentrate on "how to do it?" rather than "why does it work?".
Subsections 6.2 to 6.5 are devoted to represent the elements of $S$ by "normal forms". Interestingly, and little known, in the case of semilattices $S$ one can dispense with normal forms and moreover get $S$ in a compressed format.

\vspace{3mm}
{\bf 6.1} In 2.2 we saw that each finite sgr $\langle a\rangle$ is isomorphic to some sgr $C_{m,n}$.
What about the converse? Given say $m=3,\ n=4$, {\it is there} a semigroup $C_{3,4}$? One may be tempted to answer as follows. 

\vspace{3mm}
{\bf 6.1.1} Yes,  take any symbols $a_1,a_2,...,a_6$ and define  $a_i\ast a_j:=a_{i+j}$ if $i+j\le 6$; otherwise put
$a_i\ast a_j:=a_k$ where $k$ is the unique number in $\{3,4,5,6\}$ satisfying $k\equiv i+j\ (mod\ 4)$. Trouble is, proving the associativity of $\ast$ is awkward. 

\vspace{3mm}
{\bf 6.1.2} Here comes a better way. Let $(T_7,\circ)$ be the semigroup  of all selfmaps $a:\{1,..,7\}\to\{1,..,7\}$ under composition. Let us exhibit some ssgr of $T_7$ which is of type $C_{3,4}$. Namely, putting

$$a:={\begin{pmatrix}
    1 & 2& 3& 4& 5& 6& 7\\
    2&3&4&5& 6& 3& 1
\end{pmatrix}}$$

\noindent
one checks brute-force that all maps $a,a^2,...,a^6$ are distinct and that $a^7=a^3$. Therefore the ssgr $\langle a\rangle$ of $T_7$ is of type $C_{3,4}$.

This  beats  6.1.1, but isn't perfect either. While one may guess which $a\in T_{m+n}$ to pick for general $C_{m,n}$, it will not work when more than one "relation" $a^{m+n}=a^m$ needs to be satisfied.
\vspace{3mm} 

{\bf 6.2} Let $S=\langle a,b,c\rangle$ be a commutative semigroup whose generators satisfy the relations $a=a^2,\ b^3=ab^2,\ bc=c^2$. Recall from 2.1 that each element of $S$ can be written (possibly in several ways)  as $a^ib^jc^k$. In order to find an upper bound for $|S|$ we {\it direct} each relation in the {\it presentation} 

$$\{a=a^2,\ b^3=ab^2,\ bc=c^2 \}$$

\noindent
from military-larger to military-smaller (see 2.8):

$$a^2\ra a,\ b^3\ra ab^2,\ c^2\ra bc$$

\noindent
Below we list all words $a^ib^jc^k$ (= members of $F_3$ by 2.8.1) in military order as well, 
thus starting with $a<_M b<_M c<_M a^2<_M<\cdots$. The brackets show how some words can be "reduced", using the directed relations above, to previously listed words:

\begin{itemize}
    \item[i.] $a,\ b,\ c$
    \item[ii.] $(a^2=a), ab,\ ac,\ b^2,\ bc, (c^2=bc) $
    \item[iii] $(a^3=a),(a^2b=ab),(a^2c=a),ab^2,abc,(ac^2=abc),(b^3=ab^2),b^2c,\\(bc^2=bbc),(c^3=c\cdot bc=bc^2)$
    \item[iv:] $(a^4=a^2a^2=aa),(a^3b=ab),(a^3c=ac),(a^2xy=axy), \\
    (ab^3=a\cdot ab^2=ab^2),\quad ab^2c,\quad (abc^2=ab\cdot bc), (ac^3=abc^2),\\
    (b^4=b\cdot ab^2=ab^3),(b^3c=ab^2\cdot c), 
    (b^2c^2=b^2bc),\\
    (bc^3=b\cdot b^2c),(c^4=c^2c^2=bcbc)$
    \item[v:] $(i\ge 2: a^ib^jc^k=ab^jc^k),(j\ge 3: ab^jc^k=a^2b^{j-1}c^k),(ab^2c^2=ab^3c),\\(abc^3=ab^3c),(ac^4=ab^2c^2),(j\ge 3: b^jc^k=ab^{j-1}c^k),(b^2c^3=b^3c^2),\\
    (bc^4=b^3c^2),(c^5=b^2c^3)$
\end{itemize}

 \noindent
 Omitting all bracketed expressions we conclude that each element of $S$ can be reduced to one of these {\it normal forms}:

$$(29)\quad a,b,c;\hspace{3mm} ab,ac,b^2,bc;\hspace{3mm} ab^2,abc,b^2c;\hspace{3mm} ab^2c$$

\noindent
In particular $S$ has at most 11 elements. The  smallest qualifying $S$ is $S=\{e\}$ which satisfies all relations: $a^2=a=b^3=ab^2=c^2=bc=e$. 

\vspace{2mm}
{\bf 6.2.1}
Although several semigroups may satisfy a given set of relations, it turns out  that one of these semigroups, written  as

$$RFCS(a,b,c:\ a^2\ra a,\ b^3\ra ab^2,\ c^2\ra bc)\hspace{3mm} (=:RF_1),$$

\noindent
is {\bf the largest}\footnote{A  proof is provided in Section 10. Notice  $RFCS(..)$ can be infinite, e.g. $RFCS(a,b,c:\es)$ (no relations) is isomorphic to $F_3$ from 2.8. However, all upcoming sgr $RFCS(..)$ are tuned to be finite because the structure theory of Section 8 only applies in the finite case.} insofar that the others are epimorphic images of $RFCS(...)$. The acronym RFCS stands for {\it relatively free  commutative semigroup} (wrt the postulated relations). 

The above entails that in particular each  element of $RF_1$ can be written as some normal form listed in (29). Could it be (as it happens for $S=\{e\}$) that different normal forms yield the same element of $RF_1$? We will find out soon.

\vspace{3mm}
{\bf 6.3} Let us calculate the normal forms of the elements of

$$RF_2:=RFCS(a,b:\ b^4\xrightarrow{1} b^2,\  a^3\xrightarrow{2} b^2,\ a^4\xrightarrow{3} a).$$

\noindent
Here and henceforth the relations in the presentation $\{b^4=b^2,\  a^3= b^2,\ a^4=a \}$ are already directed from military-larger to military-smaller.

\vspace{3mm}
{\bf 6.3.1} 
Let us  list the normal forms of $RF_2$ as we did in 6.2:

\begin{itemize}
    \item[i.] $a,b$
    \item[ii.] $a^2, ab, b^2$
    \item[iii.] $(a^3=b^2),a^2b, ab^2, b^3 $
    \item[iv.] $(a^4=a), (a^3b=b^2\cdot b), a^2b^2, ab^3, (b^4=b^2)$
    \item[v.] $(a^5=a^2), (a^4b=ab),(a^3b^2=b^2\cdot b^2=b^2),a^2b^3, (ab^4=ab^2),(b^5=b^3)$
    \item[vi.] $(i\ge 4:\ a^ib^{6-i}=a^{i-3}b^{6-i}),(a^3b^3=b^5),(i\le 2:\ a^ib^{6-i}=a^ib^{4-i})$
\end{itemize}

\noindent
It follows that every element of $RF_2$ can be written as one of these normal forms:

$$(30)\quad a, b;\hspace{3mm} a^2,ab,b^2;\hspace{3mm} a^2b,ab^2,b^3;\hspace{3mm}
a^2b^2, ab^3;\hspace{3mm} a^2b^3$$

\noindent
What was dooming at the end of  6.2 takes place here. Since $ab^2=a\cdot a^3=a$, different normal forms describe the same element of $RF_2$! 

\vspace{3mm}
{\bf 6.4} This leads us to the crucial  issue of "local confluence". To begin with, if

$$ab^2\xrightarrow{4} a$$

\noindent
gets added as a new relation then the  above problem is settled since $ab^2$ ceases to be a normal form. But perhaps {\it other} problems remain. Figure 4A clarifies that our problem was a kind of incompatibility of the 2nd relation $\rho_2$ with the 3rd relation $\rho_3$, and that by adding the new relation $\rho_4$ one achieves "local confluence".

\vspace{3mm}
Specifically one calls a set of  relations {\it locally confluent} if the following takes place for each pair of distinct relations  $\rho_i:=(v_i\xrightarrow{i} w_i)$ and $\rho_j:= (v_j\xrightarrow{j} w_j)$. Let $v:=lcm(v_i,v_j)$ be the least common multiple\footnote{Here $v_i,v_k,w_i,w_k$ are words (=members of $F_n$ from 2.8) over some "alphabet" $\{a,b,..\}$. The "least common multiple" has an obvious meaning; say $lcm(a^2bc^4,a^5c^3)=a^5bc^4$.} of the "premises" $v_i$ and $v_j$. It is evident that $\rho_i$ reduces $v$ to some word $w_i'$, and $\rho_j$ reduces $v$ to some word $w_j'$. The definition of "locally confluent" requires that one can  reduce both $w_i'$ and $w_j'$ to a common word $w$. One calls $(w_i',w_j')$  a {\it critical pair}; see also 10.4.3.

\vspace{3mm}
In our example with 4 relations we must therefore check ${\binom{4}{2}}=6$ pairs of relations. That this test is successful for 4 out of 6 pairs is illustrated in Figures 4A,4B,4C,4D. The pairs $\rho_1,\rho_2$ and $\rho_1,\rho_3$ fit the hat of what happens when the premises of $\rho_i,\rho_j$ are {\it disjoint} (i.e. without common letters). In this case $lcm(v_i,v_j)=v_iv_j$, and this guarantees local confluence (see Fig 4E). 

\vspace{1cm}
\includegraphics[scale=0.64]{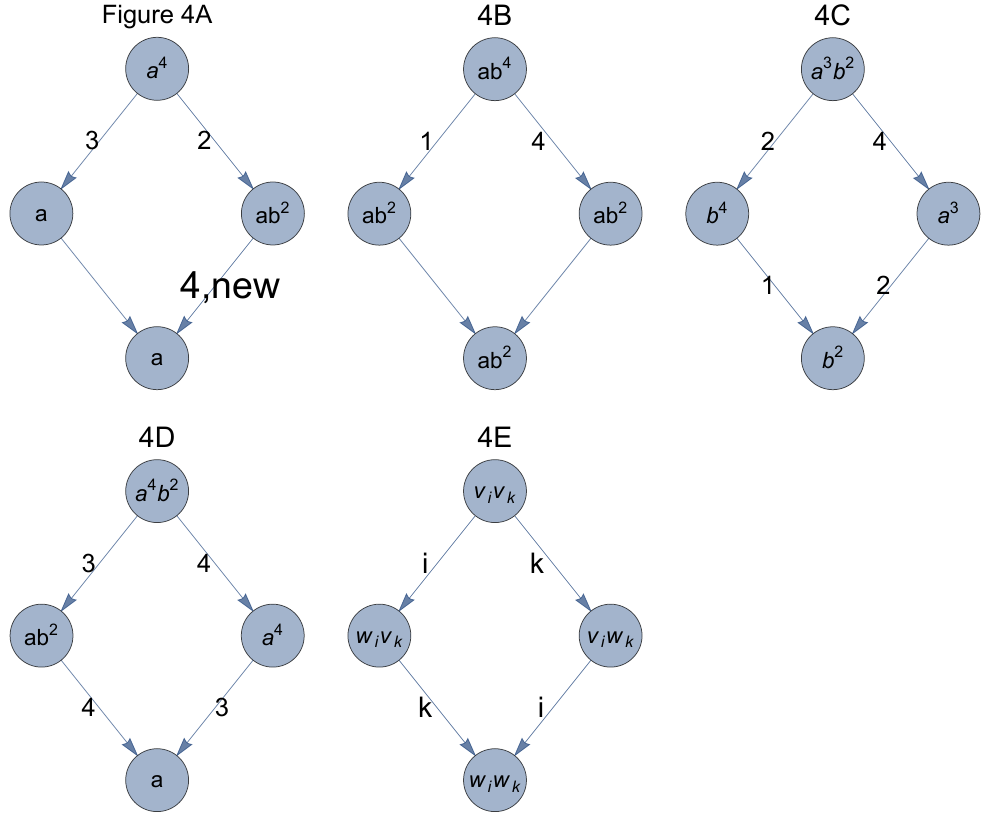}

\vspace{5mm}
{\bf Theorem 8: }{\it If the relations defining $RF:=RFCS(...)$ are locally\\ confluent, then the normal forms bijectively match the elements of $RF$.}
\vspace{5mm}

The proof is deferred to Section 10. As to $RF_2$, the list of normal forms with respect to our enlarged presentation $\{\rho_1,\rho_2,\rho_3,\rho_4\}$ is obtained by pruning  list (30) with $\rho_4$. We hence obtain

$$(31)\quad RF_2=\{a, b;\hspace{3mm} a^2,ab,b^2;\hspace{3mm} a^2b,b^3\}$$

\noindent
Because we checked local confluence, it follows from Theorem 8 that $|RF_2|=7$.

\vspace{3mm}
{\bf 6.4.1} Because the relations defining $RF_1$ had pairwise disjoint premises, and this is sufficient for local confluence, it also follows from Thm.8 that  $|RF_1|=11$, the normal forms being given in (29).

If there only is one relation, then local confluence is even more trivial.
In particular, let the single relation be $a^{m+n}=a^m$. It then follows that for any two integers $m,n\ge 1$ there is a semigroup of type $C_{m,n}$. This blows away the problems of 6.1.

\vspace{3mm}
{\bf 6.5} Whatever the fine structure of $RF_3$ below, its element $\zero$ is indeed a zero of $RF_3$, and  $|RF_3|\le 3\cdot 4\cdot 5=60$ (why?).

\begin{itemize}
\item[$RF_3$] $:=RFCS(a,b,c,\zero:\  \zero\zero\ra \zero,\
a\zero\ra \zero,\ b\zero\ra \zero,\ c\zero\ra \zero,$
\item[] $\ a^3\xrightarrow{1} \zero,\ 
b^4\xrightarrow{2} \zero,\ c^5\xrightarrow{3} \zero,\  
 a^2b^2c^3\xrightarrow{4} \zero,\ ac^4\xrightarrow{5} \zero,\  b^3c^2\xrightarrow{6} \zero,\  ab^3\xrightarrow{7} \zero)$. 
\end{itemize}

\noindent
Let us argue that the given presentation (and all of this type) is locally confluent. 
Take any two of the given $4+7$ relations, say $v_i\to\zero$ and $v_j\to\zero$. Taking $v:=lcm(v_i,v_j)$ yields $v\to w_i\zero$ and $v\to w_j\zero$. In view of the first four relations both $w_i\zero$ and $w_j\zero$ can step by step (and in many ways) be reduced to $\zero$.

We thus know from Thm. 8 that the unique normal forms bijectively match the elements of $RF_3$. What are the normal forms? If $w\in F_3=\langle a,b,c\rangle$ is such that (component-wise)

$$ a^3\le w\quad or\quad b^4\le w\quad or\ .....\ or\quad  b^3c^2\le w\quad or\quad ab^3\le w,$$

\noindent i.e. $w$ belongs to the ideal $I_1\s F_3$ in 2.10.1, then\footnote{To fix ideas, consider $w=a^3b^2c^4$, which (e.g.) is $\ge a^2b^2c^3$. Applying $a^2b^2c^3\xrightarrow{4} \zero$ to $w$ yields $ac\zero$, which further reduces to $c\zero$ and then to $\zero$.} $w$ reduces to the normal form $\zero$. On the other hand, if $w\in F_3\setminus I_1$, then no relation is applicable, i.e. $w$ is in normal form already. It follows that $RF_3$ is isomorphic to the Rees quotient $F_3/I_1$. Recall from 2.10.1 that $F_3\setminus I_1$ can be rendered in a compressed format.

\vspace{2mm}
{\bf 6.5.1} What about an arbitrary semigroup

$$RF=RF(a,b,..,c: v_1\to w_1,...,v_t\to w_t)?$$

\noindent
The set $NF$ of all normal forms consists of all words $w=a^ib^j\cdots c^k$ such that $(\forall 1\le s\le t)\ w\not\ge v_s$. As above $NF$ can be compressed, and so $|NF|$ can be calculated fast. This implies the handy  upper bound $|RF|\le |NF|$. If the set of relations is locally confluent, one even has $|RF|= |NF|$. Unfortunately, different from the Rees quotient scenario above, the {\it structure} of the sgr $RF$ remains elusive (until Section 8).


\vspace{3mm}
{\bf 6.6} Consider the semigroup $RF_4$ below. Because of the first four relations
each element of $RF_4$ can be written as
  $a^ib^jc^kd^m$ with $i,j,k,m\in\{0,1\}$.

  $$RF_4:=RFCS(a,b,c,d:\ a^2\ra a,\cdots,d^2\ra d,\ abc\xrightarrow{1} ab,\ ac\xrightarrow{2} c,\ bd\xrightarrow{3} b)$$
  
\noindent
For instance  $a^1b^1c^0d^1$ means $abd$, and  $(abd)^2=abdabd=a^2b^2d^2=abd$. Evidently all other elements of $RF_4$
 (recall $a^0b^0c^0d^0\not \in RF_4$) are idempotent as well, and so $RF_4$ is a semilattice.

\vspace{3mm}
{\bf 6.6.1}
As to local confluence, it fails for the relations $\rho_1,\rho_2$ and triggers the new relation $\rho_4:=(bc\xrightarrow{4} ab)$. Notice that $\rho_4$ and $a^2\ra a$ yield $\rho_1$ in the sense that from  $bc=ab$ and $a^2=a$ follows $abc=a\cdot bc=a\cdot ab=ab$. Dropping superfluous relations (such as $\rho_1$) shortens the calculation of the normal forms. Having checked (do it) that the presentation $\{\rho_2,\rho_3,\rho_4\}$ is locally confluent, we can be sure that the respective normal forms  represent, without repetition, the elements of $RF_4$. Here they come:

$$a,b,c,d,ab,(ac=c),ad,(bc=ab),(bd=b),cd,$$
$$(abc=bc),(abd=ab),(acd=cd),(bcd=bc)$$

\vspace{3mm}
{\bf 6.6.2}
When dealing with relatively free {\it semilattices (SL)} a trimmed notation (dropping all relations $x^2\ra x$) is preferable. Thus 

$$(32)\quad RF_4=RFSL(a,b,c,d:\ ac\xrightarrow{2} c,\ bd\xrightarrow{3} b,
\ bc\xrightarrow{4} ab)=\{a,b,c,d,ab,ad,cd\}$$

\vspace{3mm}
{\bf 6.7} Here come  the benefits of viewing type $RFSL(..)$ semilattices as join-semilattices. Let us begin by rewriting the presentation in (32) as join-semilattice presentation

 $$X:=\{a\vee c=c,\ b\vee d=b,\ b\vee c=a\vee b\}.$$
 
 \noindent
 This is equivalent (see 5.5) to 
 
 $$X':=\{c\ge a,\ b\ge d,\  b\vee c\ge a\vee b,\  a\vee b\ge b\vee c \}.$$
 
 \noindent
 Since visually $\ge$ resembles $\to$, it comes easy to move from $X'$ to the set of implications 

$$\Sigma_1':=\{\{c\}\to\{a\},\ \{b\}\to\{d\},\ \{b,c\}\to\{a,b\},\ \{a,b\}\to\{b,c\}\}$$

Since $\Sigma_1'$ is equivalent\footnote{First note that $\{c\}\to\{a\}$ and $\{b\}\to\{d\}$ are present in both $\Sigma_1$ and $\Sigma_1'$. The implication $\{a,b\}\to\{c\}$ in $\Sigma_1$ "follows from" $\{a,b\}\to\{b,c\}$  in $\Sigma_1'$. Conversely all implications in $\Sigma_1'$  follow from implications in $\Sigma_1$ (try). For a precise definition of "follows from" see [W3].} 
to $\Sigma_1$ in 4.3, it follows that $\C:=\C(\Sigma_1')$ equals $\C(\Sigma_1)$, which was shown in Fig.2A.
Look at $(\C\setminus\{\es\},\vee)$ in Fig.5, which is obtained by cutting $\{\es\}$ from Fig.2A and tilting it. Each node in Fig.5 is labeled by an element of $(\C\setminus\{\es\},\vee)$ and a corresponding element (in normal form) in $RF_4$. 

This correspondence is a semilattice isomorphism. For instance $a\cdot b=ab$ in $RF_5$ matches $\{a\}\vee\{b,d\}=\{a,b,c,d\}$ in $\C\setminus\{\es\}$. Further,
$a\cdot cd\stackrel{2}{=} cd$ matches $\{a\}\vee\{a,c,d\}=\{a,c,d\}$, and $d\cdot cd=cd$  matches $\{d\}\vee\{a,c,d\}=\{a,c,d\}$, and $b\cdot c\stackrel{4}{=} ab$ matches $\{b,d\}\vee\{a,c\}=\{a,b,c,d\}$.

\vspace{3mm}
{\bf 6.7.1}
All of this generalizes as follows [W1,Thm.5]. If $\Sigma$ is the family of implications derived from a join-semilattice presentation $X$, then $RFSL(a,b,..:\ X)$ is isomorphic to $(\C(\Sigma)\setminus\{\es\}),\vee)$. It is further shown\footnote{This is done in a direct way, i.e. without using Thue congruences as in 10.3.} in [W1] that every join-semilattice satisfying the relations in  $X$ is an epimorphic image of $(\C(\Sigma)\setminus\{\es\}),\vee)$.

It is an exercise to show that for each finite semilattice $Y$ and each $x\in Y$ there is a {\it largest} subset  $T\s S$ such that the product of all elements in $T$ is $x$. Figure 5 (and this generalizes to arbitrary $S$ of type $RFSL(..)$) gives these $T$'s explicitely. Thus if $x=b$, then $T=\{b,d\}$. If $x=ad$, then $T=\{a,d\}$. If $x=ab$, then $T=\{a,b,c,d\}$.

\vspace{1cm}
\includegraphics[scale=0.7]{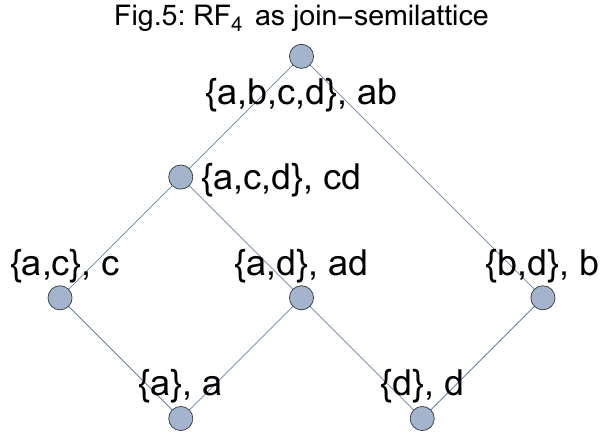}

\vspace{5mm}
{\bf 6.8} In order to showcase the advantages of the new method, let us tackle

$$ RF_5:=RFSL(a,b,c,d,e:\ a\vee b=a,\ b\vee c\vee e=b\vee c,\ a\vee b\vee d\vee e=a\vee e,\ c\vee d=d).$$

\noindent
If we were to apply the old method, 
we would have to add several\footnote{While tedious, this is certainly a good exercise to rub in the concept of local confluence.} new relations in order to achieve local confluence. 
In contrast, the new method ignores military order and local confluence. All that matters is to translate the join-semilattice relations into implications (recall, $\ge$ becomes $\to$):
 
$$\Sigma_2':=\Bigl\{ \{a\}\ra \{b\},\ \{b,c\}\ra \{e\},\ \{a,e\}\ra \{b,d\},\ \{d\}\ra \{c\}\Bigr\}$$. 

\noindent
Since $\Sigma_2'$ happens (lucky us) to be $\Sigma_2$ from 4.3.1, we find that

$$\C(\Sigma_2')=00122\ \uplus\ 00002\ \uplus\ 01002\ \uplus\ 01121\ \uplus\ 11000\ \uplus\ 11111.$$

\noindent
In particular, one reads off that $|RF_5|=12$.

\section{Archimedean semigroups}

In 7.1 we define Archimedean semigroups and investigate their kernels, then show that  c.f. semigroups have plenty Archimedean subsemigroups (7.2), then investigate the behaviour of direct products (7.3). Of special interest (7.4) are direct products of type $\Z_{p_1^{\epsilon_1}}\times\cdots\times \Z_{p_t^{\epsilon_t}}$, because  this leads (7.5) to the fine structure of $\Z_n=(\Z_n,\odot)$.

\vspace{3mm}
 {\bf 7.1} The c.f. semigroup $A$ is {\it Archimedean} if it has exactly one idempotent, i.e. $|E(A)|=1$. Hence Archimedean sgr and semilattices are two extreme types of c.f. semigroups; those with the fewest and those with the most idempotents.
In turn, c.f. groups $G$ and c.f. nilsemigroups $N$ are extreme types of Archimedean semigroups in that $E(G)=\{\one\}$ and $E(N)=\{\zero\}$.

It holds that (why?) that the Rees quotient $C_{m,n}/K(C_{m,n})$ is isomorphic to $C_{m,1} $. More generally:

\begin{itemize}
\item[(33)] {\it If $A$ is Archimedean,  then $A/K(A)$ is a nilsemigroup.}
\end{itemize}

{\it Proof.} It suffices to show that for each $a\in A$ some power $a^k$ is in $K(A)$. Indeed, the unique idempotent $e$ of $A$ sits in $K(A)$ since $K(A)$ is a subgroup (Theorem 6). We know from (5) that $a^k=e\in K(A)$ for some $k\ge 1$.

\vspace{5mm}
    {\bf 7.1.1} Take 
$A:=\{21,63,105,147,189,231,273,315,357,399,441,483\}.$
One checks brute-force that $A$ is a ssgr of $\Z_{504}$, which has a unique idempotent 441. Hence $A$ is Archimedean, but what is the fine structure of $K(A)$ and $A/K(A)$? Patience.

\vspace{5mm}
{\bf 7.2} Many c.f. semigroups $S$ are teeming with Archimedean ssgr because of the following fact. For any fixed $e\in E(S)$ let $A_e$ be the set of all $x\in S$ a power of which equals $e$, thus

$$A_e:=\{x\in S:\ e\in \langle x\rangle\}.$$

\noindent
Since each ssgr $\langle x\rangle$ of $S$ contains exactly one idempotent, $S$ is the {\it disjoint} union of the sets $A_e\ (e\in E(S)$.
It gets better: 

$$ x,y\in A_e\Ra  (\exists k,\ell)\ (x^k=y^\ell=e)\Ra(xy)^{k\ell}=(x^k)^\ell (y^\ell)^k=ee=e\Ra xy\in A_e,$$

\noindent
and so $A_e$ is a ssgr of $S$, which of course is Archimedean.
The cases $e=\one$ and $e=\zero$ are easily handled:

\begin{itemize}
    \item[(34)] If the c.f. sgr $S$ has an identity ${\bf 1}$ then $A_{\bf 1}=S^{inv}$.\\
    If $S$ has a zero $\zero$ then $A_{\bf \zero}$ is a nil ideal of $S$.
\end{itemize}

\vspace{2mm}
{\bf 7.3} Here we ask: To what extent do the operations "taking idempotents" or "taking kernels" or "taking Archimedean components" carry over to direct products? To begin with it holds that

\begin{itemize}
    \item[(35.1)] $E(S_1\times\cdots \times S_t)=E(S_1)\times\cdots\times E(S_t)$
    \item[(35.2)] $K(S_1\times\cdots \times S_t)=K(S_1)\times\cdots\times K(S_t)$
\end{itemize}

Property (35.1) is evident. As to (35.2),  $K:=K(S_1)\times\cdots\times K(S_t)$  is an ideal of $S_1\times\cdots \times S_t$, and so $K\supseteq K(S_1\times\cdots \times S_t)$. On the other hand, $K$ is a group by Theorem 6 and therefore cannot properly contain another ideal of $S$. Hence $K=K(S_1\times\cdots \times S_t)$. $\square$

\vspace{3mm}

{\bf 7.3.1} So much about $E(S)$ and $K(S)$. Let us now investigate the Archimedean components of the direct product
$S:=S_1\times\cdots \times S_t$. It's only for ease of notation that we stick to $t=2$ and switch from $S_1\times S_2$ to $S'\times S^{"}$. 
Suppose that $S'$ has the Archimedean components $A'_e\ (e\in E(S')$ and $S^{"}$ has $A_f^{"}\ (f\in E(S^{"}))$. If $A_{(e,f)}$ is the ($k$-element) Arch. component of $(e,f)\in E(S'\times S^{"})$ then for all $(x,y)\in S'\times S^{"}$ we argue similarly to 7.2: 

$$(x,y)\in A_{(e,f)}\LRa (x,y)^{k\ell}=(e,f)\LRa (x^k=e,\ y^\ell=f)\LRa (x\in A'_e,\  y\in A_f^{"}).$$

Therefore $A_{(e,f)}=A'_e\times A_f^{"}$. In particular, suppose that $S'$ has $\alpha$ Arch. components whose cardinalities sum up as $m_1+\cdots+m_\alpha=|S'|$, and likewise  $S^{"}$ has $\beta$ Arch. components whose cardinalities sum up as $n_1+\cdots+n_\beta=|S^{"}|$. Then $S'\times S^{"}$ has $\alpha\beta$ many Archimedean components $A'_e\times A_f^{"}$  whose cardinalities sum up as
$m_1n_1+m_1n_2+\cdots+m_\alpha n_\beta=|S'\times S^{"}|$.

\vspace{3mm}
{\bf 7.3.2} If $S_1,...,S_t$ are f. cyclic semigroups, then $S:=S_1\times\cdots\times S_t$ is Archimedean by (35.1). In view of (35.2) and 3.4.2 and 2.2.2 a necessary condition for  c.f. semigroups $S$ to have a $t$-element trimmed generating set is  this: $S$ is Archimedean and the Abelian group $K(S)$ is a product of $t$ cyclic groups (some of which may be trivial). 
\vspace{3mm}

{\bf Open Question 1: As compared to finite Abelian groups (recall $t_{min},\ t_{max}$), if an Archimedean semigroup happens to be a direct product of cyclic semigroups, in how many ways is this possible?}

\vspace{3mm}
{\bf 7.4} Of special interest are direct products of type
$\Z_{p_1^{\gamma_1}}\times\cdots\times  Z_{p_t^{\gamma_t}}\ (t\ge 1)$. By the above it suffices to find $Y:=E(\Z_{p^\gamma})$ and the structure of the corresponding Arch. components $A_e\ (e\in Y)$. 

So suppose $e\in \Z_{p^\gamma}=\{\zero,\one,2,\ldots,p^\gamma-1\}$ is idempotent. From $e^2\equiv e\ (mod\ p^\gamma)$ follows $e^2-e=e(e-1)\equiv 0\ (mod\ p^\gamma)$, hence $p^\gamma$ divides $e(e-1)$. Since $e$ and $e-1$ are coprime, either $p^\gamma$ divides $e$, or $p^\gamma$ divides $e-1$. Therefore we conclude:

\begin{itemize}
    \item[(36)]  $E(\Z_{p^\gamma_1})=\{\zero,\one\}$,  so $\Z_{p^\gamma}=A_\one\uplus A_\zero$, where 
    \begin{itemize}
    \item[] $A_\zero=\{p,2p,3p,...,p^{\gamma-1}p\ (=\zero) \}$ and 
    \item[] $A_\one=\Z_{p^\gamma}^{inv}=\{\one,2,3,\ldots,p-1,p+1,\ldots,2p-1,2p+1,\ldots, p^\gamma-1\}$.
      \end{itemize}
\end{itemize}

{\bf 7.4.1} To fix ideas, taking the $t=3$  prime powers $7^1,2^3,3^2$ we have  

$$\Z_7=\{0',1',...,6'\},\ \Z_8=\{0",1",...,7"\},\ \Z_9=\{0^*,1^*,...,8^*\},$$

\noindent
and so the set of idempotents of 
$\Z_7\times\Z_8\times\Z_9$ \\ is $Y_8:=E(\Z_7)\times E(\Z_8)\times E(\Z_9)=
\{0',1'\}\times\{0",1"\}\times\{0^*,1^*\}$, so $|Y_8|=8$. To spell it out:

$$(37)\quad Y_8=\Bigl\{(0',0",0^*),(0',0",1^*),(0',1",0^*),...,(1',1",0^*),(1',1",1^*)\Bigr\}.$$

\noindent
By (34) and (36) the Arch. components of $\Z_7\times\Z_8\times\Z_9$ are  direct products of groups and nilsemigroups. For instance $A_{(1',0",0^*)}=A_{1'}\times A_{0"}\times A_{0*}$ is a direct product of a group and two nilsemigroups, hence   of cardinality $6\cdot 4\cdot 3={\bf 72}$.
In view of $7\cdot 8\cdot 9=(1+6)(4+4)(3+6)$, our sgr $A_{(1',0",0^*)}$ is  one of 8 Archimedean components indexed by $Y_8$, which (in the order matching (37)) have cardinalities $12, 24,12,24, {\bf 72},144,72,144$.

\vspace{3mm}
{\bf 7.5} From 3.5 we know that $n=p_1^{\gamma_1}\cdots p_t^{\gamma_t}$ implies\\ $\Z_n \simeq \Z_{p_1^{\gamma_1}}\times\cdots\times  Z_{p_t^{\gamma_t}}$ (as monoids), and we
know how to translate back and forth. For instance, the Archimedean component

$$ A_{(0',1",0^*)}=A_{0'}\times A_{1"}\times A_{0*}=\{0'\}\times\{1",3",5" ,7"\}\times\{0^*,3^*,6^*\}$$

\noindent
of $\Z_7\times\Z_8\times\Z_9$ translates to the subset\\
$\{21,63,105,147,189,231,273,315,357,399,441,483\}$  of $\Z_{504}$, which popped up in 7.1.1, and which we now reckognize as the Archimedean component $A_{441}$ of $\Z_{504}$.
By (35.2) we have 

$$K(A_{(0',1",0^*)})=K(A_{0'})\times K(A_{1"})\times K(A_{0*})=\{0'\}\times\{1",3",5",7"\}\times\{0^*\},$$

\noindent
which translates to $K(A_{441})=\{63,189,315,441  \}\simeq C_2\times C_2$.

\vspace{2mm}
As another example, the last (wrt (37)) Archimedean component 
$A_{(1',1",1^*)}=(\Z_7\times\Z_8\times\Z_9)^{inv}$ translates to $\Z_{504}^{inv}$, which in 3.4.2 we found to be of type $C_6\times C_2\times C_2\times C_6$.

\vspace{2mm}
{\bf 7.6} We have come to understand the multiplicative semigroup $(\Z_n,\odot)$ of the particular rings $(\Z_n,+,\odot)$. 

\vspace{3mm}
{\bf Open Question 2: What is the state of affairs  for $(R,\cdot)$, where $(R,+,\cdot)$ is an  arbitrary commutative finite  ring? (Many other things are known about such rings [BF].)}

\section{The Structure Theorem}

The Structure Theorem states that each c.f. sgr $S$ is a semilattice of Archimedean semigroups. This inspires a five step recipe (I),..,(V), that achieves the following. Given the Cayley table of any c.f. sgr, its underlying semilattice and Archimedean components (including their fine structure) are identified. As to "given the Cayley table", apart from $(\Z_{n},\odot)$,  all our sgr are of type $S=RFCS(..)$, and so the Cayley table can be gleaned from the normal forms. Step (V) is the most difficult one and is dealt with in the final Subsections 8.6 and 8.7.

\vspace{4mm}
{\bf 8.1} Let $S$ be a strong semilattice $Y$ of semigroups $S_\alpha\ (\alpha\in Y)$. By definition of $S$ and (28) it holds that:

\begin{itemize}
    \item[(38)] There is a meet semilattice $Y$ such that $S$ is the disjoint union of subsemigroups $S_\alpha$ indexed by the elements of $Y$. Furthermore $S_\alpha S_\beta\s S_{\alpha\wedge\beta}$ for all $\alpha,\beta\in Y$.
\end{itemize}

Suppose now $S$ is {\it any}  semigroup that satisfies (38). We then say that $S$ is an {\it (ordinary) semilattice} $Y$ of subsemigroups $S_\alpha\ (\alpha\in Y)$. Here comes the {\bf Structure Theorem for c.f. semigroups:}

\vspace{3mm}

{\bf Theorem 9: }{\it Each commutative finite semigroup $S$ is a semilattice $Y$ of its Archimedean components $A_e\ (e\in Y)$.}

\vspace{2mm}
Before giving the proof, observe how Theorem 9 accomodates the semigroups treated in Sections 2,3 and 5; that is: (i) nilsemigroups, (ii) Abelian groups, (iii) semilattices. In case (i) the semilattice $Y$ in Thm. 9 is trivial ($Y=\{e\}$) and the unique Archimedean component $A=A_{e}$ has $K(A)=\{\zero\}$. In case (ii) again $Y=\{e\}$ but now $K(A)=A$. In case (iii) we have $Y=S$ and  $A_{e}=\{e\}$ for all $e\in Y$.

\vspace{3mm}
{\it Proof.} Let us first find the semilattice $Y$.
For all $e,f\in E(S)$ it holds that $(ef)^2=efef=eeff=ef$, and so $E(S)$ is a ssgr of $S$. Therefore $Y:=E(S)$  is a semilattice.
Recall from (26) that $ef=e\wedge f$ for all $e,f\in Y$. We use the $\wedge$ notation to emphasize the poset aspect of $Y$. Similar to the arguments in 7.2 and 7.3 we have:

\begin{itemize}
    \item[(39)] $ (x,\in A_e,\ y\in A_f) \Ra (x^k=e,\ y^\ell=f)$
    \item[] $\ \Ra (xy)^{k\ell}=(x^k)^\ell(y^\ell)^k=ef\Ra xy\in A_{ef}=A_{e\wedge f}.$
\end{itemize}

\noindent
This shows that condition (38) is satisfied. $\square$

\vspace{2mm}
From $K(A_e)\s A_e,\ K(A_f)\s A_f$ follows $K(A_e)K(A_f)\s A_{ef}$. One can show (10.6.1) that even $K(A_e)K(A_f)\s K(A_{ef})$. Furthermore, if $e'\in E(S)$ is the smallest element of $(Y,\wedge)$ then (why?) $K(S)=K(A_{e'})$. In Section 9 we persue sufficient and necessary conditions for the semilattice
$Y$ in Theorem 9 to be strong.

\vspace{3mm}
{\bf 8.2}
In 7.4.1 we found that $S:=\Z_7\times\Z_8\times\Z_9\simeq \Z_{504}$ has 8 Archimedean components $A_\alpha$ which we indexed with the elements $\alpha$ of some  set $Y_8$. We now understand that $Y_8\simeq\{\zero,\one\}^3$ is the semilattice $Y$ postulated in Theorem 10. More generally, if $S=\Z_n\simeq \Z_{p_1^{\gamma_1}}\times\cdots\times  Z_{p_t^{\gamma_t}}$, then\footnote{That relates to 5.3.2 whereby 
{\it each} f. semilattice $Y$ is a ssgr of $\{\zero,\one\}^t$.}
$Y=E(\Z_n)\simeq \{\zero,\one\}^t$.

In the same vein each Arch. component $A$ of $\Z_n$ is a 
direct product of $t$ semigroups, each one of which being a group or being nil. 
Note that for $A$ being a group it is necessary and sufficient that all
 nilsemigroups entering $A$ are {\it trivial}.
Therefore {\it all} Arch. components $A$ of $\Z_n$ are groups 
iff $\gamma_1=\cdots=\gamma_t=1$, i.e. iff $n$ is squarefree. Put another way, $(\Z_n,\odot)$ is a semilattice of groups  iff $n$ is square-free. 

\vspace{3mm}
{\bf 8.3}  So much about the peculiarities of $(\Z_n,\odot)$. Next comes a five step recipe to classify {\it any} finite commutative semigroup $S$ whose multiplication table is known:

\begin{itemize}
    \item[(I)] Determine the universes (=underlying sets) of the Archimedean components $A_e\s S$, along with their unique idempotents.
    \item[(II)] Calculate the poset structure of the semilattice $Y=E(S)$.
    \item[(III)] For each  $A_e\ (e\in Y)$ found in (I) determine the universe of its
    kernel $K(A_e)$.
    \item[(IV)] For all $e\in Y$ calculate the poset structure of the nilsemigroup $A_e/K(A_e)$.
    \item[(V)] For all $e\in Y$ calculate the structure of the Abelian group $K(A_e)$.
\end{itemize}

\vspace{3mm}

{\bf 8.3.1} Let us illustrate the details on a semigroup obtained in 6.2, i.e.

$$RF_1=\{a,b,c,ab,ac,b^2,bc,ab^2,abc,b^2c,ab^2c\}$$

\noindent
with presentation 

$$\{a^2\to a,\ b^3\to ab^2,\ c^2\to bc\}.$$

As to (I), we partition $RF_1$ into "connected pieces" by picking elements $x,y,..$ at random from $RF_1$ as follows. Starting e.g. with  $x:=a$  we get $\langle x\rangle=\{a\}$. Next $x:=b$ yields
$\langle x\rangle=\{b,b^2,b^3\}=\{b,b^2,ab^2\}$. Here (and henceforth) we rely on 6.2 where we found that $b^4=b^3$, which reduces to $ab^2$. Next $x:=c$ yields $\langle x\rangle=\{c,c^2,c^3,c^4\}=\{c,bc,b^2c,ab^2c\}$. And $y:=ab$ yields $\langle y\rangle=\{ab,(ab)^2\}=\{ab,ab^2\}$. Here $y$ is special in that some power (in fact $y^2$) coincides with some
{\it previously obtained} element. Similarly $y:=ac$ has such a power $y^3$, namely 
$\langle y\rangle=\{ac,(ac)^2,(ac)^3\}=\{ac,abc,ab^2c\}$.

\vspace{1cm}
\includegraphics[scale=0.95]{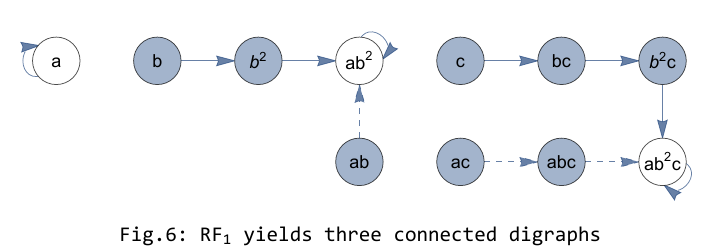}

\vspace{2mm}
Let us explain why generally the node sets $D_i$ of the connected digraphs\footnote{Recall: digraph=directed graph. Furthermore,  "connected" has an obvious meaning (i.e. concerns the underlying {\it undirected} graph) and  must not be confused with the more restricted  concept "strongly connected".} obtained this way are the (universes of the) Arch. components of the semigroup $S$ at hand. The "birth" of any fixed connected $D_i$ is some set $\langle x\rangle$. Hence there is a unique $e\in E(S)$ with $e\in\langle x\rangle\s A_e$. The first augmentation of $\langle x\rangle$ is by a set $\{y,y^2,..,y^k\}\ (k\ge 2)$ which intersects
$\langle x\rangle$ in $y^k$. If we had $y\in A_f\ (f\neq e)$ then all powers $y^i$ would remain in the ssgr $A_f$, which contradicts $y^k\in A_e$. Hence $\{y,y^2,..,y^k\}\s A_e$. The same reasoning applies to all further augmentations, and so $D_e:=D\s A_e$. And this holds for all $e\in E(S)$. Since all of $S$ gets partitioned into node sets of connected digraphs $D_e$, we have 

$$\bigcup\Bigl\{D_e:\ e\in E(S)\Bigr\}=S=\bigcup\Bigl\{A_e:\ e\in E(S)\Bigr\},$$

\noindent
which in view of $D_e\s A_e$
forces $D_e= A_e$ for all $e\in E(S)$.

\vspace{3mm}
{\bf 8.3.2} As to (II), from (I) we know that $Y=E(RF_1)=\{a,ab^2,ab^2c\}$ (the white nodes in Fig. 6). From $a\cdot ab^2=ab^2$ and $ab^2\cdot ab^2c=a^2b^4c=a^2ab^2c=ab^2c$ follows that $ab^2c<ab^2<a$, and so the semilattice $(Y,\le)$ is a 3-element chain.

\vspace{3mm}
{\bf 8.3.3} As to (III), recall that $K(A_e)$ is an ideal of $A_e$ with identity $e$. Hence for all $x\in A_e$ it holds that ($x\in K(A_e)$ iff $ex=x$). In our toy example one checks that it holds for all $e\in Y$ and all $x\in A_e$  that ($ex=x$ iff $x=e$). Therefore $K(A_e)=\{e\}$, i.e. each $A_e\ (e\in Y)$ is nil with zero $e$.
For succinctness we put

$$0_1:=a,\hspace{3mm} 0_2:=ab^2,\hspace{3mm} 0_3:=ab^2c$$

\vspace{3mm}
{\bf 8.3.4} As to (IV), let us start with the nil semigroup $A_{0_3}$. At first sight the structure of $(A_{0_3},\le_{\cal J})$ (rightmost part of Fig.7) seems obvious: $b^2c$ is below $bc$ because it is the multiple $b\cdot bc$ of $bc$, and similarly in all other cases. Trouble is, $b\not\in A_{0_3}$ and only multipliers in $A_{0_3}$ can be used (the whole of $RF_1$ is {\it not} partially ordered by $\le_{\cal J}$). 

\vspace{1cm}
\includegraphics[scale=0.95]{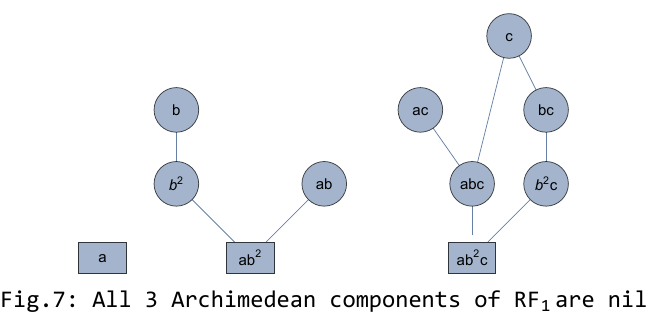}

\vspace{7mm}
\noindent
The justification of $(A_{0_3},\le_{\cal J})$ must hence rely on the multiplication  within $A_{0_3}$.
Its multiplication table (omitting $0_3$) is given below.

\vspace{3mm}
\begin{tabular}{c||c|c|c|c|c|c|}
			 & $c$ & $bc$ & $b^2c$ & $ac$ & $abc$   \\ \hline\hline
 $c$ &  $bc$ &  $b^2c$ & $0_3$ & $abc$ & $0_3$  \\ \hline
  $bc$ &  $b^2c$ &  $0_3$ & $0_3$ & $0_3$ & $0_3$  \\ \hline
   $b^2c$ &  $0_3$ &  $0_3$ & $0_3$ & $0_3$ & $0_3$  \\ \hline
    $ac$ &  $abc$ &  $0_3$ & $0_3$ & $abc$ & $0_3$  \\ \hline
     $abc$ &  $0_3$ &  $0_3$ & $0_3$ & $0_3$ & $0_3$  \\ \hline
    \end{tabular}

\vspace{3mm}
\noindent
Adopting the notation of 2.7.3  we read off that $PM(abc)=PM(b^2c)=\{0_3\}$, and so $abc,\ b^2c$ are the upper covers of $0_3$ in the poset $(A_{0_3},\le_{\cal J} )$. The remainder of the diagram ensues from $PM(ac)=\{0_3,abc\},\ PM(bc)=\{0_3,b^2c\},\ PM(c)=\{0_3,bc,b^2c, abc\}$.
In similar fashion one finds the poset $(A_{0_2},\le_{\cal J} )$. And $A_{0_1}=\{0_1\}$ is trivial.

\vspace{3mm}
{\bf 8.3.5}
As to (V), this is easy here; all three  kernels $K(A_e)\ (e\in Y)$ are 1-element groups.
 A systematic treatment of step (V) comes in 8.7.

\vspace{3mm}
{\bf 8.4} Let us apply the recipe (I) to (V) to the sgr 
$RF_2=\{a,b,a^2,ab,b^2,a^2b,b^3\}$ in (31) with presentation $\{b^4\to b^2,\ a^3\to b^2,\ a^4\to a,\ ab^2\to a\}$.

\vspace{3mm}

  As to (I), this time only one connected digraph, i.e. one Archimedean component $A=RF_2$ arises. Its unique idempotent is $b^2$ (Fig. 8A). 
  
  As to (II), we have the trivial semilattice $Y=\{b^2\}$.
  
   As to (III), upon checking that $b^2x=x$ for all $x\neq b$, we conclude that $K(A)=A\setminus\{b\}$. 
   
   Hence, (IV), the Rees quotient $A/K(A)$ is a 2-element nilsgr (Fig. 8B). 
   
   As to (V), one verifies that\\ 
   $K(A)=\langle ab\rangle=\{ab,(ab)^2,...,(ab)^6\}=\{ab,a^2,b^3,a,a^2b,b^2\}$, and so $K(A)\simeq C_{6}$.

\includegraphics[scale=0.55]{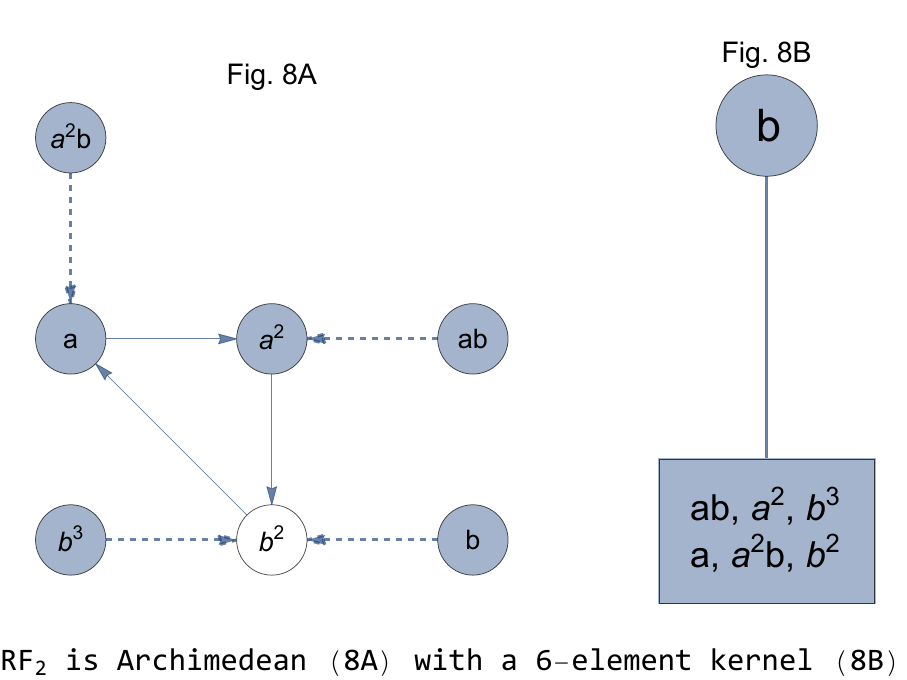}

\vspace{3mm}
{\bf 8.5} Carrying out steps (I) and (II) for $(\Z_{18},\odot)$ yields the semilattice 
$Y=\{0,1,9,10\}$ and the Arch. components

$$\{1,5,7,11,13,17\},\ \{2,4,8,10,14,16\},\ \{3,9,15\},\ \{0,6,12\}.$$

\noindent
Steps (III) to (V) show that the two large components are groups and the smaller ones are  nilsemigroups. Because $\Z_{18}$ is of type $\Z_n$, much of this was  predictable in view of Section 7.
Namely, $\Z_{18}\simeq \Z_2\times \Z_9$, and the semilattice for $\Z_2\times \Z_9$ is
$E(\Z_2)\times E(\Z_9)\simeq \{0,1\}^2$.
 Further the Arch. components of $\Z_2$ are $\{0\},\{1\}$, while the ones of $\Z_9$ are the group $G=\Z_9^{inv}\ (\simeq C_{6})$ and the nilsemigroup $N=\{0,3,6\}$. Hence  the 4 Arch. components of $\Z_2\times \Z_9$ are the groups $\{0\}\times G$ and $\{1\}\times G$, as well as the nilsgr  $\{0\}\times N$ and $\{1\}\times N$. As to "much of this was predictable", what remains is the translation from $\Z_2\times \Z_9$ to $\Z_{18}$ (see 3.5).

\vspace{3mm}
{\bf 8.6} Akin to 6.2.1 the {\it relatively free Abelian group} 
$RFAG(x,y,..,z:\{...\})$ is the largest  Abelian group generated by $x,y,..,z$ and subject to a set $\{...\}$ of postulated relations. For instance (and generalizing in obvious ways)

$$(40)\quad  RFAG(x,y: \{x^5=\one,\ y^7=\one\})\ is\ isomorhic\ to\ C_5\times C_7$$

\noindent
since $C_5\times C_7=\langle a\rangle\times\langle b\rangle$  is  (i) generated by $x:=(a,\one),\ y:=(\one,b)$, and (ii) satisfies $x^5=(\one,\one),\ y^7=(\one,\one)$, and (iii) every other Abelian group with analogous generators $x',y'$ is an epimorphic image of $C_5\times C_7$ (proven as in 10.3).
In principle each finite Abelian group $G$ is of type $G\simeq RFAG(..)$  since the whole Cayley table readily yields a (highly redundant) presentation.

\vspace{2mm}
{\bf 8.6.1} More relevant however is the opposite. Thus suppose the Abelian group $G$ is presented by generators and relations as in (40). How can one unravel the way $G$ expands as a direct product of cyclic groups?

Using additive notation (so $\zero$ is now the identity) let us illustrate\footnote{The remainder of 8.6 is based on [Ar,chapter 22]. Some apriori "miraculous" moves will be demystified in 8.6.2.} the details  on

$$RFAG(x,y,z:\ 60x-112y+94z=\zero,\ 56x-108y+92z=\zero,\ 84x-160y+136z=\zero)$$

\noindent
which we abbreviate as $RF_6$. Here $n=3$ and the $m=3$ relations we like to render in matrix form

$$(41)\quad \begin{pmatrix}60&-112&94\\ 56&-108&92\\ 84&-160&136\end{pmatrix}
\begin{pmatrix} x\\y\\z\end{pmatrix}=\begin{pmatrix} 0\\0\\0\end{pmatrix}$$

\noindent
Let us first switch to another generating set $\{x',y',z'\}\s RF_6$ which will turn out more suitable and which is implicitely defined by

$$(42)\quad \begin{pmatrix} x\\y\\z\end{pmatrix}\ =:\ \begin{pmatrix}-1&2&1\\ 2&1&2\\ 3&0&2\end{pmatrix}
\begin{pmatrix} x'\\y'\\z'\end{pmatrix}$$

\noindent
As to "implicitely defined", an explicit definition of $x',y',z'$ in terms of $x,y,z$ would involve the inverse of the above square matrix  (more on that in 8.6.2). Combining (41) and (42) we get

$$ \begin{pmatrix}-2&8&24\\ 4&4&24\\ 4&8&36\end{pmatrix}
\begin{pmatrix} x'\\y'\\z'\end{pmatrix}=
\begin{pmatrix}60&-112&94\\ 56&-108&92\\ 84&-160&136\end{pmatrix}
\begin{pmatrix}-1&2&1\\ 2&1&2\\ 3&0&2\end{pmatrix}
\begin{pmatrix} x'\\y'\\z'\end{pmatrix}$$

$$\stackrel{(42)}{=}\begin{pmatrix}60&-112&94\\ 56&-108&92\\ 84&-160&136\end{pmatrix}
\begin{pmatrix} x\\y\\z\end{pmatrix}\stackrel{(41)}{=}\begin{pmatrix} 0\\0\\0\end{pmatrix}$$

\noindent
We see (not surprising) that the new generators also satisfy new relations:

$$ \begin{matrix}-2x'+8y'+24z'=0 & (R1)\\
4x'+4y'+24z'=0 & (R2)\\
4x'+8y'+36z'=0 & (R3)
\end{matrix}
$$ 

The new relations $(R1),(R2),(R3)$  behave  better than the old ones
insofar that  now (invertible) integer combinations of $(R1),(R2),(R3)$ can be found that do the job:

$$(43)\quad \begin{matrix} 
-(R1)-2(R2)+2(R3) & is & {\bf 2x'}+0y'+0z'=0\\
-2(R1)-7(R2)+6(R3) & is & 0x'+{\bf 4y'}+0z'=0\\
2(R1)+6(R2)-5(R3) & is & 0x'+0y'+{\bf 12z'}=0
\end{matrix}$$

\noindent
Since all manipulations are invertible (see 8.6.2), we conclude that
$RF_6\simeq RFAG(x',y',z': 2x'=\zero,\ 4y'=\zero,\ 12z'=\zero)$ which, arguing as in (40), is isomorphic\footnote{Observe this the unique type of direct product that uses $t_{min}$ factors (see 3.6.3). If the initial generators $x,y,z$ of $RF_6$ are given in some concrete format (say by matrices), then also the generators $x',y',z'$ can be obtained explicitely by applying the inverse of the matrix in (42).} to $C_2\times C_4\times C_{12}$

\vspace{3mm}
{\bf 8.6.2} Some further explanations are in order. Let $A,\ B$ be the square matrices appearing in (41),(42) respectively.  Additionally let $C$ be the matrix that encodes the integer combinations of $(R1),(R2),(R3)$, i.e. $C:=\begin{pmatrix}-1&-2&2\\ -2&-7&6\\ 2&6&-5\end{pmatrix}$. Finally, if $D$ is the diagonal matrix with entries $2,4,12$, then it follows from  (41),(42),(43) that $CAB=D$.
We mentioned already that the integer-valued matrix $B$ needs to be invertible. Moreover the entries of $B^{-1}$ must be integers as well (why?). Necessary and sufficient  for such a $B^{-1}$ to exist, is the {\it unimodularity} of $B$, i.e. $det(B)\in \{1,-1\}$. Similarly, since the relations derived from 
$(R1),(R2),(R3)$ must be fit to rederive $(R1),(R2),(R3)$, the transition matrix $C$ must be unimodular as well.

 One says an $m\times n$ matrix $D=(d_{i,j})$ has {\it Smith Normal Form} if its only nonzero entries  are non-negative integers $d_i:=d_{i,i}\ (1\le i\le t)$ such that $d_i$ divides $d_{i+1}$ for all $1\le i<t\le min\{m,n\}$.
There is an algorithm (see [Ar] or [RG,ch.2]) that brings each integer-valued $m\times n$ matrix $A$  into  Smith Normal Form by applying suitable row and column operations to $A$. The row operations have the same effect as left multiplication by an unimodular $m \times m$ matrix $C$, and the column operations can be simulated by right multiplication with a unimodular $n\times n$ matrix $B$. 

\vspace{3mm}
{\bf 8.7} We are now in a position to tackle more systematically step (V) of our recipe in 8.3. Since the elements of each fixed Abelian group $G:=K(A_e)$ triggered by our c.f. sgr $S'$ are given by normal forms, one can readily calculate the order $o(x)$ for each $x\in G$. By 3.7.2 this reveals the structure of $G$.
This is nice enough yet doesn't yield a trimmed generating set  $X\s G$. 

\vspace{3mm}
{\bf 8.7.1} As glimpsed in 3.7.3, the latter is achieved in [S], which also surveys the 50 year old history of the problem. In particular [S,p.478] it is stated that most prior attempts first constructed a relation matrix (this being the time-intensive part), and then reduced it to Smith Normal Form as illustrated in 8.6. While Sutherland's method might, as he claims, usually be faster, our particular scenario might be an exception but this requires further research. 

In a nutshell, this is why. The algorithm in [FP] achieves the following. Given the generators of a {\it concrete}\footnote{Roughly speaking, the generators must be elements of a larger semigroup $T$ in which "concrete computations" can be carried out. For instance $T$ may be the translation semigroup on 7 elements (having cardinality $7^7$), or the semigroup of all $2\times 2$ matrices over $\Z_{59}$.} semigroup $S$, it finds a (usually small) semigroup presentation of $S$. For our particular concrete semigroup $S=K(A_e)$ the [FD] algorithm likely speeds up by two reasons. First, commutativity always helps. Second, recall, $K(A_e)$ is a subgroup of a semigroup of type $RFCS(...)$, i.e. of a semigroup which {\it has} already a semigroup presentation.

\section{Ideal extensions of one cyclic semigroup by another}

\vspace{3mm}
 In the best of all worlds each semilattice of semigroups would be\footnote{This is why. For each covering $\alpha\succ\beta$ in $Y$ storing the definition of  $\sigma_{\alpha,\beta}$ is more economic (and insightful) than storing the full multiplication table of $S$. In a similar vein solvable groups are more economic than the ordinary kind.}  a {\it strong} semilattice of semigroups. Unfortunately, statistically speaking most semilattices of semigroups are {\it not} strong.

But there are  beams of light. Suppose $(S,\ast)$ is a semilattice $Y$ of  {\it monoids} $A_e\s S$ (where $e\in Y\s S$). Thus $e$ is the identity of $A_e$. For all $d\ge e$ in $Y$ and $a,b\in A_d$ one has $(a\ast b)\ast e=a\ast b\ast e\ast e=(a\ast e)\ast (b\ast e)$, and so $a\sigma_{d,e}:=a\ast e$ defines a morphism $A_d\to A_e$. In particular, if $d=e$ then $a\sigma_{e,e}=a$ for all $a\in A_e$.
For the sake of readibility in $(44)$ we put $ed:=e\ast d$. Then for all $e,d\in Y$ and
$a\in A_e,\ b\in A_d$ it follows from $a\ast b\in A_{ed}$ that

$$(44)\quad a\ast b=(a\ast b)\ast ed=[a\ast ed]\ast[b\ast ed]=[a\sigma_{e,ed}]\ast [b\sigma_{d,ed}].$$

\noindent
Therefore each semilattice of monoids is "automatically" strong. For instance (see 8.2), if $n$ is squarefree, then $(\Z_n,\odot)$ is a strong semilattice of groups.
Observe that "up to its tail" each  $C_{m,n}$ is a monoid (even group) $C_n$.  This inspires  the following questions:
Is a finite semilattice of {\it cyclic} semigroups automatically a strong semilattice? If not, how far off is it? We return to this issue at the end of Section 9, having dwelled on the case "finite semilattice = 2-element semilattice" and (generalizing) on the so called Ideal Extension Problem. Much  work was done already in Section 5.

\vspace{2mm}
{\bf 9.1} In this Subsection we adopt the notation of [CL,p.137] for ease of comparison.
Let $S$ be a semigroup and $T$ a disjoint semigroup with zero $\zero$. Put $T^*:=T\setminus\{\zero\}$. There may be {\it zero-divisors} $a,b\in T^*$ in the sense that $ab=\zero$.  
A semigroup $(\Sigma,\circ)$ on the set $\Sigma:=S\cup T^*$  is an {\it ideal extension of $S$ by $T$} if $S$ is an ideal of $\Sigma$ and the Rees quotient $\Sigma/S$ (see 2.9) is isomorphic to $T$.  Roughly speaking, {\it inflating $\zero$ to $S$ inflates $T$ to $\Sigma$} (see Fig.9). 

\includegraphics[scale=0.45]{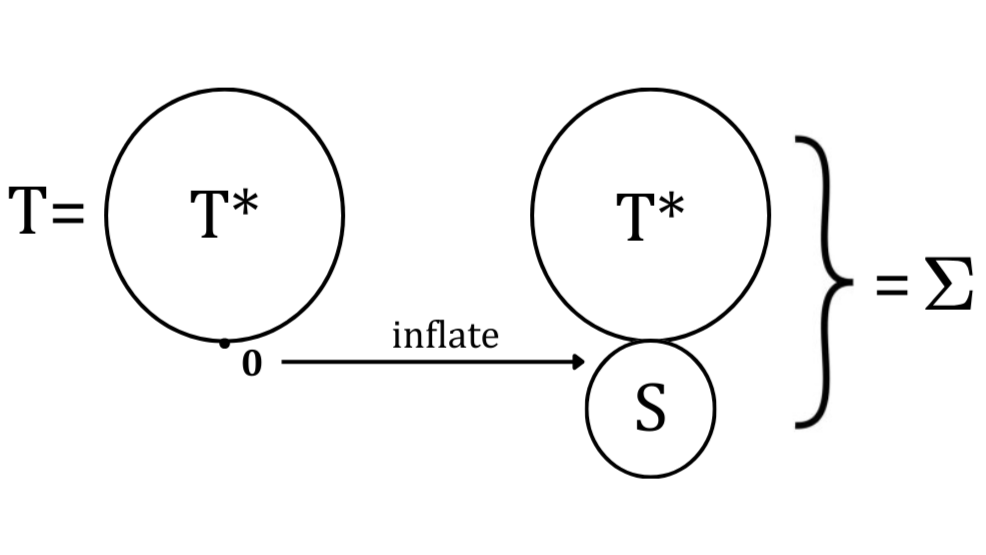}
 {\Large{ \tt Figure 9: Ideal extension of $S$ by $T$}}

\vspace{11mm}

Conversely, given disjoint semigroups $S$ and $T\ (\ni\zero)$, here comes the cheapest way to {\it obtain} an ideal extension $\Sigma$ of $S$ by $T$. Put $\Sigma:=S\uplus T^*$ and define $a\circ b=b\circ a:=b$ for all $a\in T^*,\ b\in S$ (and otherwise don't change the multiplication). It is easy to verify the associativity of $\circ$. One calls this the {\it trivial} ideal extension  of $S$ by $T$.

 More subtle, given again disjoint semigroups $S$ and $T\ (\ni\zero)$, another way to obtain an ideal extension $\Sigma$ of $S$ by $T$, is to look  for a
{\it partial morphism}, i.e. a map $\varphi:T^*\to S$ such that $(ab)\varphi=(a\varphi)(b\varphi)$ whenever $ab\neq \zero$. Having found $\varphi$ define the groupoid $(\Sigma,\circ)$ as follows.

\begin{itemize}
    \item[(i)] $\ a\circ b:=ab$ for all $a,b\in T^*$ with $ab\neq\zero$;
    \item[(ii)] $\ a\circ b:=(a\varphi)(b\varphi)$ for all $a,b\in T^*$ with $ab=\zero$;
    \item[(iii)] $\ a\circ b:= (a\varphi)b$ for all $a\in T^*$ and $b\in S$;
    \item[(iv)] $\ a\circ b:= a(b\varphi)$ for all $a\in S$ and $b\in T^*$;
     \item[(v)] $\ a\circ b:=ab$ for all $a,b\in S$.
\end{itemize}

\noindent
(In the commutative case (iii) and (iv) are equivalent.)
Theorem 4.19 in [CL] establishes, based on eight straightforward subcases, that $\circ$ is associative\footnote{A concrete calculation occurs in 9.3.}, i.e. $(\Sigma,\circ)$ is an ideal extension of $S$ by $T$.

Furthermore, Theorem 4.19 states: If $S$ happens to have an {\it  identity} then {\it each} ideal extension of $S$ by $T$ is induced\footnote{One nice consequence is that each Archimedean semigroup (Section 7) is a $\varphi$-induced ideal extension of an Abelian group by a nil semigroup.} by a partial morphism $\varphi$ as above; the proof is short and similar to (44). 

A quick definition before we can continue. If any sgr $H$ contains a zero then put $H^0:=H$. If not, then by definition $H^0:=H\uplus\{\zero\}$ is the sgr with zero $\zero$ where $xy$ is the same in $H$ and $H^0$ for all $x,y\in H$. (In likewise fashion one can "adjoin an identity" to a semigroup $H$ and obtain a monoid $H^{\one}$.)

\vspace{3mm}
{\bf 9.1.1}
Back to ideal extensions. A noteworthy special case arises if $T$ has no zero-divisors, i.e. if $T^*$ is a ssgr of $T$. Then every ideal extension $\Sigma$ of $S$ by $T$ is a 2-element semilattice $\{S,T^*\}$ of semigroups. (This semilattice has zero $S$ and identity $T^*$.)
Conversely, every 2-element 
semilattice $\{S,T^*\}$ of semigroups can be viewed as an  ideal extension of $S$ by $T:=(T^*)^0=T^*\uplus\{\zero\}$.

The named case takes place in the remainder of Section 9 with $T^*:=C_{m,n}$  and $S:=C_{m',n'}$. Specifically, in Section 9 all ideal extensions $\Sigma$ of $S$ by 
$T:=C_{m,n}\uplus\{\zero\}$ will be determined. This extends [CL,Thm.4.19] in two ways. First, we unravel when a partial morphism $\varphi$ {\it exists}. Second, all ideal extensions $\Sigma$ which are {\it not based} on  partial morphisms will be identified. To get rid of the clumsy $\zero$ of $T$ more suitable terminology will be adopted, starting in 9.2.
(Yet in 9.5 we briefly hark back to [CL] terminology to take stock.)

\vspace{2mm}
{\bf 9.2} Let $m,n,m',n'\ge 1$ and $k\in\{0,1,...,m'+n'-1\}$ be integers. We then call $Q:=(m,n,m',n';k)$ a {\it quintuple}. It is {\it realizable} if there is a semigroup $\langle a,b\rangle$ (i.e. generated by $a,b$) such that

\begin{itemize}
    \item[(45.1)] $\langle a,b\rangle=\langle a\rangle\uplus \langle b\rangle$
    \item[(45.2)] $\langle a\rangle\simeq C_{m,n}$
    \item[(45.3)] $\langle b\rangle\simeq C_{m',n'}$
    \item[(45.4)] $ab=b^{k+1}$
    \end{itemize}

By abuse\footnote{Properly speaking $\Sigma$ is an ideal extension of $\langle b\rangle$ by $\langle a\rangle\uplus\{\zero\}$. In $\Sigma$ evidently  $ab=b^{k+1}$ for some 
$k\in\{0,1,...,m'+n'-1\}$. By saying that $\Sigma$ realizes $Q=(m,n,m',n';k)$ we additionally point out which $k$ occurs. } of language we call $\Sigma:=\langle a,b\rangle$ an {\it ideal-extension} that {\it realizes} $Q$.  We stress that (45.4) uniquely determines the multiplication in $\Sigma$. For instance $a^2b=(aa)b=a(ab)=ab^{k+1}=(ab)b^k=b^{2k+1}$, and this e.g. implies
$a^2b^5=b^{2k+5}$. We postpone the (easy) inductive argument for
 $a^ib^j=b^{ik+j}$ to 9.4.
It follows that for fixed semigroups $C_{m,n}$ and $C_{m',n'}$ we can (and will) classify all ideal extensions of $C_{m',n'}$ by $C_{m,n}\uplus\{\zero\}$ as follows: 

\vspace{2mm} Determine all $k\in\{0,1,...,m'+n'-1\}$ for which the quintuple $(m,n,m',n';k)$ is realizable!

\vspace{2mm}
{\bf 9.2.1} The simplest case is $k=0$. Then (45.4) becomes $ab=b$. Therefore each ideal-extension $\Sigma$  that realizes $Q:=(m,n,m',n';0)$
 satisfies $a^i b^j=b^j$ (set $k=0$ in $a^ib^j=b^{ik+j}$) and hence is the trivial ideal extension.

Therefore each quintuple $Q$ that has $k=0$ will be called {\it trivial}. If $k\ge 1$ then $Q$ is {\it nontrivial }. 

 \vspace{2mm}
 {\bf 9.2.2} Consider the quintuple is $Q:=(m,n,m',n';k):=(m,n,m',1;m'-1)$.  By (45.4) one has $ab=b^{m'}$. But $b^{m'}$ is the zero $\zero$ of $C_{m',1}$. Take a set of symbols $\Sigma:=\{a,a^2,..,a^{m+n-1},b,b^2,..,b^{m'-1},\zero\}$, define $a^i*a^j$ and $b^i*b^j$ in the obvious way, and put $a^i*b^j:=\zero$. Checking the associativity of $*$ is as easy as in 9.2.1. It follows that $(\Sigma,*)$ realizes $Q$. Notice that
 $Q$ is nontrivial iff $m'>1$.
 
\vspace{3mm}
{\bf 9.3}
Let $f:C_{m,n}\to C_{m',n'}$ be any morphism. According to [CL,Thm.4.19] it triggers a particular kind of ideal extension $\Sigma$. Specifically, by (iii) above the multiplication in $\Sigma$ is given by

    $$(46)\quad a^i\cdot b^j:=(a^if)b^j.$$
    

\noindent
From (46) follows that $a\cdot b=(af)b$. There is a unique $k\in \{1,2,...,m'+n'-1\}$ with
$af=b^k$. Obviously  the quintuple $Q_0:=(m,n,m',n';k)$ is realizable.
Upon a change of notation (e.g. $a_\alpha:=a^i,\ \sigma_{\alpha,\beta}:=f$) equation (46) becomes $(28')$, and so we are dealing with a {\it strong} 2-element semilattice of $C_{m,n}$ and $C_{m',n'}$. 

This motivates the following definition. If $Q=(m,n,m',n';k)$ is any nontrivial quintuple and there is a morphism $f:C_{m,n}\to C_{m',n'}$ with $af=b^k$ then $Q$ is {\it strongly} realizable. In view of Theorem 1 strong realizibility  takes place iff:

\begin{itemize}
    \item[(SR1)] $m'\le mk$ 
    \item[(SR2)] $n'$ divides $nk$
\end{itemize}

For instance suppose the quintuple $(m,n,m',1;m'-1)$ in 9.2.2 is nontrivial. It is strongly realized if $m\ge 2$: (SR1) holds since $m'=k+1\le mk$, and (SR2) holds since $n'=1$. So what is the underlying morphism $f$? Obviously it is $(\forall i)\ a^if:=\zero$. (This is a morphism of the 2.3.1 kind.)  

\vspace{3mm}
{\bf 9.4} So much about strong realizability. But what about  {\it ordinary} realizability? 
Lemma 10 below proves the necessity\footnote{Observe that (SR1),(SR2) above are not only sufficient for ordinary realizability, but "almost" 
necessary as well.} of certain conditions (R1) and (R2), while Theorem 11 will establish their sufficiency. Each trivial quintuple being realizable, demanding nontriviality in Lemma 10 is hardly a restriction, but will be necessary in its proof.

\vspace{2mm}
{\bf Lemma 10: }{\it If the nontrivial quintuple $Q:=(m,n,m',n';k)$ is realizable, then
\begin{itemize}
    \item[(R1)] $m'-1\le mk$ 
    \item[(R2)] $n'$ divides $nk$
\end{itemize}  }

\vspace{2mm}
{\it Proof.}  We first show that if $\Sigma=\langle a,b\rangle$ realizes $Q$, then it holds that

$$(47)\quad a^tb=b^{tk+1}\ for\ all\ t\ge 1.$$

\noindent
Indeed, for $t=1$ this becomes $ab=b^{k+1}$, which holds by (45.4). By induction assume that $a^{t-1}b=b^{(t-1)k+1}$ for some $t\ge 2$. Then 

$$a^tb=(ab)b^{(t-1)k}=b^{k+1}b^{tk-k}=b^{tk+1}.$$

\noindent
By assumption $\langle a\rangle\simeq C_{m,n}$. Hence $a^{m+n}=a^m$, and so

$$b^{mk+1}\stackrel{(47)}{=}a^mb=a^{m+n}b\stackrel{(47)}{=}b^{(m+n)k+1}.$$

\noindent
This, together with the assumption $k\ge 1$, shows that $b^{mk+1}$ is in the body of $\langle b\rangle\simeq C_{m',n'}$, and so $m'\le mk+1$, which is (R1).

Furthermore, $mk+1\equiv(m+n)k+1$ modulo $n'$, hence\\
$mk\equiv mk+nk\ (mod\ n')$. This implies that $n'$ divides $nk$ (statement (R2)). $\square$

\vspace{2mm}
It follows from (47) that every realizing ideal-extension $\Sigma$ of $Q$ in Lemma 10 satisfies

$$(47')\quad a^ib^j=b^{ik+j}\ for\ all\ i,j\ge 1.$$

\vspace{2mm}
{\bf 9.5}
In order to see that (R1),(R2) are also sufficient for realizibility we  show that a certain groupoid is in fact a semigroup. For any groupoid $(S,\cdot)$ we say that $a\in S$ {\it associates with everybody} if
$(x\cdot a)\cdot y=x\cdot(a\cdot y)$ for all $x,y\in S$. As is well known, if $a$ and $b$ associate with everybody, then\footnote{It is likely that "everybody" can be weakened when $S$ is commutative, but how exactly?} so does $ab=a\cdot b$:

$$[x\cdot (ab)] y=[(xa)\cdot b]y=(xa)\cdot [by]=x[a\cdot(by)]=x[(ab)\cdot y]$$

\noindent
In particular, if $a,b$ generate $(S,\cdot)$, then $(S,\cdot)$ must be semigroup.

\vspace{3mm}
{\bf Theorem 11: }{\it A nontrivial quintuple $Q:=(m,n,m',n';k)$  is realizable iff (R1) and (R2) hold.}

\vspace{2mm}
{\it Proof of Theorem 11.} By Lemma 10 it remains to show that (R1),(R2) are sufficient. We define a groupoid  $(\Sigma,\ast)$ whose elements are the equivalence classes $[x]$ on a certain infinite set of independent (yet suggestively labeled) symbols $x$. Specifically, there will be $(m+n-1)+(m'+n'-1)$ classes whose "canonical" representatives, respectively, are

$$a,a^2,...,a^{m+n-1},b,b^2,...,b^{m'+n'-1}.$$

\noindent
Embracing the details, by definition $[a^i]:=\{a^i\}$ for all $1\le i<m$, and $[a^i]:=\{a^i, a^{i+n},a^{i+2n},...\}$ for all $m\le i\le m+n-1$. Likewise $[b^i]:=\{b^i\}$ for all $1\le i<m'$, and $[b^i]:=\{b^i, b^{i+n'},b^{i+2n'},...\}$ for all $m'\le i\le m'+n'-1$. Led by $(47')$ we define

$$ \Sigma:=\bigl\{[a],[a^2],...,[a^{m+n-1}],\ [b],[b^2],...,[b^{m'+n'-1}]\bigr\}\ as\ well\ as$$
$$ [a^i]*[a^j]:=[a^{i+j}],\quad  [b^i]*[b^j]:=[b^{i+j}],\quad [a^i]*[b^j]=[b^j]*[a^i]:=[b^{ik+j}].$$

That the first two cases of $*$ are well-defined is clear\footnote{As to the subgroupoids $\langle [a]\rangle$ and $\langle[b]\rangle$ of $\Sigma$ also being associative, see 6.1. It is also implied by the forthcoming argument.}. As to $[a^i]*[b^j]$ being well-defined, we fix any $[b^j]$,  let $[a^i]=[a^{i_0}]$, and strive to show that $[b^{ik+j}]=[b^{i_0k+j}]$. By assumption $i\equiv i_0\ (n)$. 
Hence $(i-i_0)k$ is divisible by $nk$. Since $n'$ divides $nk$ by (R2), it follows that
$(i-i_0)k\equiv 0\ (n')$. The latter implies 

$$(48)\quad ik+j\equiv i_0k+j\ (n')$$

\noindent
Our claim being trivial for $i=i_0$ we may assume that $i\neq i_0$. Then $i,i_0\ge m$, which together with (R1) implies $ik+j\ge mk+1\ge m'$, as well as $i_0k+j\ge mk+1\ge m'$.
Therefore $[b^{ik+j}]=[b^{i_0k+j}]$ in view of (48).

Similarly we fix any $[a^i]$, let $[b^j]=[b^{j_0}]$, and strive to show
that $[b^{ik+j}]=[b^{ik+j_0}]$. By assumption $j\equiv j_0\ (n')$, and so 

$$(49)\quad ik+j\equiv ik+j_0\ (n').$$

We can again assume that $j\neq j_0$, and so $j,j_0\ge m'$. From $ik+j,\ ik+j_0\ge m'$ and (49)
follows $[b^{ik+j}]=[b^{ik+j_0}]$.

\vspace{2mm}
Having checked that $\ast$ is well-defined, let us proceed to prove  associativity. Since the groupoid $(\Sigma,\ast)$ is generated by $[a],[b]$, it suffices to show that these two associate with everybody.

\vspace{2mm}
As to $[a]$, we thus need to show (dropping $\ast$) that $([x][a])[y]=[x]([a][y])$ for all $[x],[y]\in S$. 
Case 1: $[x]=[a^i],\ [y]=[a^j]$. Then $([a^i][a])[a^j]=[a^{i+1}][a^j]=[a^{i+1+j}]=[a^i]([a][a^j])$.

\vspace{2mm}
Case 2: $[x]=[a^i],\ [y]=[b^j]$. Then $([a^i][a])[b^j]=[a^{i+1}][b^j]=[b^{(i+1)k+j}]$, which coincides with $[a^i]([a][b^j])=[a^i][b^{k+j}]=[b^{ik+k+j}]$.

\vspace{2mm}
Case 3: $[x]=[b^j],\ [y]=[a^i]$. Then $([b^j][a])[a^i]=[b^{k+j}][a^i]=[b^{ik+k+j}]$, which coincides with $[b^j]([a][a^i])=[b^j][a^{i+1}]=[b^{(i+1)k+j}]$.

\vspace{2mm}
Case 4: $[x]=[b^i],\ [y]=[b^j]$. Then $([b^i][a])[b^j]=[b^{k+i}][b^j]=[b^{k+i+j}]$, which coincides with $[b^i]([a][b^j])=[b^i][b^{k+j}]=[b^{i+k+j}]$.

\vspace{3mm}
As to $[b]$, we  need to show  that $([x][b])[y]=[x]([b][y])$ for all $[x],[y]\in S$.

Case 1: $[x]=[a^i],\ [y]=[a^j]$. Then $([a^i][b])[a^j]=[b^{ik+1}][a^j]=[b^{kj+ik+1}]$,
which coincides with $[a^i]([b][a^j])=[a^i][b^{kj+1}]=[b^{ik+kj+1}]$.

\vspace{2mm}
Case 2: $[x]=[a^i],\ [y]=[b^j]$. Then $([a^i][b])[b^j]=[b^{ik+1}][b^j]=[b^{ik+1+j}]$, which coincides with $[a^i]([b][b^j])=[a^i][b^{j+1}]=[b^{ik+j+1}]$.

\vspace{2mm}
Case 3: $[x]=[b^j],\ [y]=[a^i]$. Then $([b^j][b])[a^i]=[b^{j+1}][a^i]=[b^{ik+j+1}]$, which coincides with $[b^j]([b][a^i])=[b^j][b^{ik+1}]=[b^{j+ik+1}]$.

\vspace{2mm}
Case 4: $[x]=[b^i],\ [y]=[b^j]$. Then $([b^i][b])[b^j]=[b^{i+1}][b^j]=[b^{i+1+j}]$, which coincides with $[b^i]([b][b^j])=[b^i][b^{j+1}]=[b^{i+j+1}]$. $\square$

\vspace{3mm}
{\bf 9.5.1}
 One may wonder whether Theorem 11 could be proven via generators and relations, akin to  9.2.1. While this can be done with little effort\footnote{If the relation $ab\to b$ in 9.2.1 gets replaced by $ab\to b^2$ (matching $k=1$), then the three relations are no longer locally confluent. But l.c. is recovered upon adding the (derivable) relation $a^{m'+n'-1}b\to a^{m'-1}b$. For $k>1$ recovering l.c. got cumbersome and no traces of pattern emerged.} for all quintuples having $k=1$ and satisfying (R1),(R2), the failure to extend the method to $k>1$ triggered the above proof of Theorem 11.

\vspace{3mm}
{\bf 9.5.2} Recall that each quintuple $Q$ has at most one realizer $\Sigma$ since (when $\Sigma$ exists) the multiplication in $\Sigma$ is uniquely determined by the $k$ in $Q$. Perhaps surprisingly, {\it different} $Q_1=(m,n,m',n';k_1)$ and $Q_2=(m,n,m',n';k_2)$ may have the {\it same} realizer $\Sigma$. Namely, putting $C_{m',n'}=\langle b\rangle$, this happens iff
$b^{k_1+1}=b^{k_2+1}$, hence iff $k_1=m'-1$ and $k_2=m'+n'-1$. In particular, if $m'=1$ then $k_1=0$ but $k_2\neq 0$. In other words, if (and only if) $C_{m',n'}=C_{1,n'}$ is a monoid, then the {\it trivial} ideal-extension of $C_{1,n'}$ by any $C_{m,n}$ can also be triggered by a {\it nontrivial} quintuple.

\vspace{2mm}
{\bf 9.5.3}
Here comes the gist of Theorem 1 and Theorem 11 phrased in the [CL] terminology of 9.1: Let $T:=C_{m,n}\uplus\{\zero\}$ and $S:=C_{m',n'}$. We leave trivial ideal extensions aside
(although by 9.5.2 they are "not so trivial"). Thus a nontrivial ideal extension of $S$ by $T$ exists iff there is a nontrivial $Q:=(m,n,m',n';k)$ satisfying (R1) and (R2).
A nontrivial $\varphi$-based ideal extension of $S$ by $T$ exists iff there is a nontrivial $Q:=(m,n,m',n';k)$ satisfying (SR1) and (SR2). 

To fix ideas, let $T=C_{3,9}\uplus\{\zero\}$ and $S=C_{13,18}$. Then there is some nontrivial ideal extension $\Sigma_1$ of $S$ by $T$ which  {\it is not}  $\varphi$-induced. 
But some other nontrivial ideal extension $\Sigma_2$ of $S$ by $T$   {\it is}  $\varphi$-induced. (Consider $Q_0=(3,9,13,18;k)$ with $k=4$, respectively $k=6$.)

 \vspace{3mm}
{\bf 9.5.4} So far divisibility concerned $n,n'$. Let us look at divisibility in relation to $m,m'$. For starters, reconsider $Q_0=(m,n,m',n';k)=(3,9,13,18;k)$.  Depending on $k$ the quintuple $Q_0$ is strongly realizable, or just realizable, or not realizable at all.
It holds that $m$ divides $m'-1$ (i.e. $3|12$).

 This is noteworthy in light of the following. Suppose the quintuple $Q:=(m,n,m',n';k)$ is realizable. If additionally $m$ does {\it not divide} $m'-1$, then $Q$ is {\it strongly} realizable.
To prove this, it suffices to verify that $m'\le mk$.  Since the assumption $m'-1=mk$ yields the contradiction $m|(m'-1)$, we conclude $m'-1\neq mk$. In view of $m'-1\le mk$ (due to (R1)) this forces $m'-1< mk$. But this implies $m'\le mk$.


\vspace{2mm}
{\bf 9.6} Let  us venture away from 2-element semilattices (=particular ideal extensions) to arbitrary finite semilattices\footnote{Now the letter $S$ conforms to the notation in 5.4, thus not the [CL] terminology of 9.1.} $S$ of cyclic semigroups $S_\alpha\ (\alpha\in Y)$. In 9.6.1 and 9.6.2 we deal with strong semilattices, in 9.6.3 with ordinary
semilattices, and 9.6.4 glimpses at the paper [AS].

\vspace{2mm}
{\bf 9.6.1} Suppose $S$ is a {\it given} semilattice of cyclic subsemigroups $S_\alpha\ (\alpha\in Y)$. How to decide whether it is a  {\it strong} semilattice? 
For each covering $\alpha\succ\beta$ let $S_\alpha=\langle a\rangle=C_{m,n}$ and $S_\beta=\langle b\rangle=C_{m',n'}$. Compute the unique\\
$k\in\{0,1,...,m'+n'-1\}$ that satisfies $ab=b^{k+1}$. (Depending on how $S$ is "given" that may be easier said than done.)

{\it Case 1:} $k\ge 1$. Then $Q:=(m,n,m',n';k)$ (being realizable) satisfies (R1) and (R2)=(SR2) by Lemma 10. If even (SR1) holds, then  there is a morphism $\sigma_{\alpha,\beta}$ with $ab=(a\sigma_{\alpha,\beta})b$.
If however (SR1) does not hold then the semilattice
$\{\langle a\rangle,\langle b\rangle\}$ is not strong. A fortiori the global semilattice $S$ is not strong. 

{\it Case 2:} $k=0$. Then closer inspection is required to decide the existence of $\sigma_{\alpha,\beta}$ (see 9.5.2). 

\vspace{2mm}
Even if  $\sigma_{\alpha,\beta}$ exists for all coverings $\alpha\succ\beta$, a lot
 of work remains. All problems pop up already for a 4-element frame $Y=\{\alpha,\beta,\gamma,\delta\}$  of cyclic semigroups $S_\alpha,S_\beta,S_\gamma,S_\delta$ (which we pick for the sake of notation). To recap, if $\alpha$ is the top and $\delta$ the bottom of $Y$, we have established that $\sigma_{\alpha,\beta},\sigma_{\beta,\delta},\sigma_{\alpha,\gamma},\sigma_{\gamma,\delta}$ exist. Next we need to check whether the morphisms $\sigma_{\alpha,\beta}\circ \sigma_{\beta,\delta}$ and $\sigma_{\alpha,\gamma}\circ \sigma_{\gamma,\delta}$ coincide.
If no, $S$ is not strong. If yes, let $\sigma_{\alpha,\delta}$ be this composed morphism.
Next one needs to check whether $\sigma_{\alpha,\delta}$ reflects the actual multiplication in $S$, i.e. whether $ad=(a\sigma_{\alpha,\delta})d$  (where $S_\delta=\langle d\rangle$). 

One still needs to verify whether our morphisms comply with the actual multiplication in the case of {\it incomparable} indices. For our $Y$ this means checking whether $bc=(b\sigma_{\beta,\delta})(c\sigma_{\gamma,\delta})$. If no, then $S$ is not strong. If yes, let us verify by induction (wlog going from $(i,j)$ to $(i+1,j)$) that it works for all powers of $b,c$ as well:

$$b^{i+1}c^j=bb^ic^j\stackrel{!}{=}(b\sigma_{\beta,\delta})(b^ic^j)\stackrel{ind.}{=}(b\sigma_{\beta,\delta})(b^i\sigma_{\beta,\delta})(c^j\sigma_{\gamma,\delta})
=(b^{i+1}\sigma_{\beta,\delta})(c^j\sigma_{\gamma,\delta})$$

\noindent
The second "=" holds because $b,\ b^ic^j$ belong to $S_\beta,\ S_\delta$ respectively, and $\beta>\delta$ are comparable.

\vspace{2mm}
{\bf 9.6.2} As to {\it constructing} a strong semilattice $S$ based on a fixed "frame" $Y$ and fixed disjoint cyclic semigroups
 $S_\alpha\ (\alpha\in Y)$, let us again  stick to $Y=\{\alpha,\beta,\gamma,\delta\}$.
We did most of the work already in 5.4.3 where 
 for a {\it specific} set of  sgr $S_\alpha,S_\beta,S_\gamma,S_\delta$ we argued that there are exactly 36 strong semilattices $S$. 
 
 \vspace{2mm}
 First let us adapt the notation $Exq(...)$ from 2.3. Namely, for $S_\alpha\simeq C_{m,n}$ and $S_\beta\simeq C_{m',n'}$ define the set $Exq(\alpha,\beta)$ as follows: $k\in Exq(\alpha,\beta)$ iff $k\ge 1$ and $(m,n,m',n';k)$ satisfies (SR1) and (SR2). Always $Exq(\alpha,\beta)\neq\es$ in view of 2.3.1.
The construction of $S$ for a {\it general} set of cyclic sgr $S_\alpha,S_\beta,S_\gamma,S_\delta$
 is similar to 5.4.3. Compute $Exq(\alpha,\beta,\delta):=Exq(\alpha,\beta)\cdot Exq(\beta,\delta)$ and $Exq(\alpha,\gamma,\delta):=Exq(\alpha,\gamma)\cdot Exq(\gamma,\delta)$.
If the intersection $IS(\alpha,\delta):=Exq(\alpha,\beta,\delta)\cap Exq(\alpha,\gamma,\delta)$ is empty, then there is no strong semilattice $S$.

If $IS(\alpha,\delta)\neq\es$, then there are $ss\ge 1$ strong semilattices and $ss$ can be calculated as follows. Each $k\in IS(\alpha,\delta)$ yields $ch(\beta,k)$ many choices $\{k_{\alpha,\beta},k_{\beta,\delta}\}$, and these choices bijectively match pairs
$\{\sigma_{\alpha,\beta},\sigma_{\beta,\delta}\}$ of morphisms. Similarly $ch(\gamma,k)$ is defined. If say $IS(\alpha,\delta)=\{k,k',k''\}$, then

$$ss=ch(\beta,k)ch(\gamma,k)+ch(\beta,k')ch(\gamma,k')+ch(\beta,k'')ch(\gamma,k'')$$

\noindent
In 5.4.3 we had $IS(\alpha,\delta)=\{k\}$ and $ss=ch(\beta,k)ch(\gamma,k)=6\cdot 6=36$.
For larger semilattices $Y$ and any fixed $\alpha>\delta$ in $Y$ there may be more than two chains  $\alpha\succ\beta\succ\cdots\succ\beta'\succ\delta$ 
and $\alpha\succ\gamma\succ\cdots\succ\gamma'\succ\delta$ and so forth. Accordingly

$$IS(\alpha,\delta):=Exq(\alpha,\beta,...,\beta',\delta)\cap Exq(\alpha,\gamma,...,\gamma',\delta)\cap\cdots$$

\noindent
Clever ways to calculate $ss$ remain to be found.

\vspace{2mm}
{\bf 9.6.3} Let us adapt the starter question of 9.6.1: Suppose $S$ is a {\it given} semilattice of cyclic subsemigroups $S_\alpha\ (\alpha\in Y)$. How to decide whether it is a  {\it ordinary} semilattice? This is a silly question; every semilattice of sgr is ordinary!
However the analogon of 9.6.2 is more demanding: How to {\it construct} an (ordinary) semilattice $\widetilde{S}$ based on a fixed frame $Y$ and fixed disjoint cyclic semigroups
 $\widetilde{S}_\alpha\ (\alpha\in Y)$?  

 \vspace{2mm}
 Let us start out as in 9.6.2. Thus say $\widetilde{S}_\alpha\simeq C_{m,n}$ and $\widetilde{S}_\beta\simeq C_{m',n'}$. Then by definition  $k\in \widetilde{Exq} (\alpha,\beta)$ iff either $k=0$ or ($k\ge 1$ and $(m,n,m',n';k)$ satisfies (R1) and (R2)). Evidently  $Exq(\alpha,\beta)\s \widetilde{Exq}(\alpha,\beta)$ and the latter contains $\zero$. 
Unfortunately, multiplying $k$'s from various sets $\widetilde{Exq}(...)$ and comparing the arising products will not work since  the $k$'s are no longer coupled to morphisms.

 Settling matters
for the two 3-element semilattices $Y=\{\beta,\gamma,\delta\}$ would be a  first step. Here comes an easy special case.  Let $\widetilde{S}_\beta,\widetilde{S}_\gamma,\widetilde{S}_\delta$ be arbitrary cyclic semigroups. Is there a semilattice $\widetilde{S}$ where each quotient $\alpha>\alpha'$ in $Y$ is coupled to a {\it trivial} ideal extension? The answer is yes when $Y=\{\beta\succ\gamma\succ\delta\}$ is a chain, but when
$Y=\{\beta\succ\delta\prec\gamma\}$ then the answer depends on the structure of $S_\delta$.

\vspace{2mm}
{\bf 9.6.4}
Attempts towards characterizing $Y$-frame semilattices $S$ of given cyclic semigroups $S_\alpha\ (\alpha\in Y)$ were also made in [AS]. However, the conditions in Theorem 3 of [AS] are wanting. For instance, certain functions $f:Y\times Y\to\N$ and $g:Y\times Y\to\Z$ are defined which enter the definition of the multiplication in $S$. But $f,g$ are incompletely defined in the sense that they must guarantee the satisfaction of some cumbersome\footnote{To quote the authors (page 5): {\it Condition (viii) says essentially  that associativity of third degree and fourth degree terms is sufficient to guarantee all associativity.}} condition (viii). However, no hint is given how to fine-tuning $f,g$ accordingly. 
There is no talk about strong semilattices either. Further, after a brief glimpse, I could not  muster enough energy\footnote{Readers are welcome to identify (and rewrite!) potential hidden insights.} to unravel whether or how all of this simplifies when $|Y|=2$. If all cyclic semigroups $S_\alpha$ are {\it infinite}, then things [AS, Cor.4] look  smoother. This is not  surprising since then neither indices nor periods of cyclic  semigroups are interfering.

\vspace{3mm}
{\bf 9.7}
The two ingredients in the proof of Theorem 11 may actually carry over to scenarios where a semilattice of two {\it non-cyclic} semigroups  needs to be built.
Recall, the first ingredient is setting up some\footnote{It helps if there is only {\it one} candidate operation (as in our scenario).} groupoid operation and checking its well-definedness. Second, one needs to verify that suitable generators of the groupoid associate with everybody. The case distinction in the proof of Thm.11 could have been trimmed a bit by exploiting commutativity. How to do that systematically, remains to be detected.

\section{Loose ends}

 Subsection 10.1 is purely graph-theoretic and 
proves that a Noetherian digraph is Church-Rosser iff it is locally confluent.
In 10.2 congruences on commutative semigroups are defined and four easy examples are given. The more subtle {\it Thue congruence} in 10.3 is the technical basis of the previously considered relatively free c. semigroups $RFCS(..)$. Subsection 10.4 exploits 10.1 in order to establish that the efforts in Sec.6 to achieve the local confluence of a presentation, always terminate.
In 10.5  we turn to {\it arbitrary} semigroups and render the highlights of the fundamental  Green equivalence relations $\cal H,L,R,D,J$, and the smallest semilattice congruence $\eta$. 
Subsection 10.6 demonstrates that matters simplify drastically
when commutativity (and/or finiteness) is added.

\vspace{3mm}
{\bf 10.1} Let $D$ be a digraph (=directed graph) with vertex set $V$ and arc-set $Arc\s V\times V$. As in Section 8 (identification of the Archimedean components) the {\it connected components} of $D$ by definition are the connected components of the underlying undirected graph $G$.
A directed path is a possibly infinite sequence of vertices $(x_1, x_2,x_3,...)$  such that always $(x_i,x_{i+1})\in Arc$ (and $x_i=x_j$ for $i\neq j$ is allowed). We call $x\in V$ {\it irreducible} if it has outdegree 0, i.e. there are no arcs of type $(x,y)$. Furthermore, call $D$ {\it Noetherian} if there are no infinite directed paths. In particular there are no (finite or infinite)  directed circuits. Evidently each connected component of a Noetherian digraph contains irreducible vertices.

A digraph $D$ is {\it Church-Rosser} if each connected component contains exactly one irreducible vertex. And $D$ is {\it locally confluent}\footnote{We leave it to the reader to show that "locally confluent in the 6.4 sense" implies "locally confluent in the above sense". Hint: When $w_1,w_2,w\in V$ are such that $w_1,w_2$ are subwords of $w$, then $lcm(w_1,w_2)$ is a subword of $w$ as well.} if for any arcs $a\to b$ and $a\to c$ in $D$ there is a vertex $d$ such that there are finite directed paths $b\to\cdots\to d$ and $c\to\cdots\to d$.
Let $D$ be Noetherian. Trivially, if $D$ is Church-Rosser, then it is locally confluent. Surprisingly the converse is true as well (Newman 1942):

\vspace{3mm}
{\bf Theorem 12: }{\it Let $D$ be a Noetherian digraph. Then $D$ is Church-Rosser iff it is locally confluent.}

\vspace{3mm}
{\it Proof.} In order to show the nontrivial directon, we proceed in two steps (50) and (51). First another definition. A digraph is {\it confluent} if for any directed paths  $v\to\cdots\to v_1$ and $v\to\cdots\to v_2$  there is a vertex $w$ such that there are (finite) directed paths $v_1\to\cdots\to w$ and $v_2\to\cdots\to w$.

$$(50)\quad Let\  D\ \hbox{\it be Noetherian. Then locally confluent implies confluent.} $$

To verify this, consider directed paths  $v\to\cdots\to v_1$ and   $v\to\cdots\to v_2$.
We must exhibit some vertex $w$ and directed paths $v_1\to\cdots\to w$ and $v_2\to\cdots\to w$.

Let $v\to v_1'$ and $v\to v_2'$ be the first arcs in these directed paths.
By local confluence there exists $w'$ with $v_1'\to\cdots\to w'$ and $v_2'\to\cdots\to w'$
(see Fig.10A). According to Noetherian induction (e.g. see [Co,p.61]) we can assume that confluence takes place at vertex $v_1'$. Specifically, given the directed paths $v_1'\to\cdots\to w'$ and $v_1'\to\cdots\to v_1$, there must be a vertex $u$ such that there exist directed paths $w'\to\cdots\to u$ and $v_1\to\cdots\to u$. Similarly, applying Noetherian induction to the directed paths $v_2'\to\cdots\to w'\to u$ and $v_2'\to\cdots\to v_2$ yields a vertex $w$ and directed paths $u\to\cdots\to w$ and $v_2\to\cdots \to w$. A look at Fig.10A confirms that the desired
directed paths $v_1\to\cdots\to w$ and $v_2\to\cdots\to w$ exist.

\begin{itemize}
    \item[(51)] {\it Let $y,z$ be any distinct vertices in the same connected component of $D$. Then there are directed paths $y\to\cdots\to w\ and\ z\to\cdots\to w$. (Here $w\in\{y,z\}$ is allowed.)}
\end{itemize} 

 \noindent
 In particular (51) precludes the existence of two irreducible vertices $y,z$ in the same connected component, i.e. (51) implies Church-Rosser. 

We prove (51) by induction on the distance $dist(y,z)$, i.e. length of the shortest length $n$ of an undirected path between $y$ and $z$. If $n=1$, then either $w=y$ or $w=z$ does the job. Let $n>1$. 

{\it Case 1:} The last arc in the path from $y$ to $z$ is of type $z\to z'$, see Fig.10B.
Since $dist(y,z')\le n-1$, by induction $y,z'$ have a common  bound $w$. This also is a common  bound of $y,z$.

{\it Case 2:} The last arc in the path from $y$ to $z$ is of type $z'\to z$, see Fig.10C.
Again by induction $y,z'$ have a common upper bound $w'$. Since by (45) we have confluence at $z'$, for some vertex $w$ there are directed paths $z\to\cdots\to w$ and $w'\to\cdots\to w$. Therefore $w$ is a common  bound of $y,z$. $\square$

\vspace{3mm}
\begin{center}
\includegraphics[scale=0.67]{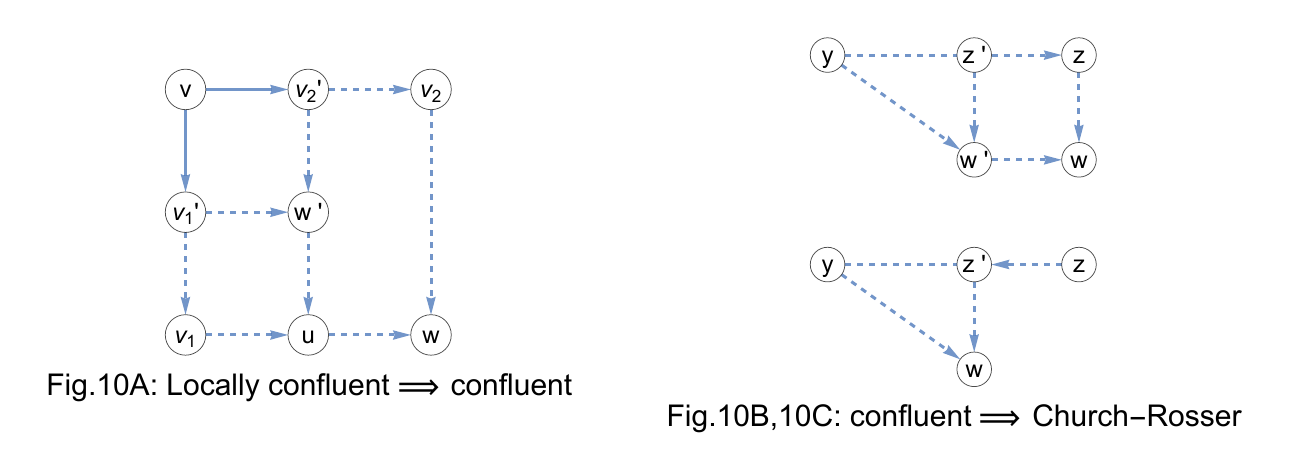}
\end{center}

\vspace{3mm}
{\bf 10.2} We emphasize that $S$ will be a {\it commutative} semigroup (albeit  generalizations are possible) in the remainder of Section 10, except for Subsection 10.5. An equivalence relation $\theta\s S\times S$ is called a {\it congruence (relation)} if for all $a,b,c\in S$ it holds that

$$(52)\quad (a,b)\in\theta\ \Ra\ (ac,bc)\in\theta$$

 For any congruence $\theta$ of $S$ define the {\it quotient} $S/\theta$ as the set of all congruence classes $a\theta:=\{b\in S:\ (a,b)\in\theta\}$. This quotient becomes a semigroup itself  by setting  $(a\theta)(b\theta):=(ab)\theta$. (Check that  this operation is well-defined and  associative.) Here come four easy kinds of congruences.

 First, if $f:S\to T$ is a morphism, then $(a,b)\in ker(f):\LRa af=bf$ defines a congruence
 $ker(f)\s S\times S$, called the {\it kernel} of $f$. Conversely, let $\theta$ be any congruence of $S$. Then a morphism (even epimorphism) $g:S\to S/\theta$ is obtained by setting $ag:=a\theta$.

 Second, because of (39) the Arch. components $A_e\ (e\in Y)$ of each c.f. sgr $S$ are the $\eta$-classes of a congruence $\eta$ of $S$ that satisfies $S/\eta\simeq Y$.

 Third, let $I\s S$ be an ideal and let $\theta$ be the equivalence relation whose $\theta$-classes are $I$ and all singletons $\{x\}\ (x\not\in I)$. One verifies that $\theta $ is a congruence.  The associativity issue of the Rees quotient $S/I$ in 2.9 now vaporizes since $S/I$ is reckognized as an instance of a quotient sgr $S/\theta$. 
 
 Fourth, the set $Con(S)$ of all congruences of $S$ is easily seen to be a closure system and hence is a lattice (see 5.5.). This implies that for each set $X\s S\times S$ there is a {\it smallest} congruence $\theta(X)$ that contains $X$. Namely $\theta(X)$, called the congruence {\it generated} by $X$, is the intersection of all congruences containing $X$.
 Always $\bigtriangleup,\bigtriangledown\in Con(S)$, where $\bigtriangleup:=S\times S$ is the largest, and  $\bigtriangledown:=\{(x,x):\ x\in S\}$ is the smallest element of the lattice $Con(S)$.

 \vspace{3mm}
 {\bf 10.3} Let $S=\langle a,b\rangle$ be  any commutative (possibly infinite) semigroup that satisfies $a^3=a$ and $ab^2=ab$. Viewing $F_2$ as the semigroup of all words
 $\alpha^i\beta^j$ over the alphabet $\{\alpha,\beta\}$, the map $f:\ F_2\to S$ defined by $(\alpha^i\beta^j)f:=a^ib^j$ is an epimorphism:

 $$(\alpha^i\beta^j\cdot \alpha^s\beta^t)f=(\alpha^{i+s}\beta^{j+t})f=a^{i+s}b^{j+t}=
 a^ib^j\cdot a^sb^t=(\alpha^i\beta^j f)\cdot (\alpha^s\beta^t f)$$

\noindent
By the {\it First Isomorphism Theorem} $S\simeq F_2/\theta$, where $\theta:=ker(f)$. From $(\alpha\beta^2)f=(\alpha\beta)f$ follows that $(\alpha\beta^2,\alpha\beta)\in\theta$. Similarly $(\alpha^3,\alpha)\in\theta$. This indicates how to formalize the postulated "largest" semigroup $RFCS(a,b:\ a^3=a,\ ab^2=ab)$ of 6.2.1. Namely, the {\it Thue congruence} induced by our presentation is defined as the congruence $\theta_0$ 
generated by $\{(\alpha\beta^2,\alpha\beta),(\alpha^3,\alpha)\}$. Evidently $\theta_0\s\theta$. Hence the {\it Second Isomorphism Theorem} implies that $F_2/\theta$ (which is $\simeq S$) is an epimorphic image of $F_2/\theta_0$, which we hence take as the formal definition of $RFCS(a,b:\ a^3=a,\ ab^2=ab)$. 

One can also obtain epimorphic images of $RF:=RFCS(a,b:\ a^3=a,\ ab^2=ab)$  by adding further relations. For instance $RF':=RFCS(a,b:\ a^3=a,\ ab^2=ab, b^3=a^2b)$ is an epimorphic image of $RF$.
All of this has nothing to do with local confluence, nor with the finiteness of $RFCS(...)$. In fact one can show (try) that $RF$ contains the infinite set $\{b,b^2,b^3,..\}$, whereas $|RF'|=6$.

\vspace{3mm}
{\bf 10.4} But now we {\it do} turn to local confluence, and for this it suits us to reconsider $RF_2=RFCS(a,b:\ b^4=b^2,\ a^3=b^2,\ a^4=a)$ from 6.3. Recall from 10.3 that formally $RF_2:=F_2/\theta_0$, where $\theta_0=\theta(X)$ is the congruence generated by $X:=\{(\beta^4,\beta^2),(\alpha^3,\beta^2),(\alpha^4,\alpha)\}$. Hence the elements of $RF_2$ are $\theta_0$-classes.

In addition to 10.3 we  view these as the connected components of a certain digraph $D(X)$. 
Recall from 9.1 the concept of adjoining an identity.
Thus, if $F_2=\langle \alpha,\beta\rangle$, then $F_2^1$ additionally contains the "empty word" $\one$; say $\alpha^7\beta^9\one=\alpha^7\beta^9$. By definition the digraph $D(X)$ has vertex set $V:=F_2$ and
if $w,w'\in V$, then by definition there is an arc between these vertices (i.e. $w\to w'$) iff e.g. $w=\alpha^3v$ and $w'=\beta^2v$ for some $v\in F_2^1$ (instead of $(\alpha^3,\beta^2)$ any other ordered pair in $X$ can be taken). 

\vspace{2mm}
{\bf 10.4.1} Looking closer at $D(X)$, because $\theta_0$ is symmetric and transitive, each $\theta_0$-class is a union of connected components of $D(X)$. It turns out\footnote{For this it isn't enough  that $\theta_0$ is the smallest congruence containing $X$; one needs a deeper understanding of $\theta_0$, as provided in any book about Universal Algebra.} that actually each $\theta_0$-class is {\it one}  connected component.

If $w\to w'$ in $D(X)$ then (why?) $w>_M w'$. It follows, viewing that  $(F_2,\ge_M)$ is a Noetherian poset by 2.8.1, that $D(X)$ is a Noetherian digraph (and this is independent of the particular presentation $X$ of $RF_2$). Therefore, starting a directed path at any vertex $w\in V$, and extending it in arbitrary fashion as far as possible, one will end up with an irreducible vertex $v$ after finitely many steps. Trouble is, $v$ needs not be unique. Indeed, recall from 6.3 that for $w:=\alpha^4$ one may end up in $v=\alpha$ or in $v'=\alpha\beta^2$.

Fortunately, by Theorem 12, {\it if} such a faulty presentation $X$ can be replaced by a locally confluent $X'$, {\it then} the new\footnote{As a digraph $D(X')$ may be much different from $D(X)$, yet it remains Noetherian and of course still has the $\theta_0$-classes as connected components.} digraph $D(X')$ is Church-Rosser, and so the unique irredundant vertices in the connected components of $D(X')$ can serve as the normal forms for the elements of $RFCS(..)$. In 10.4.2 the digraph  $D(X)$ more generally has vertex set $V=F_k$ and we show that the representation $X$ {\it can} indeed be replaced by a locally confluent $X'$. A crucial ingredient will be "Dickson's Lemma" which states:

\begin{itemize}
    \item[(53)] {\it The poset $(F_k,\le_c)$ (see 2.8.1) has only finite antichains.}
\end{itemize}

\noindent
Here an {\it antichain} is a set of pairwise incomparable elements. A proof of (53) can be found in [Wi,p.186] or [RG,p.48].

\vspace{2mm}
{\bf 10.4.2} It is a priori plausible that the process  of adding new relations to a given presentation $X$ (in order to make it locally confluent) never stabilizes. Suppose $X$ consists of the relations $\rho_i:v_i\to w_i\ (1\le i\le n)$. Say $\rho_1,\rho_4$ is the "left most" pair that violates local confluence wrt $\{\rho_1,...,\rho_n\}$. We know from Section 6 how to find a (derivable) relation 
$\rho_{n+1}:v_{n+1}\to w_{n+1}$ that establishes local confluence. By construction the vertex $v_{n+1}$ has outdegree 0 in $D(X)$, and so $v_i\not\le_c v_{n+1}$ for all $1\le i\le n$. 
Suppose now all $\rho_i,\rho_j\ (1\le i<j\le n)$ are locally confluent
but $\rho_7,\rho_{n+1}$ is not. Then add a suitable relation $\rho_{n+2}:v_{n+2}\to w_{n+2}$ to fix that. Arguing as above it holds that $v_i\not\le_c v_{n+2}$ for all $1\le i\le n+1$.

By way of contradiction suppose that continuing in this fashion we {\it never} establish local confluence for all pairs of relations. Then the sequence\\ $v_1,...,v_n,v_{n+1},...$ is infinite and such that 

$$(54)\quad v_i\not\le_c v_j\ for\ all\ n\le i<j$$

\noindent
By (53) it suffices to exhibit and infinite antichain $\{v_{i_1},v_{i_2},...\}$ in $(F_k,\le_c)$. We put $v_{i_1}:=v_n$ and by induction assume that $\{v_{i_1},...,v_{i_t}\}$ is an antichain with $t\ge 1$ and $i_1<i_2<\cdots<i_t$. By (54) it suffices to pinpoint an index $i_{t+1}>i_t$ such that $v_{i_{t+1}}\not\le_c v_{i_1},...,v_{i_t}$.
Writing (say) $v_{i_1}\!\downarrow\ :=
\{v\in F_k:\ v\le_c v_{i_1}\}$ it is clear\footnote{To spell it out, if say $v_{i_1}=a_1^{\alpha_1}\cdots a_k^{\alpha_k}\in F_k$, then $|v_{i_1}\!\downarrow|=(\alpha_1+1)(\alpha_2+1)\cdots(\alpha_k+1)$.} that $v_{i_1}\!\downarrow\cup\cdots\cup v_{i_t}\!\downarrow$ is finite, and so there are (infinitely many) indices $j>i_t$ with $v_j\not\in v_{i_1}\!\downarrow\cup\cdots\cup v_{i_t}\!\downarrow$. Take any such $j$ and put $i_{t+1}:=j$. 

(Let us mention that an argument along the lines of 10.4.2 would have made the proof of [BL,Lemma 2] more intelligible.)

\vspace{3mm}
{\bf 10.4.3} Let $K$ be any field. The whole "business" of local confluence can be raised from commutative semigroups to the level of polynomial rings $R:=K[x_1,...,x_k]$. Then one e.g. can decide the following: Given a finite basis of some ideal $I\s R$, when do  elements $f+I$ and $g+I$ of the factor ring $R/I$ represent the same element? (Equivalently: Is $f-g\in I$?) Trouble is, the handling of critical pairs of {\it poly}nomials $f,g\in R$ gets more complicated than the straightforward (Section 6) handling of critical pairs of {\it mono}mials $a_1^{\alpha_1}\cdots a_k^{\alpha_k}$ and $a_1^{\beta_1}\cdots a_k^{\beta_k}$. The crucial insight is in Buchberger's PhD thesis of 1965. As a gentle introduction to these matters  (key word: Gr\"obner bases) we recommend [Wi].

Historically however the arrow does not simply go from semigroups to polynomial rings. A few remarks must suffice; [Bu] provides a broader picture. Critical pairs were introduced in a crucial  1969 article of Knuth-Bendix in the context of ordinary term rewriting systems. (The results of Newman 1942 (=Theorem 12) and Dickson 1913 (see (53)) are older still.)  Exploiting critical pairs for commutative semigroups can\footnote{This opinion is e.g. supported in [Bu,p.20].} be attributed to [BL]. This notwithstanding  Lankford and Ballantyne  acknowledging that similar ideas (on the level of rings) were, unbeknownst to them, used in [Be] three years earlier. Bergman briefly mentions in [Be,Sec.9.1] that his results carry over to commutative semigroups. The details however are spelled out only in [BL] (and partly improved in our own Sections 6 and 10). Bergman in turn seems to have obtained his results unbeknownst of Buchberger\footnote{Quoting from [Bu, p.20]:
{\it Apparently independently of my own work, Bergman (1978) rediscovered essentially the same algorithm, however, in a slightly more general form, namely...} But later on [Bu,p.20]: {\it However the approach is not broad enough to encompass the case of integer polynomial ideals because...}}, whom he does not cite.

\vspace{3mm}
{\bf 10.5 } Let us  glance\footnote{All omitted proofs in 10.5, and much more about the Green relations (e.g. their interplay with "regular" elements), can be found in reader-friendly form in [Go].} at the five {\it Green equivalence relations}, as well as  $\eta$, in {\it arbitrary} semigroups $S$. Afterwards (in 10.6) we get stunned how things collapse in the commutative case.

\vspace{3mm}
{\bf 10.5.1}
As opposed to 2.9, in arbitrary semigroups one has to distinguish between {\it left-ideals,  right-ideals}, and {\it (2-sided) ideals} (the definitions being obvious). 
Recall from 9.1 that $S^\one$ is the semigroup obtained from $S$ by adjoining an identity.
Thus one says that $x,y\in S$ are {\it $\cal J$-related} iff they generate the same ideal, i.e.  $S^1xS^1=S^1yS^1$. They are {\it $\cal L$-related} iff they generate the same left-ideal, i.e.  $S^1x=S^1y$. They {\it $\cal R$-related} iff they generate the same right-ideal, i.e.  $xS^1=yS^1$. Obviously ${\cal L,R,J}$ are equivalence relations.
The difference between ${\cal L}$ and ${\cal R}$ can be drastic.
 Consider say a {\it left-zero} sgr $S$ where by definition  $xy=x$ for all $x,y\in S$. Then ${\cal L}=\bigtriangleup$ but ${\cal R}=\bigtriangledown$. Nevertheless,  Green discovered that always ${\cal L}\circ {\cal R}={\cal R}\circ {\cal L}=:{\cal D}$. Evidently ${\cal L,R}\s {\cal D}$. Furthermore 

 \begin{itemize}
     \item[] $x({\cal L}\circ {\cal R})y\Ra (\exists z)(x{\cal L}z{\cal R}y)\Ra S^1x=S^1z,\
 zS^1=yS^1$
 \item[] $\Ra S^1xS^1=S^1zS^1=S^1yS^1\Ra x{\cal J}y,$
 \end{itemize}
 
\noindent
and so ${\cal D}\s {\cal J}$. For $|S|<\infty$ it holds that ${\cal D}={\cal J}$. One can show that all ${\cal L}$-classes contained in a ${\cal D}$-class have the same cardinality, and likewise for the ${\cal R}$-classes. If ${\cal H}:={\cal L}\cap {\cal R}$, then all ${\cal H}$-classes contained in a ${\cal D}$-class have the same cardinality. If a ${\cal D}$-class $D$ contains an idempotent, then each ${\cal L}$-class and each ${\cal R}$-class contained in $D$ has at least one idempotent. Furthermore the ${\cal H}$-classes in $D$ that happen to contain an idempotent are mutually isomorphic subgroups of $S$.

\vspace{3mm}
{\bf 10.5.2} A semigroup $S$ is $\cal J$-{\it trivial} if ${\cal J}=\bigtriangledown$. In stark contrast (but unfortunately with similar name) one says $S$ is $\cal J$-{\it simple} if ${\cal J}=\bigtriangleup$. Thus, in the latter case, the only ideal of $S$ is $S$ itself.

The set $S/{\cal J}$ of all ${\cal J}$-classes $[x]$ becomes partially ordered by putting \\$[x]\le_{\cal J} [y]$ iff $S^1xS^1\s S^1yS^1$. The smallest element of the poset $S/{\cal J}$ is the kernel $K(S)$. If $\cal J$ happens to be a congruence then the semigroup $S/{\cal J}$ is $\cal J$-trivial.

\vspace{3mm}
{\bf 10.5.3} 
One calls $\theta\in Con(S)$  a {\it semilattice congruence} if $S/\theta$ is a semilattice. One checks that $\theta$ is a semilattice congruence iff $(ab)\theta (ba)$ and $a\theta a^2$ for all $a,b\in S$. It follows that the intersection $\eta$ of all semilattice congruences is itself a semilattice congruence, and evidently the smallest one.

In order to show that ${\cal J}\s\eta$, let $\le$ be the partial ordering of the semilattice $S/\eta$. It follows from $x{\cal J}y$ that $x=syt$ and $y=s'xt'$. Hence
$x\eta=(s\eta)(y\eta)(t\eta)$ and $y\eta=(s'\eta)(x\eta)(t'\eta)$, hence $x\eta\le y\eta$ and $y\eta\le x\eta$, hence $x\eta=y\eta$, hence $x\eta y$.
 To summarize 
 
 $${\cal H}\s {\cal L}\s {\cal D}\s {\cal J}\s \eta\hspace{1cm} (one\ can\ replace\ {\cal L}\ by\ {\cal R}).$$    
 
\vspace{3mm}
{\bf 10.6}
In the remainder of the article $S$ is again {\it commutative}. Then (why?)\\
${\cal H}={\cal L}={\cal R}={\cal D}={\cal J}$!
As a perk, in contrast to the general case ${\cal J}$ is a congruence:

$$a{\cal J}b\Ra (\exists t,s\in S^1)(a=bt,\ b=as)\Ra (ac=bct,\ bc=acs)\Ra ac{\cal J}bc.$$

Observe that $a{\cal J}b\ \LRa\ (a\le_{\cal J} b\ and\ b\le_{\cal J} a)$, where $\le_{\cal J}$ is the preorder defined in (9). This preorder is a partial order on $S$ iff $S$ is $\cal J$-trivial. Recall that e.g. semilattices and all semigroups $F_k$ are $\cal J$-trivial. If $S$ is not $\cal J$-trivial, it can be condensed to the $\cal J$-trivial semigroup $S/{\cal J}$ (viewing that ${\cal J}$ is a congruence).

\vspace{2mm}
{\bf 10.6.1}
If additionally $S$ is {\it finite}, it gets better still. Recall from Theorem 4 that all c.f. nilsemigroups are $\cal J$-trivial (i.e. partially ordered by $\le_{\cal J}$). Generally it holds that each kernel $K(A_e)$ is a $\cal J$-class ($e\in E(S)$), and all other $\cal J$-classes are singletons.

Two immediate consequences. First, a c.f. sgr is $\cal J$-trivial iff it is a semilattice of nilsemigroups. 
Second, a c.f. sgr $S$ satisfies ${\cal J}=\eta$ iff $S$ is a semilattice of Abelian groups. For instance (see 8.2), the sgr $(\Z_n,\odot)$ satisfies
${\cal J}=\eta$ iff $n$ is square-free. As a special case, a c.f. sgr is $\cal J$-simple iff it is an Abelian group. To further specialize, a c.f. sgr $S$ is {\it congruence-simple} (i.e. $Con(S)=\{\bigtriangleup,\bigtriangledown\}$) iff\footnote{More generally: The (not necessarily commutative) finite sgr $S$ is congruence-simple iff it is a congruence-simple group, i.e. one without proper normal subgroups. The classification of the latter class of groups is still being finalized and constitutes the greatest mathematical collaboration ever.} $S\simeq C_p$ for some prime $p$.

Recall from Section 8 that in the c.f. scenario all $\eta$-classes are  Archimedean subsemigroups $A$, i.e. having unique\footnote{In the non-commutative case an $\eta$-class can have several idempotents.} idempotents. The  ${\cal J}$-class  within  $A$ that catches the idempotent is the kernel $K(A)$. (Recall from 10.5.1 that even in the non-commutative case each $\cal H$-class with an idempotent is a group.)

\vspace{3mm}
{\bf 10.6.2} Observe that $\eta$ is a {\it retract} congruence, i.e. there is a set $Y$ of {\it representatives} of the $\eta$-classes such that $Y$ is a subsemigroup of $S$. Indeed, take 
 $Y:=E(S)$. 
 
Consider now the four Archimedean components of $\Z_{18}$ listed in 8.5. The two 6-element components  are the $\cal J$-classes [1] and [10] respectively. As to the component $\{3,9,15\}$, it splits into the $\cal J$-classes $[3]=\{3,15\}$ and $[9]=\{9\}$. As to the component $\{0,6,12\}$, it splits into the $\cal J$-classes 
$[12]=\{6,12\}$ and $[0]=\{0\}$. The structure of the poset $\Z_{18}/{\cal J}$ is rendered in Fig.1B. The set $\{0,1,3,9,10,12\}$ of representatives happens to be a  ssgr $Y'$ of $\Z_{18}$
 as well, and so $\cal J$ is a retract congruence. (As for any retract congruence, $Y'\simeq \Z_{18}/{\cal J}$.)

\vspace{3mm}
{\bf Open Question 3: For which finite commutative\\ semigroups $S$ is ${\cal J}$ a retract congruence? }


\section{ References}
\begin{itemize}
    \item[\bf A] R.B.J.T. Allenby, Rings, fields and groups, Butterworth-Heinemann 2001.
    \item[\bf Ar] M.A. Armstrong, Groups and symmetry, Springer-Verlag 198S.
     \item[\bf AS] B.D. Arendt, C.J. Stuth, On the structure of commutative periodic semigroups, Pacific J. of Math 35 (1970) 1-6.
 \item[\bf Be] G.M. Bergman, The diamond lemma for ring theory,
     Advances in Mathematics  29 (1978) 178-218.
     \item[\bf Bu] B. Buchberger, History and basic features of the critical-pair/completion procedure, J. Symbolic Computation 3 (1987) 3-38.
     \item[\bf B] P. Bundschuh, Einf\"uhrung in die Zahlentheorie, Springer-Verlag 1988.
    \item[\bf BC] B. Baumslag, B. Chandler, Theory and problems of Group Theory, Schaum's outline series, McGraw-Hill 1968.
    \item[\bf BF] G. Bini, F. Flamini, Finite commutative rings and their applications, Springer 2002.
    \item[\bf BL] A.M. Ballantyne, DS. Lankford, New decision algorithms for finitely presented commutative semigroups, Comp. Math. with Applications 7 (1981) 159-165.
    \item[\bf CL] A.H. Clifford, G.B. Preston, The algebraic theory of semigroups, AMS 1961.
    \item[\bf CLM] N. Caspard, B. Leclerc, B. Monjardet, Finite ordered sets, Cambridge University Press 2012.
    \item[\bf Co] P. Cohn, Basic Algebra, Springer 2005.
   
    \item[\bf FP] V. Froidure, JE. Pin, Algorithms for computing finite semigroups, Foundations of Computational Mathematics, 1997, Rio de Janeiro, Brazil. pp.112-126.
    \item[\bf Go] V. Gould, Semigroup Theory, a lecture course. Freely available on the internet.
    \item[\bf G] P.A. Grillet, Commutative Semigroups, Kluwer Academic Publishers 2001.
    \item[\bf Gr] G. Gr\"atzer, Lattice Theory: Foundation, Birkh\"auser 2011.
    \item[\bf H] J. Howie, An introduction to semigroup theory, Academic Press, London 1976.
    \item[\bf MKS] W. Magnus, A. Karrass, D. Solitar, Combinatorial Group Theory, Dover 1976.
    \item[\bf RG] J.C. Rosales, P.A Garc\'ia S\'anchez, Finitely generated commutative monoids, Nova Science Publishers Inc. 1999.
     \item[\bf S] A.V. Sutherland, Structure computation and discrete logarithms in finite Abelian $p$-groups, Mathematics of Computation 80 (2011) 477-500.
    \item[\bf W1] M. Wild, Implicational bases for finite closure systems,Technische Universit\"at Clausthal, Informatik-Bericht 89/3. (To get a pdf, search on Google Scholar: Marcel Wild, Implicational bases for finite closure systems)
    \item[\bf W2] M. Wild, The groups of order sixteen made easy, Amer. Math. Monthly 112 (2005) 20-31.
    \item[\bf W3] M. Wild, The joy of implications aka pure Horn formulas: Mainly a survey, Theoretical Computer Science 658 (2017) 264-292.
     \item[\bf Wi] F. Winkler, Polynomial algorithms in computer algebra, Springer 1996.
\end{itemize}

\end{document}